%% file: main.tex
\documentclass[twocolumn]{svjour3nojournal}

\smartqed  
\usepackage{microtype}
\usepackage{graphicx}
\graphicspath{ {./figures/} {./} }
\usepackage[labelfont=bf,font=small,labelsep=space]{caption}

\usepackage{amsmath}
\usepackage{amssymb}
\usepackage[normalem]{ulem}
\usepackage{mathtools}
\usepackage[dvipsnames]{xcolor}
\usepackage{hyperref}
\usepackage{fancyhdr}
\usepackage{pgfplots}
\usepackage{standalone}
\usepackage{subfiles,cases,empheq}
\usepackage{caption}
\usepackage{subcaption}
\usepackage{commands,MnSymbol}
\usepackage{svg}
\usepackage{picture}
\usepackage[]{mdframed}

\usepackage{changes}
\definechangesauthor[name={Olga Mula}, color=blue]{OM}
\usetikzlibrary{decorations.pathreplacing,calligraphy}
\usetikzlibrary{
        pgfplots.fillbetween
    }

\pgfplotsset{compat=1.16}

\usepackage{tikz}

\patchcmd\subequations
 {\theparentequation\alph{equation}}
 {\subequationsformat}
 {}{}

\newcommand{\subequationsformat}{\theparentequation.\arabic{equation}}

\usepackage[ruled,vlined]{algorithm2e}

\makeatletter
\newcommand{\customlabel}[2]{%
   \protected@write \@auxout {}{\string \newlabel {#1}{{#2}{\thepage}{#2}{#1}{}} }%
   \hypertarget{#1}{#2}
}
\makeatother

\makeatletter
\newcommand{\customlabelnoprint}[2]{%
   \protected@write \@auxout {}{\string \newlabel {#1}{{#2}{\thepage}{#2}{#1}{}} }%
   \hypertarget{#1}{}
}
\makeatother

\newcommand\numberthis{\addtocounter{equation}{1}\tag{\theequation}}

\journalname{Journal of Mathematical Imaging and Vision}

\begin{document}
\setlength\emergencystretch{3em}
\input{Review_JMIV/TitleAndAbstract}
\input{Review_JMIV/1Introduction_OlgaNicky}

\input{Review_JMIV/2_1LieGroups}
\input{Review_JMIV/2_2ConnectedComponents}
\input{Review_JMIV/3MorphologicalDilations}
\input{Review_JMIV/4_1ConnectedComponentAlgorithm}
\input{Review_JMIV/4_2ChoiceOfDelta}
\input{Review_JMIV/4_3AffinityMatrices}
\input{Review_JMIV/5Experiments}
\input{Review_JMIV/6ConclusionFutureWork}

\section*{Acknowledgments}

We gratefully acknowledge Gijs Bellaard for the code underlying Figure~\ref{fig:20230822:VisualizationDistanceBallsSE2:Exact} and \ref{fig:20230822:VisualizationDistanceBallsSE2:Approximate}, and the basis for the code for Fig.~\ref{fig:20230822:VisualizationDistanceBallsSO3}.

\input{Appendices}

\input{Review_JMIV/AppendixSAM}

\input{Review_JMIV/Declarations}

\bibliographystyle{siam}
\bibliography{references}
\clearpage
\end{document}


\begin{tikzpicture}
    \filldraw[blue!100, very thick, opacity = 0.2, rotate around={45:(5,5)}] (-2.07,4.5) rectangle (14.14,5.5);
    \draw[black!100, very thick](10,10)node {\Huge$[g_1]$};
    \filldraw[red!100, very thick, opacity = 0.2] (6.5,4.5) rectangle (9,5.5);
    \draw[black!100, very thick](7.75,5)node {\Huge$[g_3]$};
    \filldraw[green!100, very thick, opacity = 0.2] (9,7) rectangle (11.5,8);
    \draw[black!100, very thick](10.25,7.5)node {\Huge$[g_2]$};
\end{tikzpicture}

%% file: Review_JMIV/TitleAndAbstract.tex
\title{Connected Components on Lie Groups and Applications to Multi-Orientation Image Analysis
}


\author{Nicky~J.~van den Berg \and Olga~Mula \and Leanne Vis \and \mbox{Remco Duits}
}

\authorrunning{N.J. van den Berg et al.} 

\institute{Nicky J. van den Berg \and Leanne Vis \and Remco Duits \at
              Department of Mathematics and Computer Science, Eindhoven University of Technology, The Netherlands \\
              \email{n.j.v.d.berg@tue.nl, l.vis@tue.nl, 
              r.duits@tue.nl}      
           \and
           Olga Mula \at
           Faculty of Mathematics, University of Vienna, Vienna, Austria\\ \email{olga.mula.hernandez@univie.ac.at}
}

\date{\today}

\maketitle

\begin{abstract}
We develop and analyze a new algorithm to find the connected components of a compact set $I$ from a Lie group $G$ endowed with a left-invariant Riemannian distance. For a given $\delta>0$, the algorithm finds the largest cover of $I$ such that all sets in the cover are separated by at least distance $\delta$. We call the sets in the cover the $\delta$-connected components of I (closely related to $\check{\text{C}}$ech complexes of radius $\delta/2$). The grouping relies on an iterative procedure involving morphological dilations with Hamilton-Jacobi-Bellman kernels on $G$ and notions of $\delta$-thickened sets. We prove that the algorithm converges in finitely many iteration steps. We find the optimal value for $\delta$ using persistence diagrams. We also propose to use specific affinity matrices that allow for grouping of $\delta$-connected components based on their local proximity and alignment. 

Among the many different applications of the algorithm, 
in this article, we focus on illustrating that the method can efficiently identify (possibly overlapping) branches in complex vascular trees on retinal images. This is done by applying an orientation score transform to the images that allows us to view them as functions from $\mathbb{L}_2(G)$ where $G=SE(2)$, the Lie group of roto-translations. By applying our algorithm in this Lie group, we illustrate that we obtain $\delta$-connected components that differentiate between crossing structures and that group well-aligned, nearby structures. This contrasts standard connected component algorithms in $\mathbb{R}^2$.

\keywords{
Connected Components \and Lie Groups \and Vessel Tree Identification \and Image Analysis 
\and Morphological Operators
\and Medical Image Analysis 
}

\end{abstract} 

%% file: Review_JMIV/1Introduction_OlgaNicky.tex
\section{Introduction}\label{sec:Introduction}
\subsection{Context and Goals}
\begin{figure*}[ht!] \label{fig:intro}
    \centering
    \begin{subfigure}[t]{0.3\textwidth}
        \includegraphics[width=\textwidth]{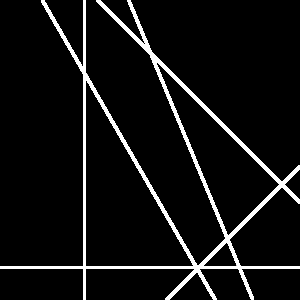}
        \caption{Reference Image.}
        \label{fig:20230922:comparisonR2vsSE2:ReferenceImage}
    \end{subfigure}\hfill
    \begin{subfigure}[t]{0.3\textwidth}
        \includegraphics[width=\textwidth]{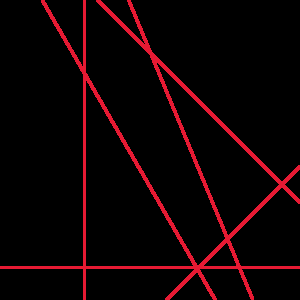}
        \caption{Connected components in $\mathbb{R}^2$.}
        \label{fig:20230922:comparisonR2vsSE2:ClassicalR2}
    \end{subfigure}\hfill
    \begin{subfigure}[t]{0.3\textwidth}
        \includegraphics[width=\textwidth]{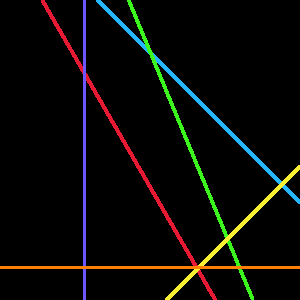}
        \caption{Connected components in \SE{2}.
        }
        \label{fig:20230922:comparisonR2vsSE2:deltaCCinSE2}
    \end{subfigure}
    \begin{subfigure}[t]{0.45\textwidth}
        \begin{picture}(160,175)
            \includegraphics[width=\textwidth]{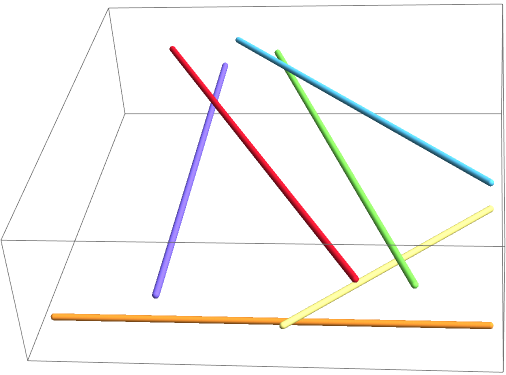}
             \put(-225,35){\small\textcolor{black}{$\theta$}}
        \end{picture}
        \caption{Connected components in $\SE{2}$ of continuous set.}
        \label{fig:20230922:comparisonR2vsSE2lines:ContinuousSet}
    \end{subfigure}\hfill
    \begin{subfigure}[t]{0.45\textwidth}
        \begin{picture}(160,175)
            \includegraphics[width=\textwidth]{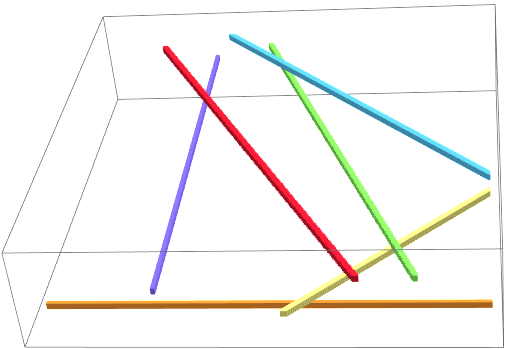}
             \put(-225,25){\small\textcolor{black}{$\theta$}}
        \end{picture}
        \caption{Connected components in $\SE{2}$ of discretized set using 16 orientations in a $\pi$-periodic framework.
        }
        \label{fig:20230922:comparisonR2vsSE2lines:DiscreteSet}
    \end{subfigure}
    
    \caption{Visualization of the connected component algorithm in $\mathbb{R}^2$ and in the Lie group \SE{2} on an image of straight lines. The classical connected component algorithm in $\mathbb{R}^2$ cannot differentiate between the different line structures, unlike the extension to \SE{2}. Here we applied the algorithm in Sec.~\ref{sec:20230822:Algorithm} using parameters $(w_1,w_2,w_3)=(0.1,1,4)$  for the left-invariant metric \eqref{eq:MTFweights}. }
    \label{fig:20230922:comparisonR2vsSE2:Lines}
\end{figure*}

\begin{figure*}[ht!] \label{fig:introChip}
    \centering
    \begin{subfigure}[t]{0.3\textwidth}
        \includegraphics[width=\textwidth]{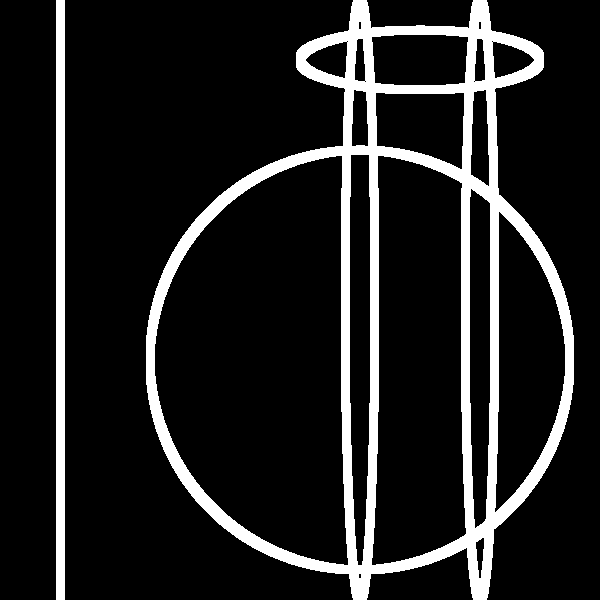}
        \caption{Reference Image.}
        \label{fig:20230922:comparisonR2vsSE2chip:ReferenceImage}
    \end{subfigure}\hfill
    \begin{subfigure}[t]{0.3\textwidth}
        \includegraphics[width=\textwidth]{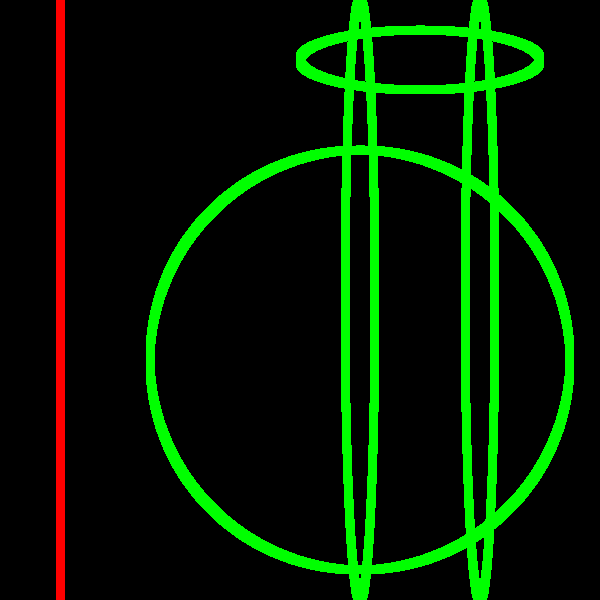}
        \caption{Connected components in $\mathbb{R}^2$.}
        \label{fig:20230922:comparisonR2vsSE2chip:ClassicalR2}
    \end{subfigure}\hfill
    \begin{subfigure}[t]{0.3\textwidth}
        \includegraphics[width=\textwidth]{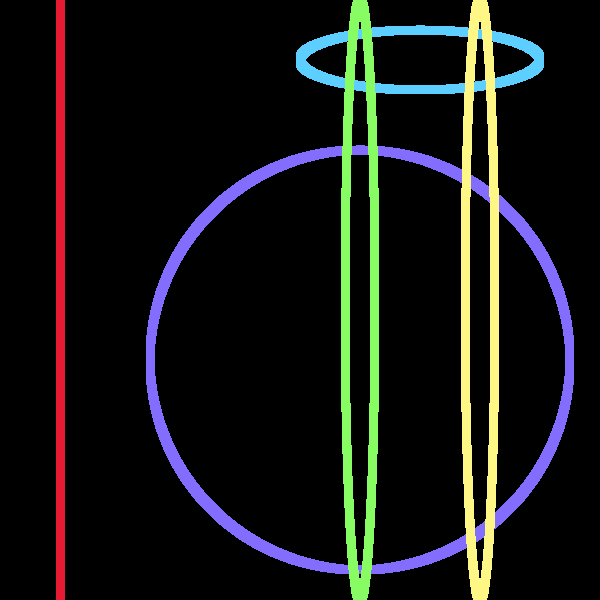}
        \caption{Connected components in \SE{2}.
        }
        \label{fig:20230922:comparisonR2vsSE2chip:deltaCCinSE2}
    \end{subfigure}
    \begin{subfigure}[t]{0.45\textwidth}
        \begin{picture}(160,175)
            \includegraphics[width=\textwidth]{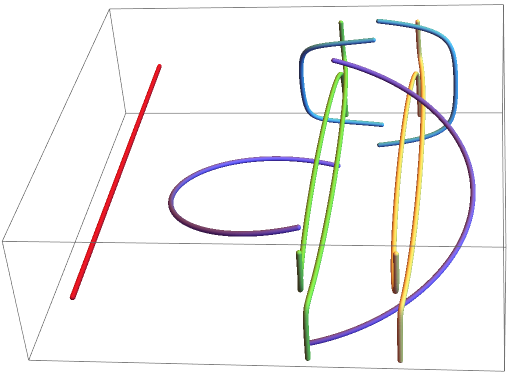}
             \put(-225,35){\small\textcolor{black}{$\theta$}}
        \end{picture}
        \caption{Connected components in $\SE{2}$ of continuous set.}
        \label{fig:20230922:comparisonR2vsSE2chip:ContinuousSet}
    \end{subfigure}\hfill
    \begin{subfigure}[t]{0.45\textwidth}
        \begin{picture}(160,175)
            \includegraphics[width=\textwidth]{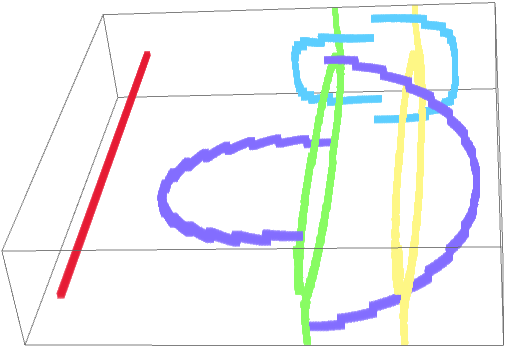}
             \put(-225,25){\small\textcolor{black}{$\theta$}}
        \end{picture}
        \caption{Connected components in $\SE{2}$ of discretized set using 16 orientations in a $\pi$-periodic framework.
        }
        \label{fig:20230922:comparisonR2vsSE2chip:DiscreteSet}
    \end{subfigure}
    \caption{Visualization of the effect of performing the connected component algorithm in $\mathbb{R}^2$ and in the Lie group \SE{2} on an image of ovals and lines. The classical connected component algorithm in $\mathbb{R}^2$ is not able to differentiate between the different line structures, unlike the extension to \SE{2}.  Here we applied the algorithm in Sec.~\ref{sec:20230822:Algorithm} using parameters $(w_1,w_2,w_3)=(0.2,1,4)$ for the left-invariant metric \eqref{eq:MTFweights}.}
    \label{fig:20230922:comparisonR2vsSE2:chip}
\end{figure*}

In image analysis, the task of connecting components plays an important role in numerous applications such as image segmentation \cite{Chudasama2015ImageSegmentation,Meyer1990Morphological,DeepthiMurthy2014Brain,Wu2017Improved}, object recognition \cite{Shih1991Object,brendel1995knowledge,Zhao2005Medical,foresti1999object,Gu1998Semiautomatic}, motion tracking \cite{foresti1999object,Gu1998Semiautomatic,Murray1994Motion,Figueroa2003Flexible,Pawar2017Morphology}, and data compression \cite{Salembier1996Video,Won1998BlockBased}. 
In this article, we develop a novel method for identifying connected components in images with line structures that possibly intersect with each other, as in Fig.~\ref{fig:20230922:comparisonR2vsSE2:Lines},~\ref{fig:20230922:comparisonR2vsSE2:chip}. When applying a naive connected component algorithm in $\mathbb{R}^2$ on these images, the algorithm finds only one or two connected components, as it cannot differentiate between the different structures at crossings (cf. Fig.~\ref{fig:20230922:comparisonR2vsSE2:ClassicalR2},~\ref{fig:20230922:comparisonR2vsSE2chip:ClassicalR2}).

We present a strategy that does not automatically merge different crossing structures. The image data is lifted from $\mathbb{R}^2$ to the space of positions and orientations, where crossing structures are disentangled based on their local orientation (cf. Fig.~\ref{fig:20230922:comparisonR2vsSE2lines:ContinuousSet},~\ref{fig:20230922:comparisonR2vsSE2chip:ContinuousSet}). 

After the lifting step, we apply our connected component algorithm on the Lie group $\SE{2}$.
This algorithm is built for any Lie group $G$ and uses theory involving (logarithmic approximations of) Riemannian distances. This results in the practical advantage of grouping different, possibly overlapping, anisotropic connected components based on their local alignment (measured by a Riemannian distance).

\begin{figure*}[ht!]
    \centering
    \begin{subfigure}[t]{0.49\textwidth}
        \includegraphics[width=\textwidth]{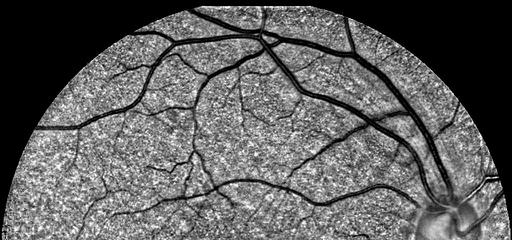}
        \caption{Reference Image.}
        \label{fig:20230922:comparisonR2vsSE2vessels:ReferenceImage}
    \end{subfigure}\hfill
    \begin{subfigure}[t]{0.49\textwidth}
        \includegraphics[width=\textwidth]{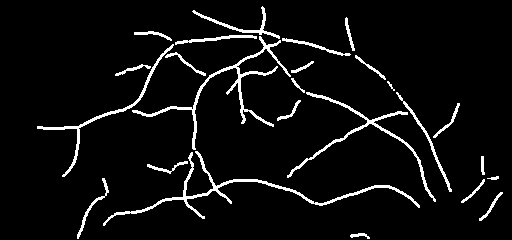}
        \caption{Reference Segmentation.}
        \label{fig:20230922:comparisonR2vsSE2vessels:ReferenceSegmentation}
    \end{subfigure}
    \begin{subfigure}[t]{0.49\textwidth}
        \includegraphics[width=\textwidth]{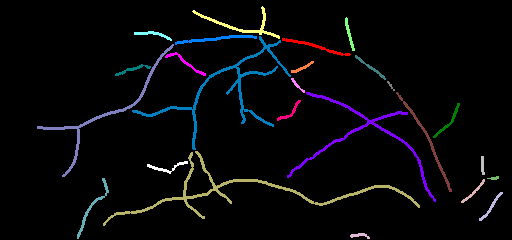}
        \caption{Connected components in $\mathbb{R}^2$.}
        \label{fig:20230922:comparisonR2vsSE2vessels:ClassicalR2}
    \end{subfigure}\hfill
    \begin{subfigure}[t]{0.49\textwidth}
        \includegraphics[width=\textwidth]{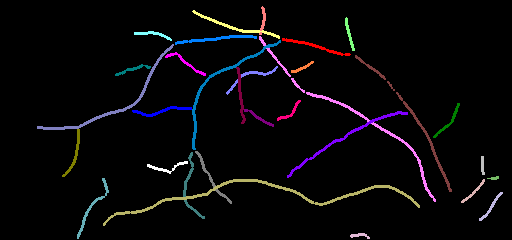}
        \caption{Connected components in \SE{2}.
        }
        \label{fig:20230922:comparisonR2vsSE2vessels:deltaCCinSE2}
    \end{subfigure}
    \caption{Visualization of the effect of the connected component algorithm in $\mathbb{R}^2$ and in the Lie group \SE{2}. The classical connected component algorithm is not able to differentiate between crossing vessels, and additionally breaks up not perfectly connected vessels, into different components. The algorithm presented in Sec.~\ref{sec:20230822:Algorithm} using parameters $(w_1,w_2,w_3)=(0.2,1,4)$ for the left-invariant metric \eqref{eq:MTFweights}, can better differentiate between different crossing structures and can additionally group well-aligned structures resulting in a more intuitive result.}
    \label{fig:comparisonR2vsSE2vessels}
\end{figure*}

Once all connected components have been identified, the output is projected back onto the input image in $\mathbb{R}^2$, as visualized in Fig.~\ref{fig:20230922:comparisonR2vsSE2:deltaCCinSE2},~\ref{fig:20230922:comparisonR2vsSE2chip:deltaCCinSE2}.
The new algorithm cannot only deal with overlapping structures, but can also deal with small interruptions of lines, still assigning them to the same connected component, as shown in Fig.~\ref{fig:comparisonR2vsSE2vessels}.

\subsection{Fundamentals of Connected Component Algorithms}

Now that we have illustrated
the main geometric idea on basic examples, let us zoom out, and consider the vast literature on connected components analysis.
There are typically 
two main approaches to define connectivity: 1) the graph-theoretic approach and 2) the topological approach. In the first, one considers images as a set of pixels connected through a graph \cite{diestel1997Graph,roerdink2000watershed,VanDeGronde2015PathBased,Heijmans1999connected}; in the second, images are rather seen as continuous functions on a smooth manifold \cite{dugundji1966topology,angulo2014Riemannian,velascoForero2013nonlocal,Burgeth2004Morphological}. This work belongs to the second category. There are two main subclasses of topology-based algorithms: grey-level mathematical morphology \cite{Heijmans1991Theoretical,WEICKERT2001Efficient} and binary morphology \cite{vandenBoomgaard1992methods,deNatale2017Detecting,BekkersNilpotent,abbasi2016geometric}. In this article, we use techniques from grey-level mathematical morphology, such as dilation, to solve binary morphology applications.

Connected components have received a lot of attention over the years, resulting in the development of various techniques and algorithms.
\begin{itemize}
    \item \underline{Pixel connectivity:} Initially, people were interested in identifying connected components in images and sets. They tackled the problem by looking at the direct neighbors of a pixel or voxel in the set, i.e. either sharing an edge or corner in 2D, or sharing a face, edge, or corner in 3D \cite{Salembier2009Connected}. This results in a very local analysis of a given input.
    \item \underline{Distance inclusion:} This local connectivity approach is not always desirable depending on the application. Therefore, techniques were developed where the connectivity was defined in terms of the distances between neighboring points, using mathematical morphology \cite{haas1967morphologie}
    . This kind of connectivity allows to assign non-adjacent pixels to the same connected component. This is useful when one wants to identify multiple, previously separate connected components, as one, e.g. in the identification of paw prints or words \cite{georgios2009generalized}.
    \item \underline{Symmetry inclusion:} The grouping of line elements should be equivariant under roto-translations. This is guaranteed if the metric is left-invariant. Note that for isotropic distances, the equivariance under roto-translations is always satisfied.
\end{itemize}

\subsection{Topological Data Analysis}
Although we focus on morphological (PDE-based) data analysis in this article, similar techniques were developed in the field of topological data analysis \cite{dey2022computational,boissonnat2018geometric,edelsbrunner2010computational}. They developed a profound theory and a variety of methods to identify clusters in point clouds, often relying on distances between points \cite{Chazal2013persistenceBased,Cao2025kmeans,Beksi20163dpoint}. In topological data analysis, clusters are created based on features of the point cloud, and the optimal threshold is chosen based on a persistence diagram \cite{Edelsbrunner2002topological,Zomorodian2005computing,Chazal2013persistenceBased,Skraba2010persistenceBased}. We apply a similar approach in Section \ref{sec:20230922:Persistence}. Topological data analysis typically focuses on the Euclidean and Riemannian settings on manifolds in general \cite{niyogi2008finding,Tinarrage2023recovering,attali2024tight}, but has not yet optimized its methods to Lie groups.

In the last decade, inclusion of efficient Lie group techniques and PDE-analysis has considerably improved geometric tracking \cite{bekkers2015pde,duits2018optimal}, denoising and image enhancement \cite{portegies2015improving,CittiSanguinetti,sherry2025diffusion}, inpainting \cite{Gauthier}, and geometric deep learning \cite{Smets2023PDEBased,bekkers2018rotoTranslation,cohen2016group}.
As we show in the Appendix~\ref{app:ToMATo}, topological data analysis techniques, such as ToMATo \cite{Chazal2013persistenceBased}, do not include (sub-Riemannian) methods of neurogeometric perceptual organization on Lie groups \cite{CittiSanguinetti,Barbieri,DuitsAssociationFieldsViaCusplessSRGeodesicsinSE2,Bellaard2023analysis,petitot2017elements}, which we show to be crucial for tackling our specific applications.

\subsection{Connectivity on the Lie Group \SE{2}}
In this article, we focus on dealing with overlapping structures by lifting an image from $\mathbb{R}^2$ to the Lie group $\SE{2}$. We highlight some of the works that have focused on this in the past.

Van de Gronde et al \cite{VanDeGronde2015PathBased} suggested a graph-based approach. There, a local orientation tensor was identified for each vertex in the graph. This allowed for the identification of all possible paths containing every vertex while satisfying a set of constraints (local orientations of two adjacent pixels sufficiently aligned). The resulting set of possibly connected vertices represented one connected component. Since this method only considers adjacent pixels for connectivity, this is a local approach in the lifted space of rototranslations.

Besides the graph-based approach in the lifted space, different algorithms based on perceptual grouping (grouping of regions and parts of the visual scene to get higher order perceptual units such as objects or patterns \cite{Brooks2015PerceptualGrouping}) have been introduced to identify vascular trees (in the lifted space) \cite{BekkersNilpotent,abbasi2016geometric}. In \cite{abbasi2016geometric}, the connectivity is learned from retinal images, after which the learned convection-diffusion kernel is used to determine the connectivity of the vessel fragments. On the other hand, in \cite{BekkersNilpotent}, a set of key points located on blood vessels are used as input for the algorithm. Then, the key points are connected based on their mutual distance. It is important to note that these algorithms do not rely on mathematical morphology.

\subsection{Our approach: Connected Components by PDE-Morphology in the Lifted Space $\SE{2}$}
We aim for an approach that relies on mathematical morphology techniques to identify connected components. We consider Lie groups $G$ equipped with a left-invariant metric tensor field $\mathcal{G}$
.
This metric tensor field induces a left-invariant distance $d_\mathcal{G}$ on the Lie group. The precise formulas for the metric $d_\mathcal{G}$ and metric tensor field $\mathcal{G}$ will follow in Section~\ref{sec:LieGroupsLeftInvariantDistances} (cf. \eqref{eq:distance},\eqref{eq:MTFweights}).

\begin{figure*}[ht!]
    \centering
    \begin{subfigure}[t]{0.3\textwidth}
        \includegraphics[width=\textwidth]{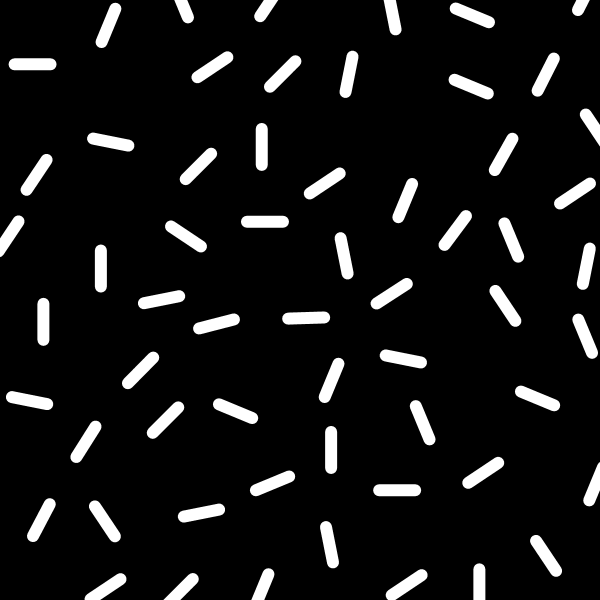}
        \caption{Reference Image containing \textbf{50} line segments.}
        \label{fig:20230922:InflMTF:ReferenceImage}
    \end{subfigure} \hfill
    \begin{subfigure}[t]{0.3\textwidth}
        \begin{picture}(160,135)
            \put(10,20){\includegraphics[width=.85\textwidth]{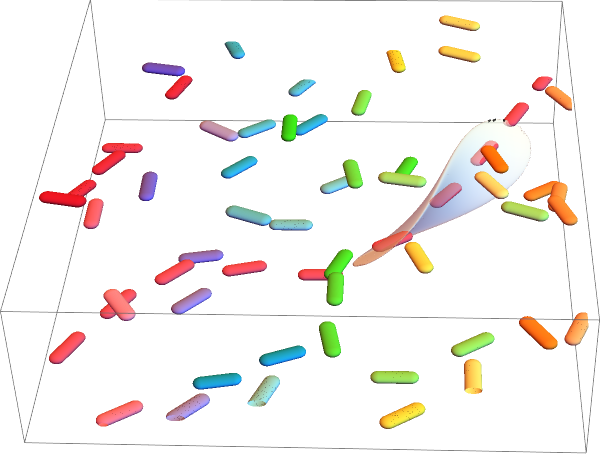}}
             \put(0,30){\small\textcolor{black}{$\theta$}}
        \end{picture}
        \caption{Connected components when using an anisotropic metric tensor field with parameters \mbox{$(w_1,w_2,w_3)=(0.05,1,2)$}, anisotropy factor 20.\\
        1-connected component count: \textbf{38}.}
        \label{fig:20230922:InflMTF:AnisotropicMetricTheta}
    \end{subfigure}\hfill
    \begin{subfigure}[t]{0.3\textwidth}
        \includegraphics[width=\textwidth]{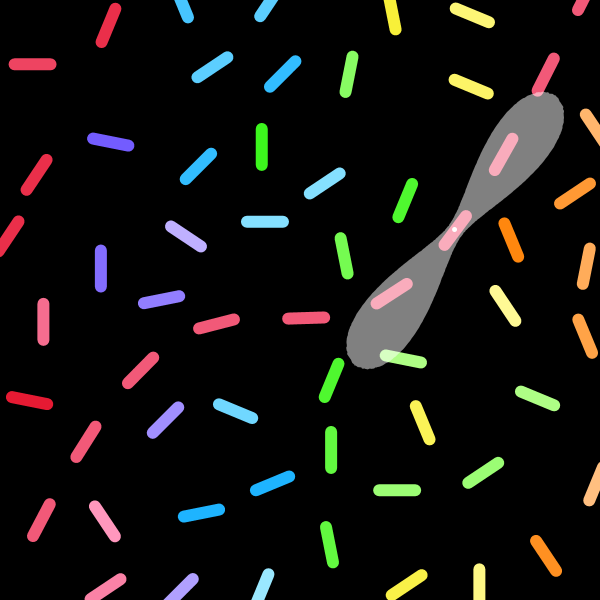}
        \caption{Connected components when using an anisotropic metric tensor field with parameters \mbox{$(w_1,w_2,w_3)=(0.05,1,2)$}.\\
        1-connected component count: \textbf{38}.}
        \label{fig:20230922:InflMTF:AnisotropicMetricThetaR2}
    \end{subfigure}
    \begin{subfigure}[t]{0.3\textwidth}
        \caption*{}
    \end{subfigure}
    \hfill
    \begin{subfigure}[t]{0.3\textwidth}
        \begin{picture}(160,135)
            \put(0,0){\includegraphics[width=\textwidth]{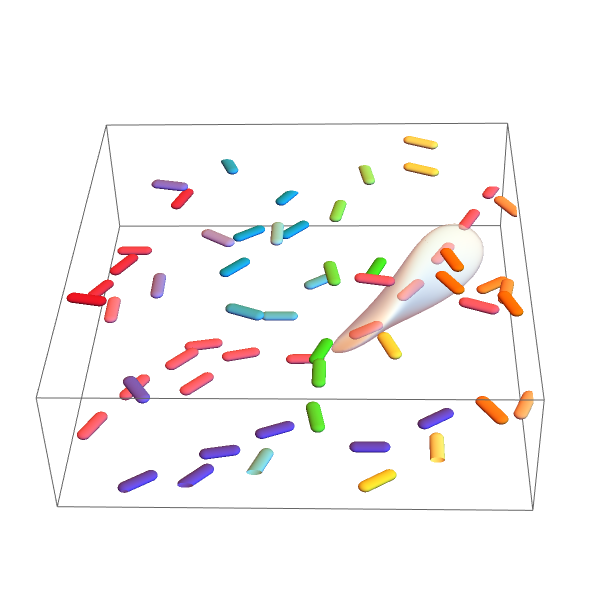}}
             \put(0,30){\small\textcolor{black}{$\theta$}}
        \end{picture}
        \caption{Connected components when using an anisotropic metric tensor field with parameters \mbox{$(w_1,w_2,w_3)=(0.05,0.3,2)$}, anisotropy factor 6.\\
        1-connected component count: \textbf{32}.}
        \label{fig:20230922:InflMTF:AnisotropicA2}
    \end{subfigure}
    \hfill
     \begin{subfigure}[t]{0.3\textwidth}
        \includegraphics[width=\textwidth]{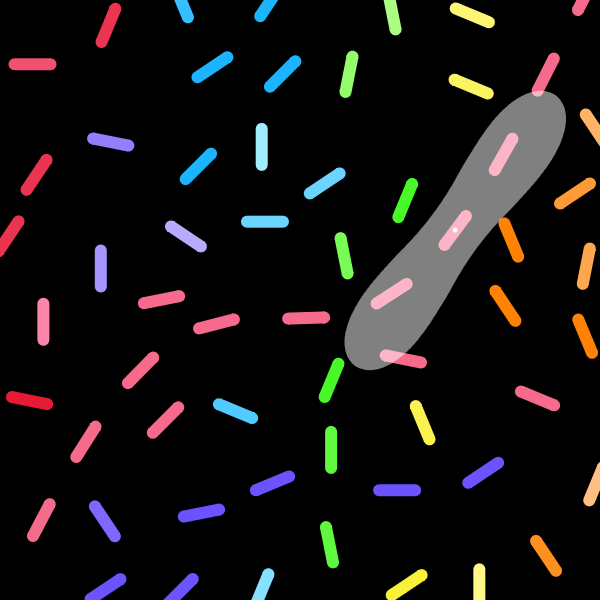}
        \caption{Connected components when using an anisotropic metric tensor field with parameters \mbox{$(w_1,w_2,w_3)=(0.05,0.3,2)$}.\\
        1-connected component count: \textbf{32}.}
        \label{fig:20230922:InflMTF:AnisotropicA2R2}
    \end{subfigure}
    \begin{subfigure}[t]{0.3\textwidth}
        \caption*{}
    \end{subfigure}
    \hfill
    \begin{subfigure}[t]{0.3\textwidth}
        \begin{picture}(160,135)
            \put(0,0){\includegraphics[width=\textwidth]{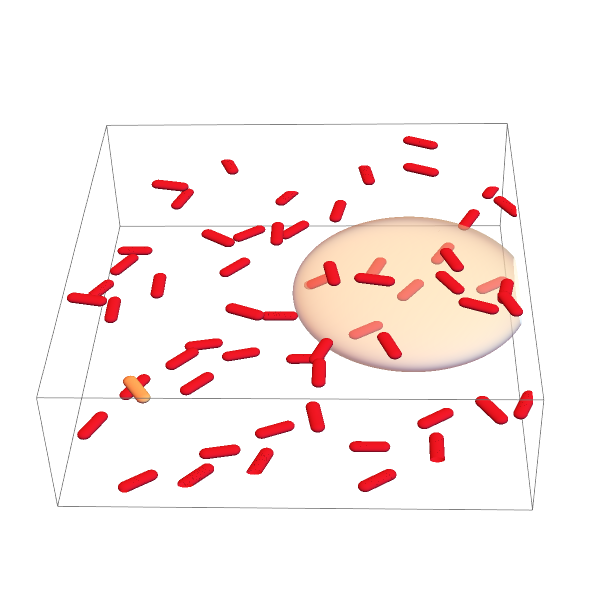}}
             \put(0,30){\small\textcolor{black}{$\theta$}}
        \end{picture}
        \caption{Connected components when using an isotropic metric tensor field with parameters \mbox{$(w_1,w_2,w_3)=(0.05,0.05,2)$}, anisotropy factor 1.\\
        1-connected component count: \textbf{2}.}
        \label{fig:20230922:InflMTF:IsotropicA2}
    \end{subfigure}
    \hfill
     \begin{subfigure}[t]{0.3\textwidth}
        \includegraphics[width=\textwidth]{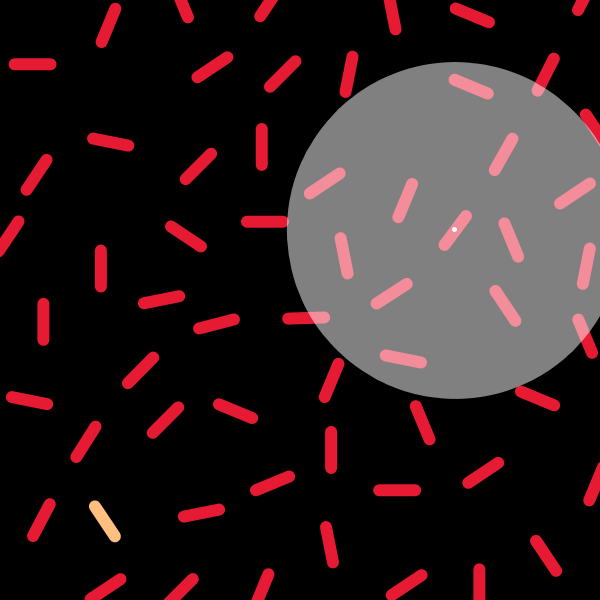}
        \caption{Connected components when using an isotropic metric tensor field with parameters \mbox{$(w_1,w_2,w_3)=(0.05,0.05,2)$}.\\
        1-connected component count: \textbf{2}.}
        \label{fig:20230922:InflMTF:IsotropicA2R2}
    \end{subfigure}
    \caption{Visualization of the effect of the metric tensor field on the connected component algorithm in the Lie group \SE{2}. The output heavily relies on the chosen distance metric $d_\mathcal{G}$ introduced in \eqref{eq:distance}, with metric tensor parameters as defined in \eqref{weights}.}
    \label{fig:InfluenceMetricTensorField}
\end{figure*}

In the specific case $G=SE(2)$, the group of roto-translations, the metric tensor field $\mathcal{G}$ is determined by three parameters $w_1$, $w_2$, and $w_3$. Intuitively, these parameters represent tangential, lateral, and angular movement costs respectively.
The influence of different choices for these parameters is illustrated in Fig.~\ref{fig:InfluenceMetricTensorField}, where the reference image is lifted to the space of roto-translations. One single Riemannian ball is plotted in each figure (lifted to $\SE{2}$ and projection on $\mathbb{R}^2$) to illustrate the influence of the metric tensor field. 

The introduced algorithm will identify, for a given threshold $\delta>0$, the largest cover of $I$ such that all sets are separated by at least distance $\delta$. We call every set in this cover a `\deltaConnected\ component'. The threshold distance $\delta>0$ is a parameter to be tuned by the user. 

\begin{remark}[$\delta$-connected components and \v{C}ech/Vietoris-Rips complexes of radius $\delta/2$]\label{rem:deltaCCvsCechComplexes}\\
\deltaConnected\ components of discrete point clouds arise as connected components of the vertices of non-singleton elements in \v{C}ech (and Vietoris-Rips) complexes of radius $\delta/2>0$. Formally, non-singleton elements in Cech complexes are simplices (not points). The identification of \deltaConnected\ components and \v{C}ech complexes also generalizes from discrete point clouds to (para)compact sets \cite{edelsbrunner2010computational,boissonnat2018geometric,dey2022computational,Kim2020homotopy} (via nerves of $\delta$-balls centered at points of the set). In our setting, where we only consider compact subsets of Lie groups, and in our analysis, the basic notion of \deltaConnected\ component suffices.
\end{remark}

In calculating these \deltaConnected\ components, we rely on formal morphological convolutions via Hamilton-Jacobi-Bellman (HJB) equations and Riemannian distance approximations, allowing for fast, parallelizable morphological convolution algorithms.

In the experimental section of this article, we have restricted ourselves to the Lie group $G=\SE{2}$. We consider two-dimensional images which we view as functions $f:\mathbb{R}^2\to \mathbb{R}$ belonging to the space $\mathbb{L}_2(\mathbb{R}^2, \mathbb{R})$. We next lift the images to the space of roto-translations $\SE{2}$ with a mapping
\begin{align*}
W_\phi:\mathbb{L}_2(\mathbb{R}^2, \mathbb{R})\to\mathbb{L}_2(\SE{2}, \mathbb{R}),
\end{align*}
where, for every $f\in \mathbb{L}_2(\mathbb{R}^2, \mathbb{R})$,
\begin{align}
(W_\phi f)(g)\coloneqq \int_{\mathbb{R}^2}\phi\left(R_\theta^{-1}(\y-\x)\right)f(\y)\mathrm{d}\y\qquad\quad \nonumber\\ \forall g=(\x,\theta)\in \SE{2},\label{eq:orientationscoreWphif}
\end{align}
and where $\phi\in \mathbb{L}_1(\mathbb{R}^2)\cap\mathbb{L}_2(\mathbb{R}^2)$ is a  rotating anisotropic wavelet. 
In this article, we use the real-valued cake wavelets (for more information, see \cite{InvertibleOrientationScores,InvertibleOrientationScores2}). The resulting function $W_\phi f$ is usually called the orientation score of the image $f$, and the operator $W_\phi$ is called the orientation score transform.

The \deltaConnected\ components are calculated on suitably binarized orientation scores $W_\phi f$, after which the resulting components are ``projected back'' onto $\mathbb{L}_2(\mathbb{R}^2,\mathbb{R})$ by inverting $W_\phi$ \cite{PhDDuits}.

In the experimental section, we see that the \deltaConnected\ component algorithm allows us to identify aligned structures in images. In the experiments, we look at both artificial images (cf. Figs.~\ref{fig:20230922:comparisonR2vsSE2:Lines},\ref{fig:20230922:comparisonR2vsSE2:chip},\ref{fig:InfluenceMetricTensorField},\ref{fig:PersistenceConnectedComponent},\ref{fig:PersistenceConnectedComponentonG}) and images of the human retina (Figs.~\ref{fig:comparisonR2vsSE2vessels},\ref{fig:CCIm5}-\ref{fig:CCIm24}). On both datasets, we can identify connected components that do not suffer from the crossing structures. We have also performed a stability analysis on the \deltaConnected\ components as illustrated in Fig.~\ref{fig:PersistenceConnectedComponent},\ref{fig:PersistenceConnectedComponentonG}. Lastly, we look at grouping different components in images of the human retina based on their mutual affinity (linked to perceptual grouping). By doing this, we can find even more complete vascular trees.

As the approach generalizes to other Lie groups and applications, cf. \cite{Barbieri,DuitsACHA,mashtakov2017tracking,PhDForssen,Duits2007Image}, we formulate the theory beyond the $\SE{2}$ setting. Crucial ingredients are a good logarithmic norm approximation (valid for $G=H(d)$, $SO(d)$, $SE(d)$) and a lift via a unitary group representation so that a left action on the image is a left action on the score \cite{PhDDuits}.

\subsection{Main Contributions}
The main contributions of this paper are:
\begin{enumerate}
    \item We develop a new method to find the connected components of a binary function defined on a Lie group G equipped with a left-invariant (sub)-Riemannian metric. The method is roto-translation equivariant and relies on standard tools from topological data analysis (e.g. \v{C}ech complexes of radius $\delta/2$ \cite{edelsbrunner2010computational,boissonnat2018geometric,dey2022computational,Chazal2013persistenceBased,niyogi2008finding,attali2024tight,Tinarrage2023recovering} are closely related to our `\deltaConnected\ components’). 
    We compute them via morphological PDEs on Lie groups. Here, we employ efficient left-invariant solvers, relying on iterative morphological group convolutions with analytic PDE kernels. 
    Such morphological group convolutions are parallelizable and employ the group structure on the Lie group $G$.

    \item We mathematically analyze our method:
    \begin{enumerate}
        \item We prove that the provided algorithm always converges to the correct \deltaConnected\ component in a finite number of steps in Theorem~\ref{th:20230605:ConvergenceConnectedComponentAlgorithm}.
        \item We show reflectional symmetries of our connected component algorithm (Cor.~\ref{cor:ReflectionalSymmetries} in App.~\ref{app:ReflectionalSymmetries}) and how this generalizes to other Lie groups $G$ of dimension $n<\infty$ (Lemma~\ref{lemma:normInvarianceReflections}). This fundamental property is due to the invariance of both the Riemannian distance and its logarithmic norm approximation under the $2^n$ reflectional symmetries in the Lie algebra. For $G=SE(2)$ and $G=SO(3)$ inclusion of these symmetries (cf.~ Fig.~\!\ref{fig:ReflectionalSymmetries}) is desirable in the connected component algorithm. 
    \end{enumerate}
    \item Along with this method, we publish code (in Mathematica and Python \cite{Berg2024connectedNotebooks}) for the morphological convolutions and the connected component algorithm.  
    \item The iterative morphological convolutions are parallelizable, fast and flexible, thanks to intuitively parameterized analytic kernels on the Lie group $G$. They do not require more expensive state-of-the-art anisotropic fast-marching schemes \cite{MirebeauPortegies,mirebeauFM} for computing the Riemannian distance maps.
    \item We show how our method, combined with (variants of) standard methodology from topological data analysis (affinity matrices \cite{abbasi2016geometric,AbbasiSureshjani}, persistence homology-based clustering \cite{Skraba2010persistenceBased,Chazal2013persistenceBased}), is very beneficial in multi-orientation image analysis of complex vascular trees in retinal imaging. We present experiments of (improved) grouping and segmentation of blood vessels.
\end{enumerate}
Generally, our algorithm performs well. It can be applied to multi-orientation analysis in flat and spherical images, relating to resp. Lie group case $SE(2)$ \cite{berg2023geodesicEnhancements} and Lie group case $SO(3)$ \cite{berg_sherry2024SO3}. 
In both cases, we rely on analytic distance approximations as explained in Appendix~\ref{app:20230922:LogarithmicNormApproximation}. 
A basic comparison to existing approaches can be found in Appendix~\ref{app:comparisonSAMToMATo}.

Our method also has two main limitations: 1)  we did not fully employ the possible parallelization of our Lie group processing, and 2) we did not perform automatic optimization of the parameters.
Therefore, in future work, we aim for: 1) faster GPU implementations via TaiChi (as done for geodesic tracking in \cite{berg_sherry2024SO3}) and 2) training of the $n=\textrm{dim}(G)$ parameters of the connected component algorithm (controlling the ball shapes in $G$, recall Fig.~\!\ref{fig:InfluenceMetricTensorField} on SE(2)).

\subsection{Structure of the Article}
First, we briefly introduce Lie groups and left-invariant distances in Sec.~\ref{sec:LieGroupsLeftInvariantDistances}. In Sec.~\ref{sec:DeltaConnectedComponents}, we introduce our notion of \deltaConnected\ components and elaborate on our approach towards identifying those components. Then, we introduce the morphological convolutions (Sec.~\ref{sec:MorphologicalDilations}) that are used in the presented \deltaConnected\ component algorithm (Sec.~\ref{sec:20230822:Algorithm}). We prove it converges to the desired components in a finite number of steps (Sec.~\ref{sec:20230822:Analysis}). Additionally, we explain the choice of $\delta$ using persistence diagrams in Sec.~\ref{sec:20230922:Persistence}. We define the affinity between different \deltaConnected\ components in Sec.~\ref{sec:20230922:Affinity}, which allows us to quantify how well different components are aligned. Once all these theoretical concepts have been introduced, we move on to Sec.~\ref{sec:20230922:Experiments}, where we present the experimental part of the article. We summarize and conclude in Sec.~\ref{sec:20230922:DiscussionConclusion}.

%% file: Review_JMIV/2_1LieGroups.tex
\section{Background on Lie Groups and Left-Invariant Distances}\label{sec:LieGroupsLeftInvariantDistances}
The theory presented in this document applies to all finite-dimensional Lie groups $G$ of dimension $n\in\mathbb{N}$, with a left-invariant Riemannian distance $d$. We recall that a Lie group is a smooth manifold $G$ equipped with a group structure such that the group multiplication $\mu: G\times G\to G$ given by $\mu(x,y)=xy$ and the inverse map $\nu: G\to G$ given by $\nu(x)=x^{-1}$ are smooth maps.

The tangent space at $h\in G$ of Lie group $G$ is the vector space of all tangent vectors of curves passing through $h\in G$, which we denote by $T_h(G)$.

Let $g\in G$. The left-action on $G$ is denoted by $L_g: G\to G$ and defined as $L_g(h)=gh$, and its corresponding push-forward $(L_g)_*:T_h(G)\to T_{gh}(G)$ is defined as $(L_g)_* A_h u = A_h (u\circ L_g)$ for all smooth functions $u: G\to \mathbb{C}$, where $A_h\in T_h(G)$. In particular if we take $h=e$, the unit element, and a coordinate basis $\left\{\left.\partial_{x^i}\right|_e\right\}_{i=1}^n$ at $T_e(G)$, then 
\begin{align*}
    (L_g)_* \left.\partial_{x^i}\right|_{e} u = \left.\partial_{x^i}\right|_{e} (u\circ L_g).
\end{align*}
The Riemannian metric tensor field is a field of inner products $\mathcal{G}_g: T_g(G)\times T_g(G)\to\mathbb{R}$ for each $g\in G$. We always consider a left-invariant Riemannian metric tensor field $\mathcal{G}$, i.e. for all $g,h\in G$ and for all $\dot{h}\in T_h(G)$
\begin{align}
   \mathcal{G}_{gh}\left((L_g)_* \dot{h},(L_g)_*\dot{h}\right)=\mathcal{G}_{h}\left(\dot{h},\dot{h}\right).\label{eq:20230822:LeftInvariantMetricTensorField}
\end{align}
Then, the Riemannian left-invariant distance $d_\mathcal{G}$ is defined by
    \begin{align*}
        d_\mathcal{G}(g,h)\coloneqq\inf\left\{\left.\int_0^1\sqrt{\mathcal{G}_{\gamma(t)}\left(\dot{\gamma}(t),\dot{\gamma}(t)\right)}\mathrm{d}t\;\right|\right.\qquad\qquad\\
        \left.\;\gamma\in\Gamma_1,\gamma(0)=g,\gamma(1)=h\right\},\numberthis\label{eq:distance}
    \end{align*}
where $\Gamma_1=\text{PC}([0,1],G)$ denotes the family of piece-wise continuously differentiable curves.
Then from \eqref{eq:20230822:LeftInvariantMetricTensorField}, it follows readily that  $d_{\mathcal{G}}(gh_1,gh_2)=d_{\mathcal{G}}(h_1,h_2)$ for all $g,h_1,h_2\in G$. Intuitively, this means that the distance is invariant to left-actions, and one has $d_{\mathcal{G}}(g,h)=d_{\mathcal{G}}(h^{-1}g,e)$, where $e\in G$ is the unit element.
\begin{remark}
    The distance $d_{\mathcal{G}}$ always relies on a metric tensor $\mathcal{G}$. Henceforth, we omit the label $\mathcal{G}$, and write $d\coloneqq d_{\mathcal{G}}$ for the distance whenever we do not need to stress the metric tensor used to calculate the distances.
\end{remark}

\begin{remark}[Explanation weights of metric tensor field]
    From \eqref{eq:20230822:LeftInvariantMetricTensorField}, one can deduce that all left-invariant (recall Eq.~\!(\ref{eq:20230822:LeftInvariantMetricTensorField})) metric tensor fields are given by 
    \begin{align}
        \mathcal{G}_g(\dot{g},\dot{g})=\sum\limits_{i,j=1}^n c_{ij}\,\dot{g}^i\dot{g}^j\label{eq:MTFweights}
    \end{align}
    with \emph{constant} coefficients $c_{ij}\in\mathbb{R}$, and where $\dot{g}= \sum_{i=1}^n \dot{g}^i \mathcal{A}_{i}|_{g} \in T_g(G)$ for a basis $\{\mathcal{A}_{i}\}_{i=1}^n$ of left-invariant vectors $G \ni g \mapsto \mathcal{A}_{i}|_g \in T_g(G)$, $i=1,\ldots,n$. Note that in general the matrix $[c_{ij}]$ is positive symmetric definite, but 
    in this article
    we constrain ourselves to the diagonal case \begin{equation} \label{weights}c_{ij}=w_i \delta_{ij}, 
    \end{equation}
    with positive weights $w_i>0$.
\end{remark}

\begin{remark}
    The theoretical results in this article also hold when the distance denotes a sub-Riemannian distance, though this requires different logarithmic approximations than \eqref{eq:logapprox} below, for details see \cite{terElst3,DuitsAMS,Bellaard2023analysis}.
\end{remark}
\begin{figure*}[ht!]
    \centering
    \begin{subfigure}[t]{0.3\textwidth}
         \centering
         \begin{picture}(160,160)
            \includegraphics[width=\textwidth]{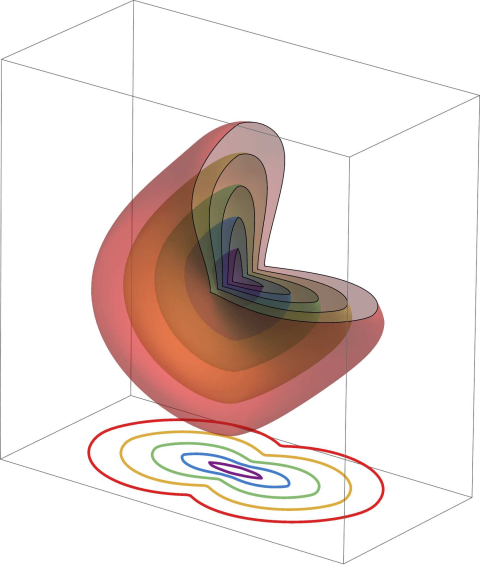}
            \put(-15,5){\textcolor{black}{$y$}}
            \put(-110,5){\textcolor{black}{$x$}}
            \put(-155,100){\textcolor{black}{$\theta$}}
         \end{picture}
        \caption{Exact Riemannian distance $d$ in G=\SE{2}  with $(w_1,w_2,w_3)=(1,4,1)$, based on \cite[Fig.5a]{Bellaard2023analysis}.}
        \label{fig:20230822:VisualizationDistanceBallsSE2:Exact}
    \end{subfigure}
    \hfill
    \begin{subfigure}[t]{0.3\textwidth}
         \centering
         \begin{picture}(160,160)
            \includegraphics[width=\textwidth]{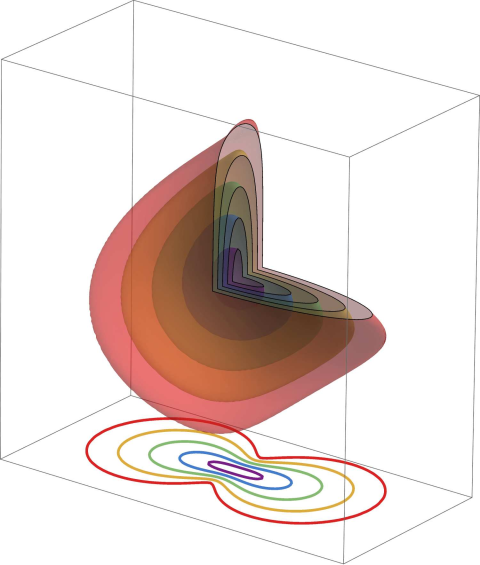}
            \put(-15,5){\textcolor{black}{$y$}}
            \put(-110,5){\textcolor{black}{$x$}}
            \put(-155,100){\textcolor{black}{$\theta$}}
         \end{picture}
        \caption{Logarithmic norm approximation $\rho_c$ of the Riemannian distance $d$ in G=\SE{2} with $(w_1,w_2,w_3)=(1,4,1)$, based on \cite[Fig.5]{Bellaard2023analysis}.}
        \label{fig:20230822:VisualizationDistanceBallsSE2:Approximate}
    \end{subfigure}
    \hfill
    \begin{subfigure}[t]{0.3\textwidth}
         \centering
         \begin{picture}(250,170)
             \put(-15,0){\includegraphics[width=1.2\textwidth]{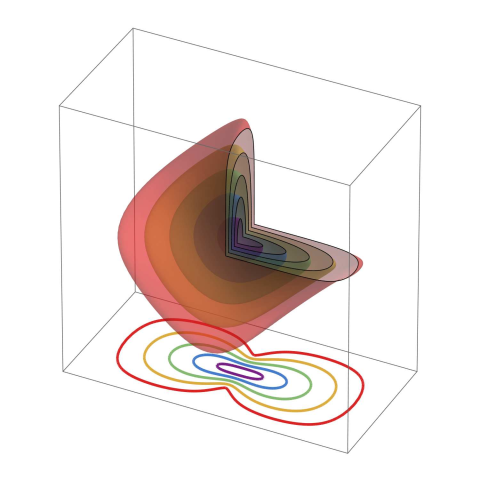}}
            \put(130,20){\textcolor{black}{$\alpha$}}
            \put(50,12){\textcolor{black}{$\beta$}}
            \put(-5,90){\textcolor{black}{$\phi$}}
         \end{picture}
        \caption{Logarithmic norm approximation in G=SO(3) with metric parameters $(\tilde{w}_1,\tilde{w}_2,\tilde{w}_3)=(1,4,1)$.}
        \label{fig:20230822:VisualizationDistanceBallsSO3}
    \end{subfigure}
    \caption{Visualization of distance balls and their logarithmic approximations in G=\SE{2} and G=SO(3). We see isocontours of $d(p_0,\cdot)$ in $G$, and on the bottom, we see the min-projection over the orientation $\theta$ of these contours. The visible contours are $d=0.5,1,1.5,2,2.5$ and the metric parameters are $(w_1,w_2,w_3)=(\tilde{w}_1,\tilde{w}_2,\tilde{w}_3)=(1,4,1)$.
    Parameter $w_1$ controls costs for tangential  motion, $w_2$ controls costs for lateral  motion
    on the base manifold (resp. $\mathbb{R}^2$ and $S^2$)
    and $w_3$ controls costs for changing the orientation along a geodesic. Thereby, together, they control the anisotropic shape of the Riemannian ball and the identification of connected components, recall white balls in column 2 of Fig.~\ref{fig:InfluenceMetricTensorField}. 
    }
    \label{fig:20230822:VisualizationDistanceBalls}
\end{figure*}
We have a particular interest in the cases of the roto-translation group
\begin{align*}
    G=SE(2)\coloneqq \mathbb{R}^2\rtimes SO(2)
\end{align*}
and the rotation group 
\begin{align*}
    G=SO(3)\coloneqq\{X\in\mathbb{R}^{3\times 3}\;|\;
    X^\top X=I,\;\det(X)=1\}.
\end{align*}
The group product of $G=SE(2)$ is given by
\begin{align}
    g_1 g_2=(x_1,R_1)(x_2,R_2)=(x_1+R_1 x_2, R_1 R_2)
\end{align}
for all $(x_1,R_1),(x_2,R_2)\in \mathbb{R}^2\rtimes SO(2)$, and the left-actions are given by roto-translations. Similarly, the group product $R_1,R_2\in SO(3)$ is given by the ordinary matrix product $R_1 R_2\in SO(3)$, where we note that
$(R_1 R_2)(R_1 R_2)^\top=R_1 R_2 R_2^\top R_1^\top=R_1 I R_1^\top= I$ and
$\det(R_1 R_2)=\det(R_1)\det(R_2)=1$.
\begin{remark}
    In the case where $G=SE(2)$, we have in \eqref{eq:MTFweights},
    $\dot{g}= \sum_{i=1}^3 \dot{g}^i \mathcal{A}_{i} \in T_g(G)$, 
    where
    \begin{equation}
    \mathcal{A}_1=\cos \theta \, \partial_x +\sin \theta \, \partial_y, \mathcal{A}_2=-\sin \theta \,\partial_x +\cos \theta \, \partial_y, \mathcal{A}_3=\partial_\theta,\label{eq:MetricTensorFrameA} 
    \end{equation}
    where we identify 
    \begin{align*}
        \theta\in \mathbb{R}/(2\pi\mathbb{Z})\leftrightarrow R_\theta=\begin{pmatrix}
            \cos\theta&-\sin\theta\\
            \sin\theta&\cos\theta
        \end{pmatrix}\in SO(2).
    \end{align*}
    In that case, $w_1$, $w_2$, and $w_3$ in \eqref{weights} describe costs for forward, sideways, and angular movement in $SE(2)$, respectively.
\end{remark}

The groups $G=SE(2)$ and $G=SO(3)$ have the special property that the exponential map $\text{Exp}: T_e(G)\to G$ is surjective and we have a complete logarithmic norm approximation 
\begin{equation} \label{eq:logapprox}d_{\mathcal{G}}(g,h)=d_{\mathcal{G}}(h^{-1}g,e)\approx \|\log h^{-1} g\|_{\mathcal{G}}
\end{equation}
for $g,h\in G$ close enough to each other, cf. Fig.~\ref{fig:20230822:VisualizationDistanceBalls}. For details on the quality of this approximation for the $G=\SE{2}$ case, see \cite{Bellaard2023analysis}. For a concise self-contained summary and explicit formulas for the logarithmic norm approximation in \eqref{eq:logapprox} in coordinates, see Eq.~\eqref{eq:logarithmicNormApproxSE2} for $G=SE(2)$ and Eq.~\eqref{eq:logarithmicNormApproxSO3} for $G=SO(3)$ in Appendix~\ref{app:20230922:LogarithmicNormApproximation}.

\begin{remark}[Inclusion of Symmetries]
    It is well-known in Lie group theory that left-invariant Riemannian metrics have a constant Gramm matrix in any left-invariant frame. So in \eqref{eq:MTFweights}, $c_{ij}\in\mathbb{R}$ is independent of $g\in G$ if $(\mathcal{A}_i)_{i=1}^n$ is a left-invariant frame. The associated distance maps of such metrics and their logarithmic norm estimates carry $2^{\textrm{dim}(G)}$ reflectional symmetries.
    In the specific Lie group case of roto-translations $G=SE(2)$ there are $2^3$ reflectional symmetries (for intuitive illustration see App.~\ref{app:ReflectionalSymmetries}) which also become reflectional symmetries of our connected component algorithm.
\end{remark}

\begin{remark}
    The Lie group $SE(2)$ can be identified with the homogeneous manifold of positions and orientations in two dimensions $\mathbb{M}_2$. However, this is not true for $n>2$: $SE(n)\not\cong \mathbb{M}_n$.
\end{remark}

%% file: Review_JMIV/2_2ConnectedComponents.tex
\begin{figure*}[ht!]
    \centering
    \begin{subfigure}[b]{0.49\textwidth}
         \centering
        \scalebox{.6}{\input {Figures/VisualizationICompact}}
        \caption{Covering of compact set $I$.}
        \label{fig:20230605:VisualizationConnectedComponentsI}
    \end{subfigure}
    \begin{subfigure}[b]{0.49\textwidth}
         \centering
        \scalebox{.6}{\input {Figures/VisualizationDeltaConnectedPoints}}
        \caption{Identification of connected components of $I$.}
        \label{fig:20230605:VisualizationConnectedComponentsPoints}
    \end{subfigure}
    \caption{The set of connected components $I$ (in blue in Fig.~\ref{fig:20230605:VisualizationConnectedComponentsI}) has a finite covering. In this case, it is covered by 10 balls with radius $\delta$ (in red), and has covering number $n_\delta(I)=10$. The distance between different connected components is larger than the radius of a ball $\delta$. The example shows three different \deltaConnected\ components, indicated by color in Fig.~\ref{fig:20230605:VisualizationConnectedComponentsPoints}. 
    }
    \label{fig:20230605:VisualizationConnectedComponents}
\end{figure*}
\section{An Algorithm to Find the \texorpdfstring{$\delta$-Connected}{delta-Connected} Components of the Compact Set \texorpdfstring{$I$}{I}}\label{sec:DeltaConnectedComponents}
To explain our algorithm, we first introduce the notion of \deltaConnected ness of a compact set $I$. This concept is closely related to the \v{C}ech complexes of radius $\delta/2$ in computational geometry and topology \cite{edelsbrunner2010computational,dey2022computational,boissonnat2018geometric}, see Remark~\ref{rem:deltaCCvsCechComplexes}.

\subsection{The Notion of  $\delta$-Connectedness}\label{sec:NotionOfDeltaConnectedness}
We start by recalling the notion of the covering number of $I$, for which we use the definition of a ball.
\begin{definition}[Riemannian Ball]
    The ball around $g\in G$ with radius $\delta>0$ is given by
    \begin{align*}
        \Ball{g}{\delta}=\left\{h\in G\;|\; d(g,h)<\delta\right\}.
    \end{align*}
\end{definition}
\begin{definition}[$\delta$-covering]
    A set $C\subset I$ is a $\delta$-covering of $I$ if $I\subseteq\bigcup \limits_{g\in C}\Ball{g}{\delta}$.
\end{definition}

\begin{definition}[Covering number]
The covering number, denoted by $n_\delta(I)$, is the smallest cardinality of all possible $\delta$-coverings of the set $I$, i.e.
\begin{align}
    n_\delta(I)\coloneqq \min\left\{\left|C\right|:\; C\text{ is a $\delta$-covering of }I\right\}.\label{eq:20230822:CoveringNumber}
\end{align}    
\end{definition}
Since we assume that the set $I$ is compact, its covering number $n_\delta(I)$ is finite. Next, we introduce the notion of \deltaConnected ness between points $g,\, h\in I$.
\begin{definition}[\deltaConnected ness of points]\label{def:20230605:Connectedness}
    Let $\delta>0$. We say that two elements $g$ and $h\in I$ are $\delta$-connected, and we denote it as $g\simdelta h$, if and only if
    there exists a finite number $m \in \mathbb{N}$ of elements $\{q_i\}_{i=0}^{m}\subset I$ such that $q_0=g$, $q_{m}=h$ and
    \begin{align*}
        d(q_{i+1}, q_i)\leq \delta, \quad \forall i \in \{0,\dots, m-1\}.
    \end{align*}
\end{definition}
\begin{remark}\label{remark:MdependsOngAndh}
    Note that $m\in\mathbb{N}$ may depend on $g,h\in G$ and $\delta$, so $m\coloneqq m_\delta(g,h)\in\mathbb{N}$.
\end{remark}
One can readily check that $\simdelta$ is an equivalence relation on $I$, and the equivalence class of a $g\in I$ under $\simdelta$, denoted by $[g]$, is defined as
\begin{align*}
    [g]\coloneqq \{ q\in I \; :\; q\simdelta g \}.
\end{align*}
These preliminary concepts allow us to give the following definition to the notion of `\deltaConnected\ components'.
\begin{definition}[\deltaConnected\ components]\label{def:20230922:deltaConnectedComponents}
The set of \deltaConnected\ components of a given set $I$ is defined as its equivalence class $\tildeIdelta:=I/\simdelta$, and each element $[g]\in \tildeIdelta$ is called a `\deltaConnected\ component'.
\end{definition}

The following result guarantees that there are at most $n_\delta(I)$ $\delta$-connected components to search for when $I$ is a compact set.

\begin{lemma}\label{lemma:20230605:maxNrdeltaConnectedComponents}
For every compact set $I$ one has that the number of \deltaConnected\ components is bounded by the covering number, i.e. $|\tildeIdelta|\leq n_{\delta}(I)$.
\end{lemma}
\begin{proof}
We argue by contradiction. Suppose that the covering number is $n_\delta(I)$, and suppose there are $m>n_\delta(I)$ connected components in $I$. Then there are $m$ points $\{g_i\}_{i=1}^m$ in $I$ such that $d(g_i, g_j)>\delta$ for all $1\leq i, j\leq m$. Therefore, we need at least $m$ balls of radius $\delta$ to cover the set $I$, so $n_\delta(I) \geq m$. However, by assumption $m>n_\delta(I)$. 
\qed 
\end{proof}

\begin{remark}
It follows from Lemma \ref{lemma:20230605:maxNrdeltaConnectedComponents} that the distance between the closest pair of points between two $\delta$-connected components is larger than $\delta$, that is,
\begin{align*}
    \inf_{(g, h) \in [g]\times [h]} d(g, h) >\delta\quad\text{if }[g]\neq[h].
\end{align*}
\end{remark}

We next define several notions of distances that will be needed in subsequent developments.
\begin{definition}[Distance from point to set]
    The distance from a point $g\in G$ to a non-empty set $A\subset G$ is given by
    \begin{align*}
        d(g,A)\coloneqq\inf\limits_{a\in A} d(g,a).
    \end{align*}
\end{definition}
\begin{definition}[Distance between two sets]\label{def:distanceBetweenTwoSets}
    The distance from a non-empty set $A\subset G$ to a non-empty set $B\subset G$ is given by
    \begin{align*}
        d(A,B)\coloneqq\inf\limits_{a\in A,b\in B} d(a,b).
    \end{align*}
\end{definition}
Additionally, we will need the concept of an $\epsilon$-thickened set when we explain the \deltaConnected\ component algorithm:
\begin{definition}[$\epsilon$-thickened set]\label{def:20230922:epsilonThickenedSet}
    For every closed subset $A\subset G$, we define its $\epsilon$-thickened version $A_\epsilon$ as the set
    \begin{align*}
        A_\epsilon \coloneqq \{ g\in G \;:\; d(g, A)\leq \epsilon \} = \bigcup\limits_{a\in A} \overline{\Ball{a}{\epsilon}}.
    \end{align*}
\end{definition}

\subsection{A General Algorithm for (Disjoint) Connected Components}\label{sec:GeneralAlgorithm}

Given a threshold $\delta>0$, the following strategy finds the $\delta-$connected components of a compact set $I$. We defer the discussion on how to select an optimal value of $\delta$ to Sec.~\ref{sec:20230922:Persistence}, and focus here on explaining the algorithm once this parameter is fixed. The strategy consists of a main algorithm \ref{alg:findAllComponents} which identifies all \deltaConnected\ components relative to the compact set $I$. It relies on an algorithm that identifies the \deltaConnected\ component $[g]\in\tildeIdelta$ of a given $g\in I$ (algorithm \ref{alg:findFullComponent}).


In practice, $\C{g}{n}$ is computed by means of
\begin{align*}
\C{g}{n}
&= \Cdelta{g}{n-1} \cap I 
=  \left(\bigcup \limits_{c\in \C{g}{n-1}}  \Ball{c}{\delta}\right) \cap I.
\end{align*}
As a result, it suffices to compute $\Ball{c}{\delta} \cap I$ for every $c\in \C{g}{n-1}$ to find the new set $\C{g}{n}$.

\begin{mdframed}
\textbf{Algorithm} \customlabel{alg:findFullComponent}{\texttt{find\_full\_component}} \vspace{0.2cm}\\
\textbf{Find $\delta$-Connected Component $[g]$ for a $g\in G$}\\
To find $[g]\in \tildeIdelta$ for a given $g$, in practice, we again proceed iteratively:
\begin{itemize}
\item $n=0$: We set $\C{g}{0}=\{g\}$.
\item $n\geq1$: Given $\C{g}{n-1}$, find
\begin{align}
\C{g}{n} = \Cdelta{g}{n-1} \cap I,\label{eq:20230822:defC_g^n+1}
\end{align}
where $\Cdelta{g}{n-1}\coloneqq\left(\C{g}{n-1}\right)_\delta$ is the $\delta$-thickened set of 
$C(g,n-1)$
(see definition~\ref{def:20230922:epsilonThickenedSet}).\\[6pt]
As we prove next in Lemma \ref{lem:convergence-test}, if $\C{g}{n} = \C{g}{n-1}$, then $\C{g}{n} = [g]$ and the algorithm terminates. Otherwise, we go to step $n+1$.
\end{itemize}
\end{mdframed}
\begin{mdframed}
    \textbf{Algorithm} \customlabel{alg:findAllComponents}{\texttt{find\_all\_components}} \vspace{0.2cm}
    \\\textbf{Find \textit{all} $\delta$-Connected Components $[g]\in \tildeIdelta$ for a Compact $I \subset G$}\\
    To find all $[g_i]\in \tildeIdelta$ in $I$, we proceed iteratively:
    \begin{itemize}
        \item $k=1$:
        \begin{enumerate}
            \item Pick any $g_1 \in G$.
            \item Compute $[g_1]$ with \ref{alg:findFullComponent} to compute \deltaConnected\ components, which is introduced in the next box.
            \item Set $\tildeIdelta(1) \coloneqq \{ [g_1] \}$.
            \item Define \mbox{$\Idelta(1) \coloneqq \{ g\in G \;:\; [g] \in \tildeIdelta(1) \}$}.\vspace{6pt}
        \end{enumerate}
        \item $k>1$: Given $\tildeIdelta(k-1)=\{[g_1], \dots, [g_{k-1}]\}$ and its associated set $\Idelta(k-1)$.
        \begin{enumerate}
            \item Search an element $g_k \in I \setminus \Idelta(k-1)$
            \begin{enumerate}
                \item if $g_k$ exists: Update 
                \begin{align*}
                    \tildeIdelta(k)\coloneqq\tildeIdelta(k-1) \cup \{[g_k]\}
                \end{align*}
                and $\Idelta(k)$ accordingly.
                \item if $g_k$ does not exist: all $k$ \deltaConnected\ components are identified, i.e. $\tildeIdelta(k-1)=\tildeIdelta$. The algorithm ends.
            \end{enumerate}            
        \end{enumerate}
    \end{itemize}
\end{mdframed}

In practice, the thickened sets $\Cdelta{g}{n}$ in \ref{alg:findFullComponent} are computed using morphological convolutions which we explain in Sec.~\ref{sec:MorphologicalDilations}.
\begin{lemma}\label{lem:convergence-test}
Let $g\in G$ and $n\geq 1$. If $\C{g}{n} = \C{g}{n-1}$, then $\C{g}{n} = [g]$.
\end{lemma}
\begin{proof}
We prove the statement by contradiction: 
Suppose $\C{g}{n}=\C{g}{n-1}$  and $\C{g}{n-1}\subsetneq[g]$. Then, there exists a $h\in [g]$ such that $h\not\in \C{g}{n-1}$ and $d(h,\C{g}{n-1})\leq \delta$ (due to the definition of the equivariance classes $[g]$), so $h\in \Cdelta{g}{n-1}$. Since we also have that $h\in I$, it follows that $h\in \Cdelta{g}{n-1}\cap I = \C{g}{n}$. Hence $h\in \C{g}{n-1}$, which is a contradiction.
\qed
\end{proof}

\begin{figure*}[ht!]
    \centering
    \begin{picture}(0,520)
        \put(-205,434){\includegraphics[width=0.118\textwidth]{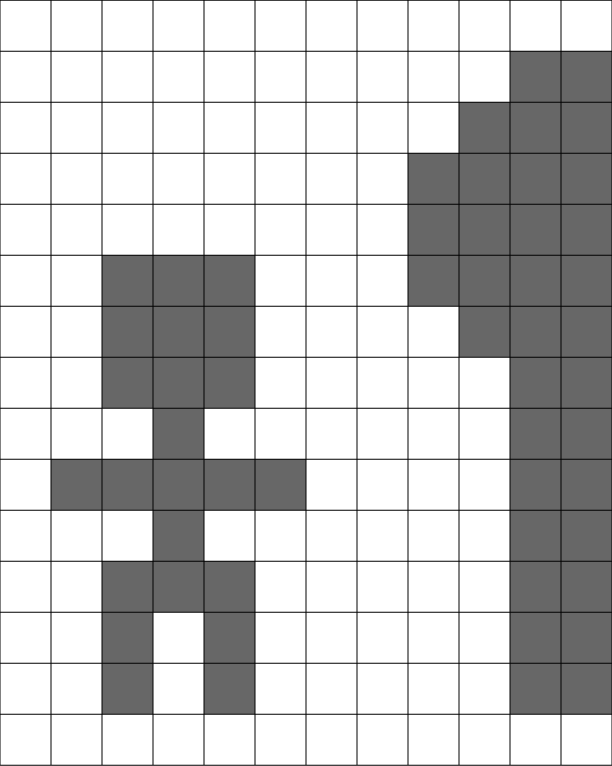}}
        \put(-191,473){\textcolor{white}{$I$}}
        \put(-157,482){\textcolor{white}{$I$}}

        \put(-205,347){\includegraphics[width=0.118\textwidth]{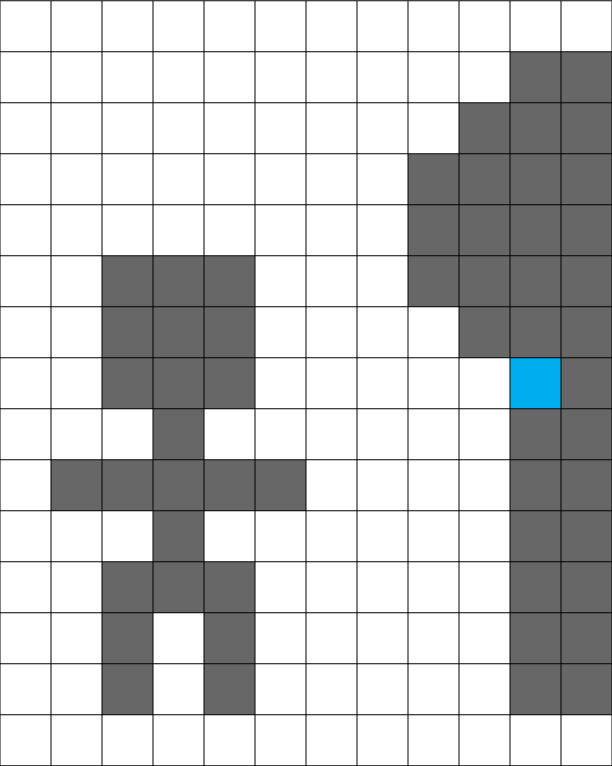}}
        \put(-142,347){\includegraphics[width=0.118\textwidth]{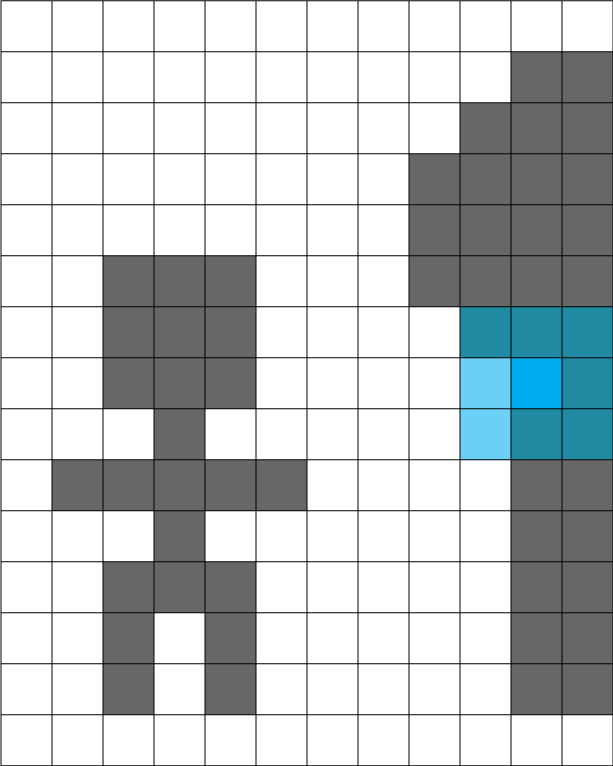}}
        \put(-80,347){\includegraphics[width=0.118\textwidth]{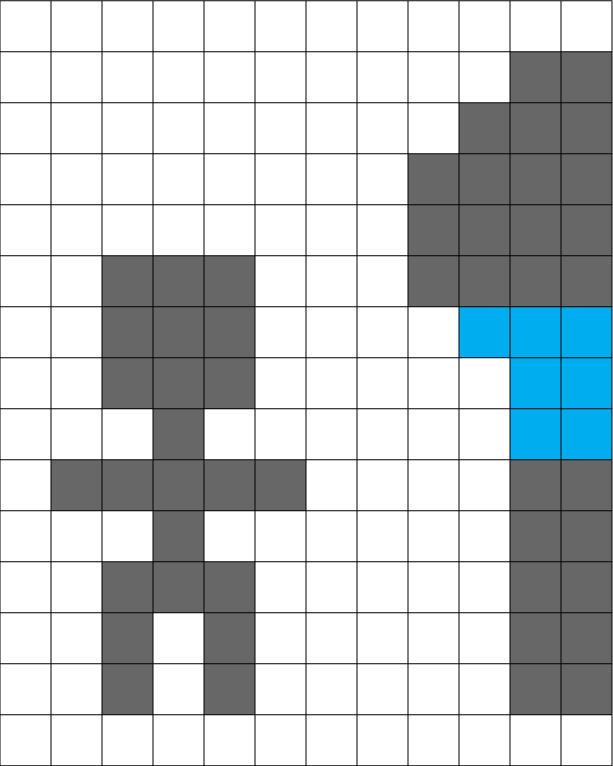}}
        \put(-17,347){\includegraphics[width=0.118\textwidth]{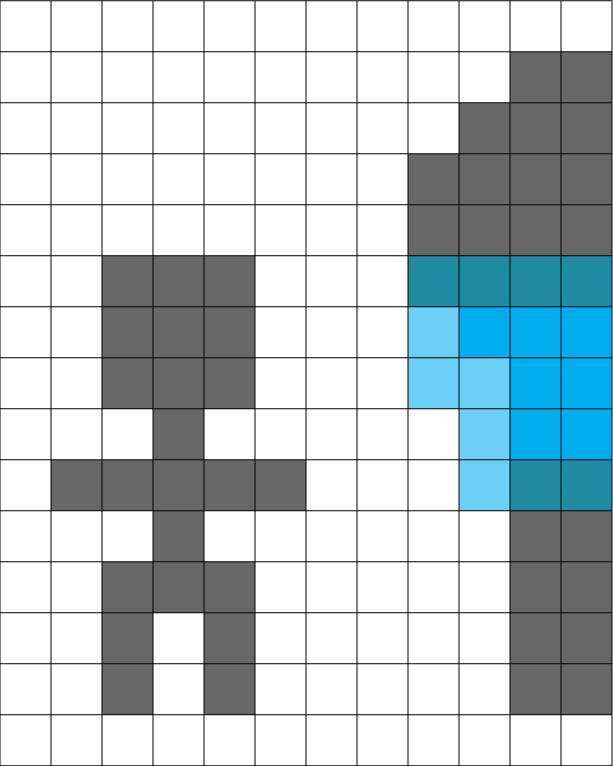}}
        \put(45,347){\includegraphics[width=0.118\textwidth]{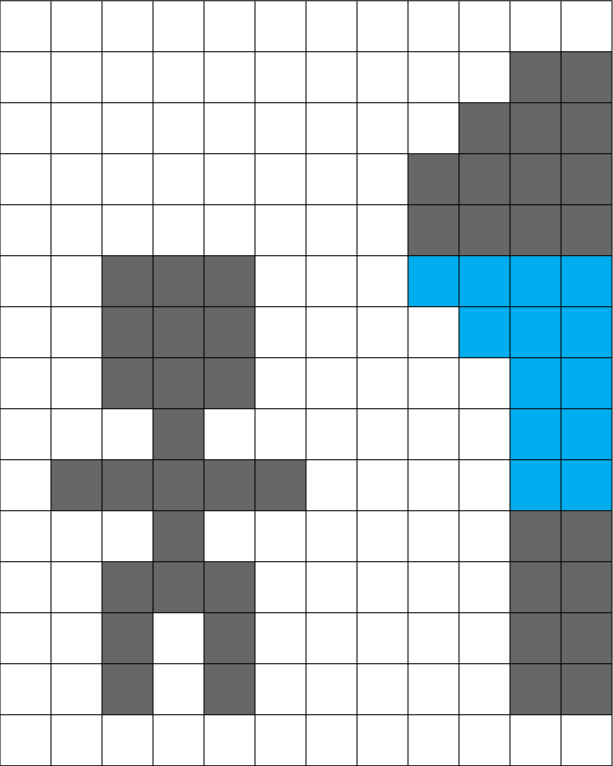}}
        \put(108,347){\includegraphics[width=0.118\textwidth]{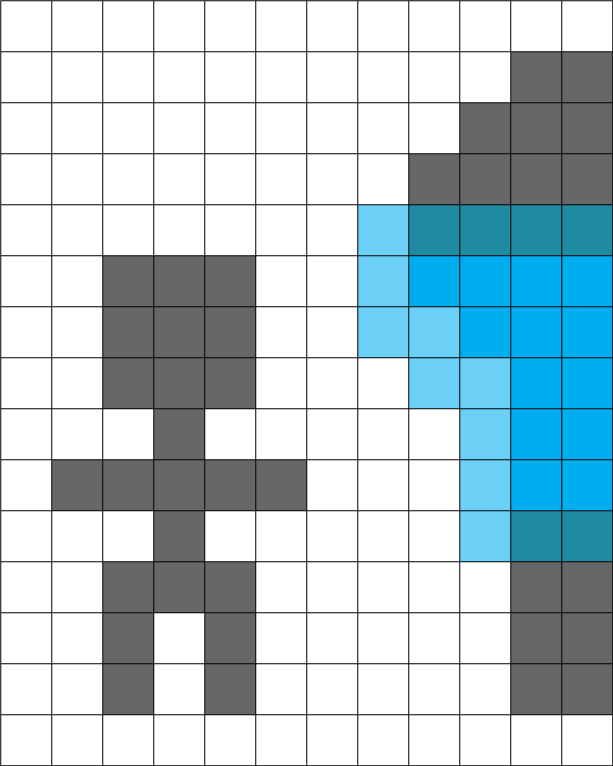}}
        \put(170,347){\includegraphics[width=0.118\textwidth]{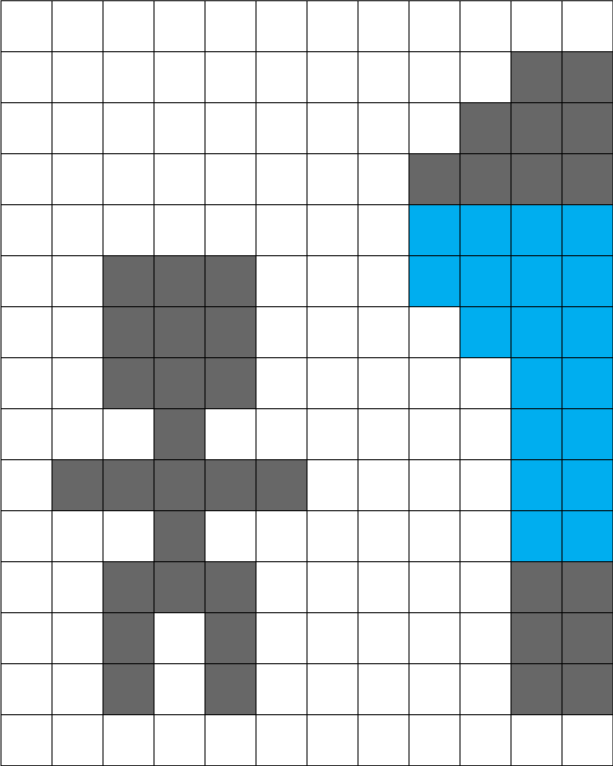}}

        \put(-142,260){\includegraphics[width=0.118\textwidth]{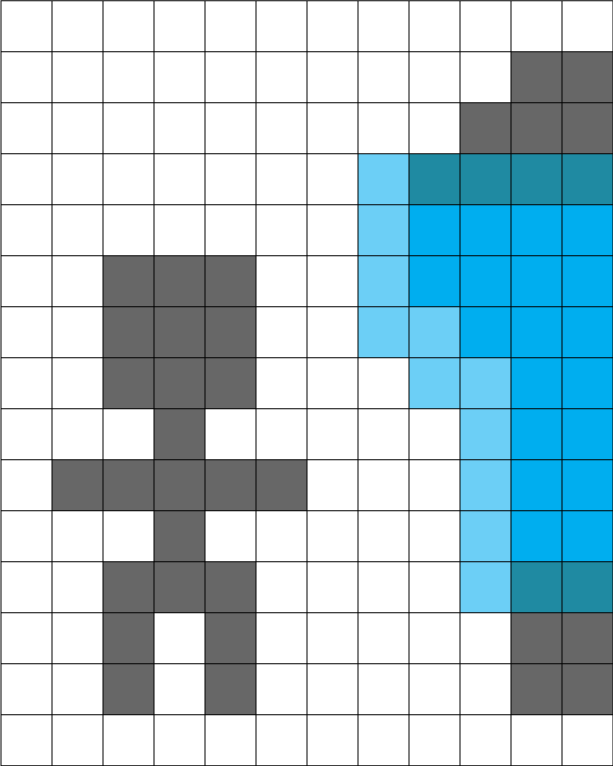}}
        \put(-80,260){\includegraphics[width=0.118\textwidth]{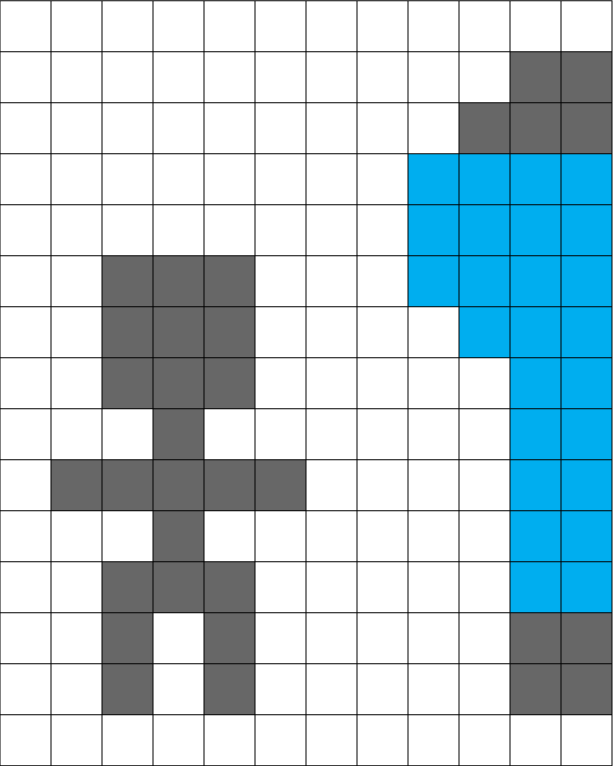}}
        \put(-17,260){\includegraphics[width=0.118\textwidth]{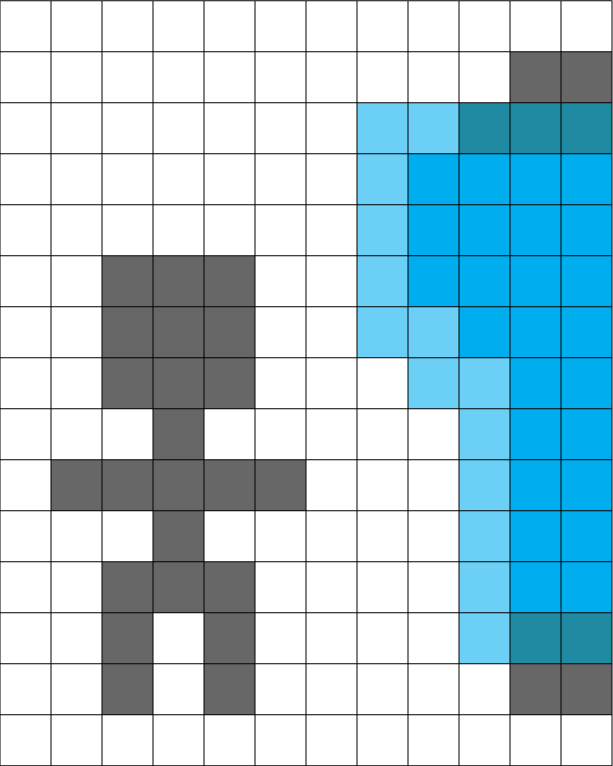}}
        \put(45,260){\includegraphics[width=0.118\textwidth]{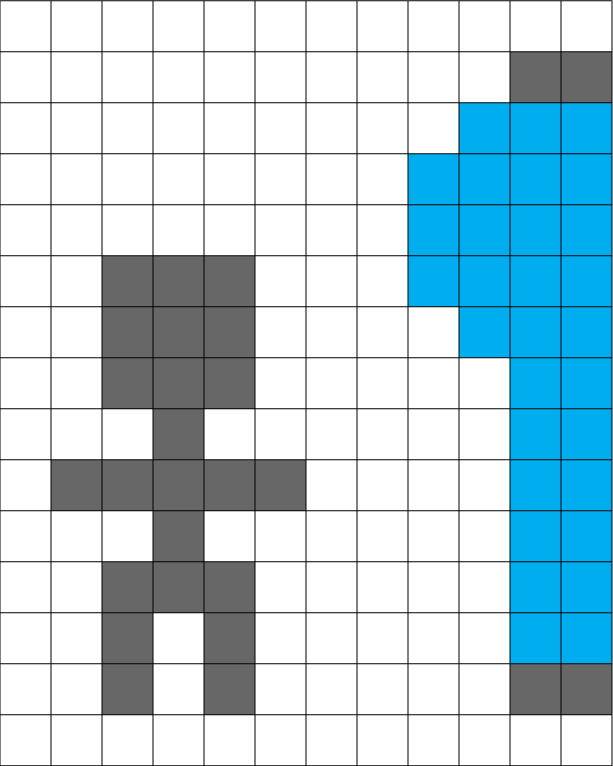}}
        \put(108,260){\includegraphics[width=0.118\textwidth]{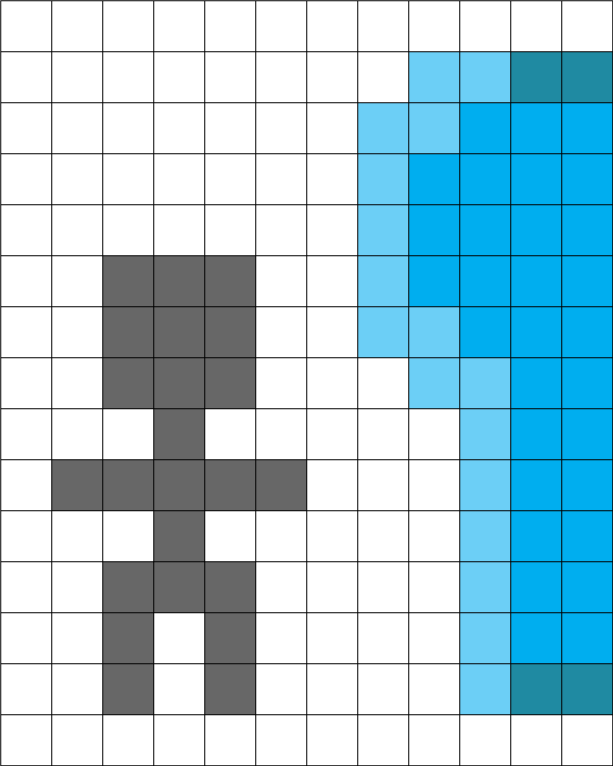}}
        \put(170,260){\includegraphics[width=0.118\textwidth]{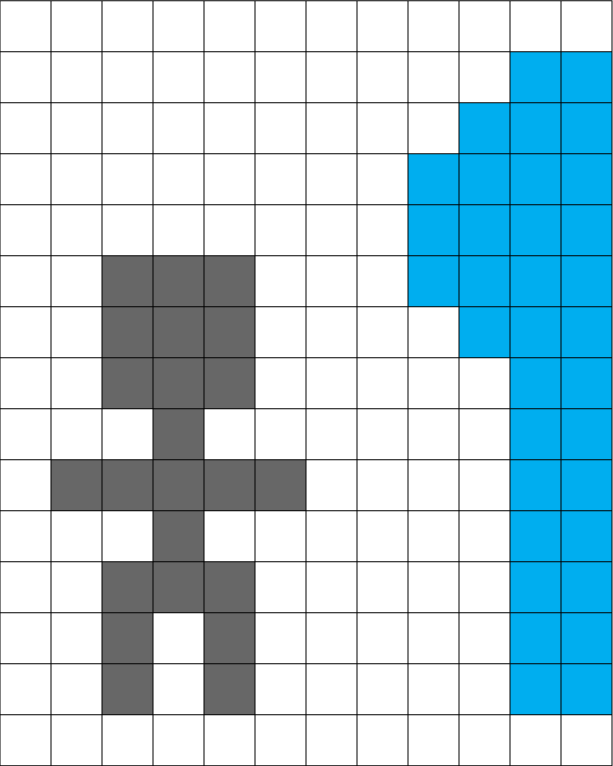}}
        
        \put(-142,173){\includegraphics[width=0.118\textwidth]{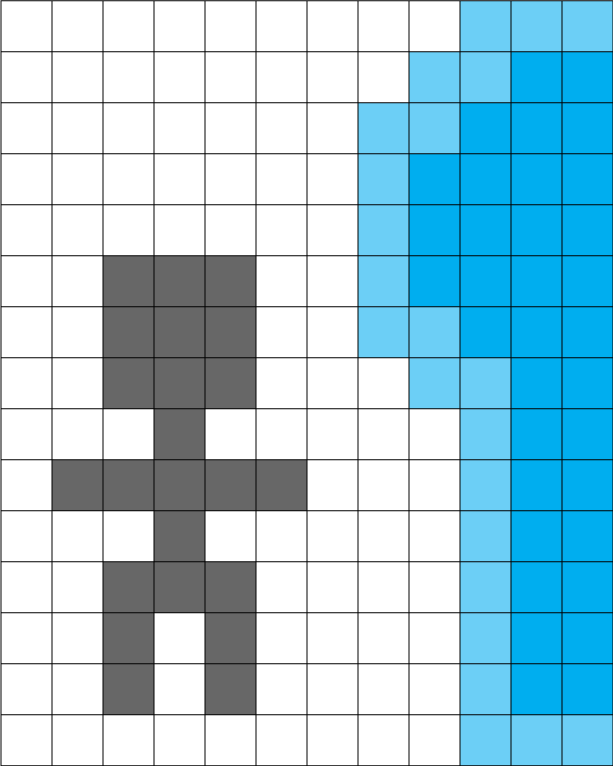}}
        \put(-80,173){\includegraphics[width=0.118\textwidth]{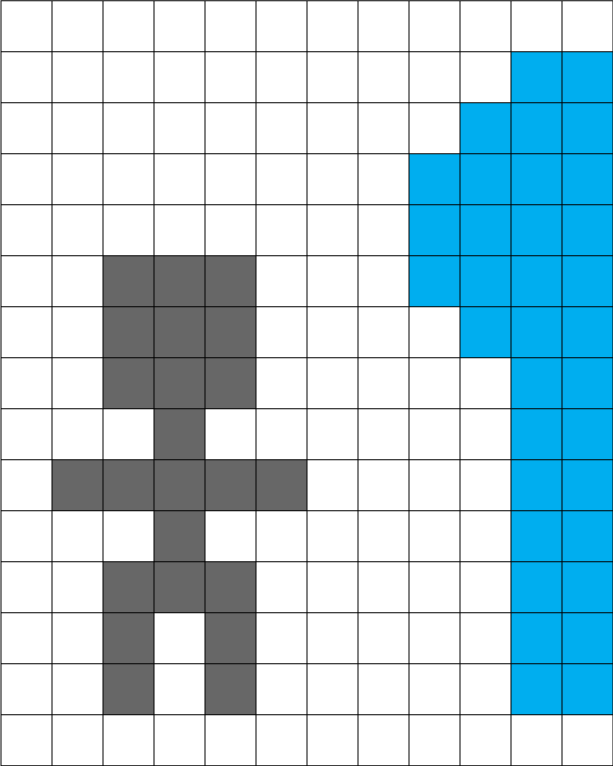}}

        \put(-205,83){\includegraphics[width=0.118\textwidth]{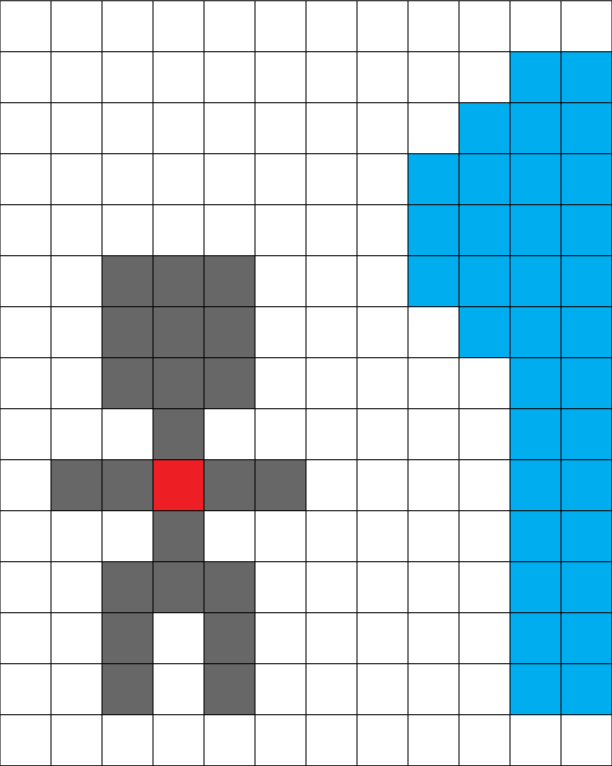}}
        \put(-142,83){\includegraphics[width=0.118\textwidth]{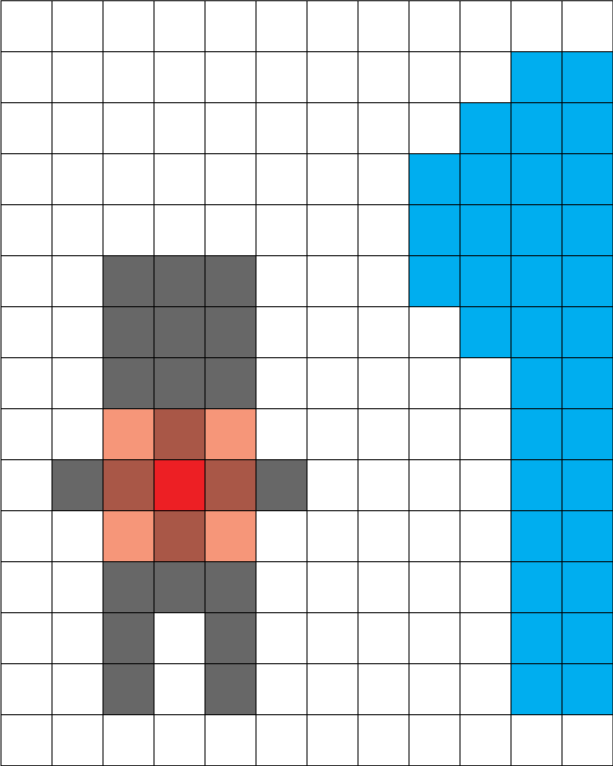}}
        \put(-80,83){\includegraphics[width=0.118\textwidth]{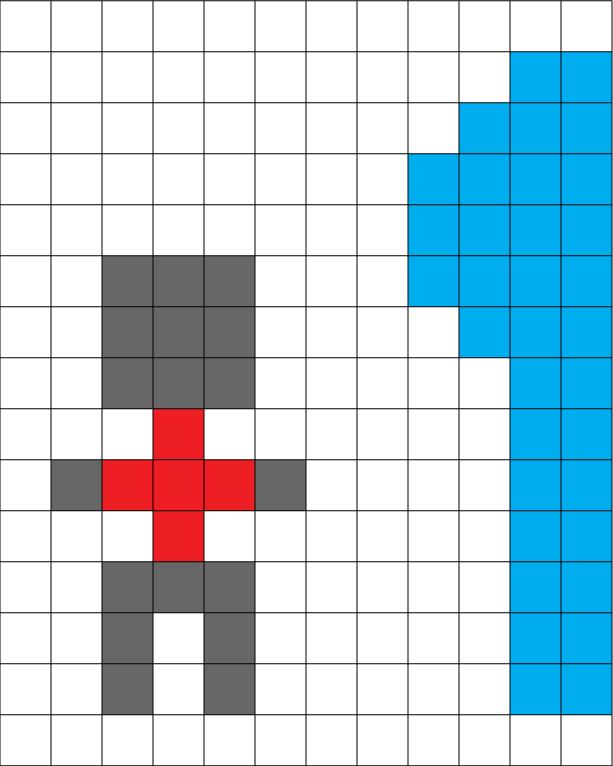}}
        \put(-17,83){\includegraphics[width=0.118\textwidth]{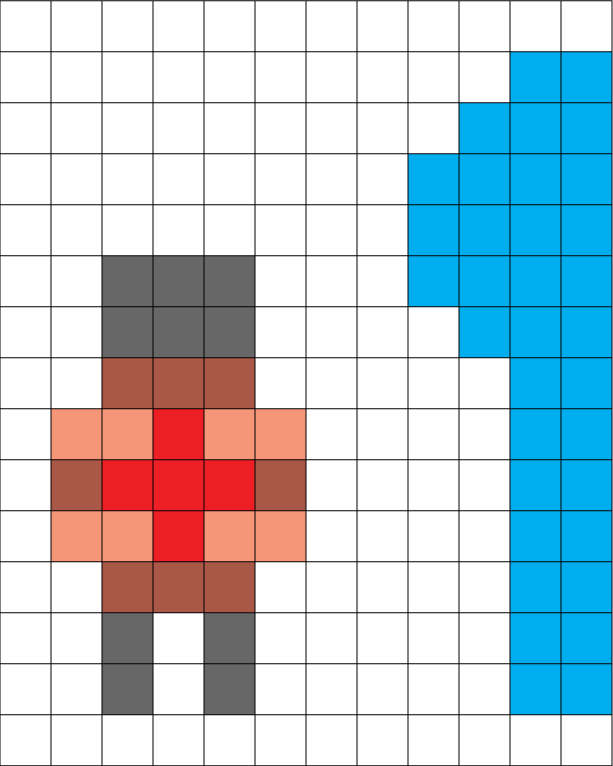}}
        \put(45,83){\includegraphics[width=0.118\textwidth]{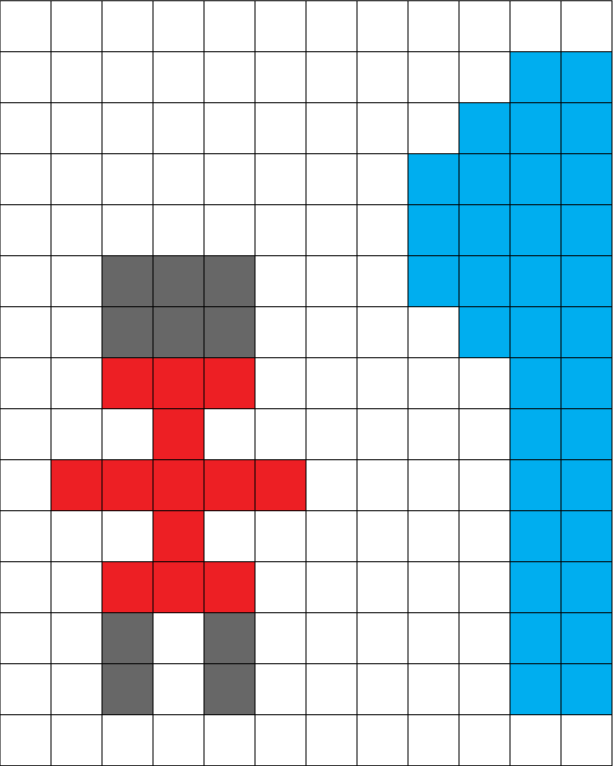}}
        \put(108,83){\includegraphics[width=0.118\textwidth]{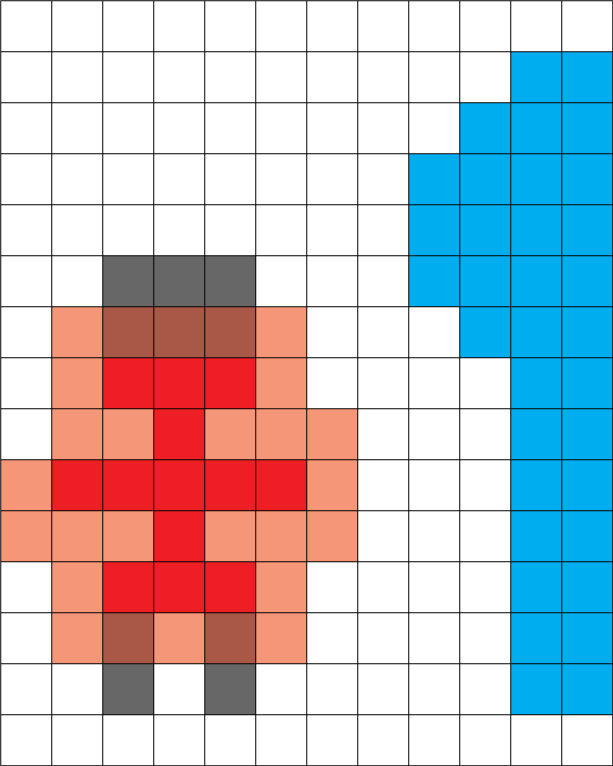}}
        \put(170,83){\includegraphics[width=0.118\textwidth]{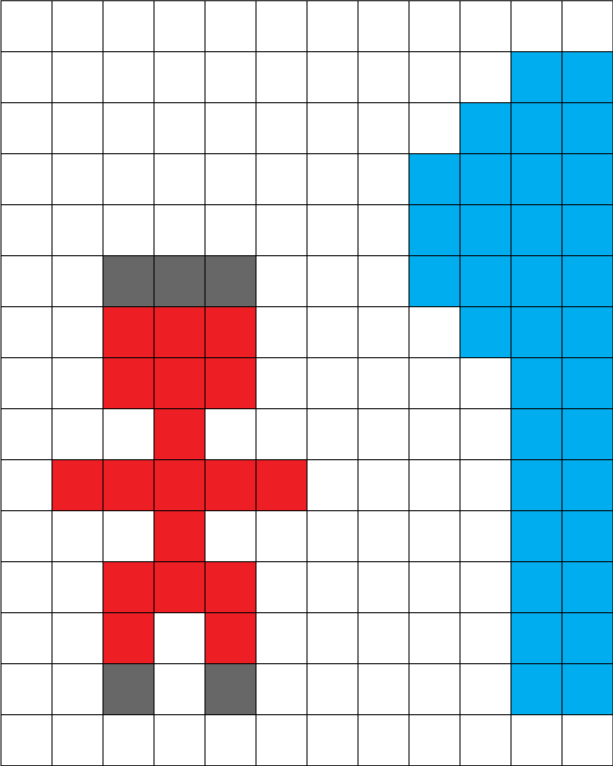}}
        
        \put(-142,-4){\includegraphics[width=0.118\textwidth]{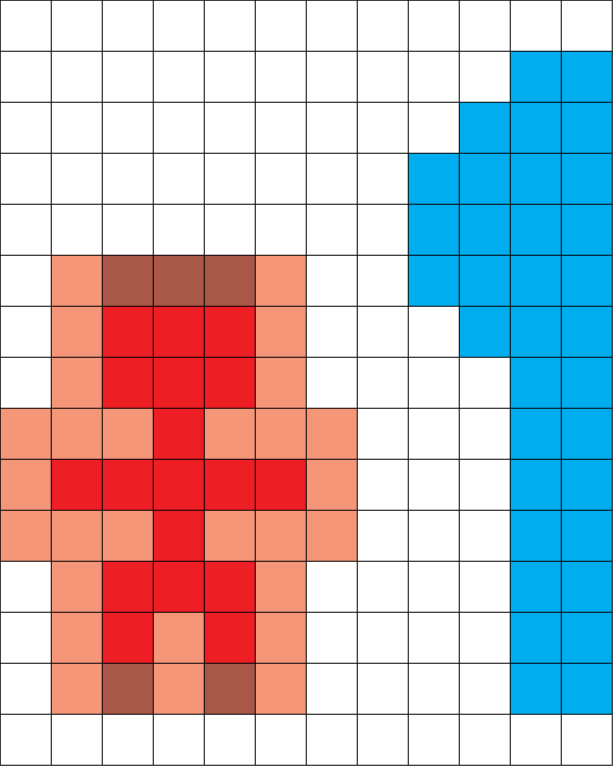}}
        \put(-80,-4){\includegraphics[width=0.118\textwidth]{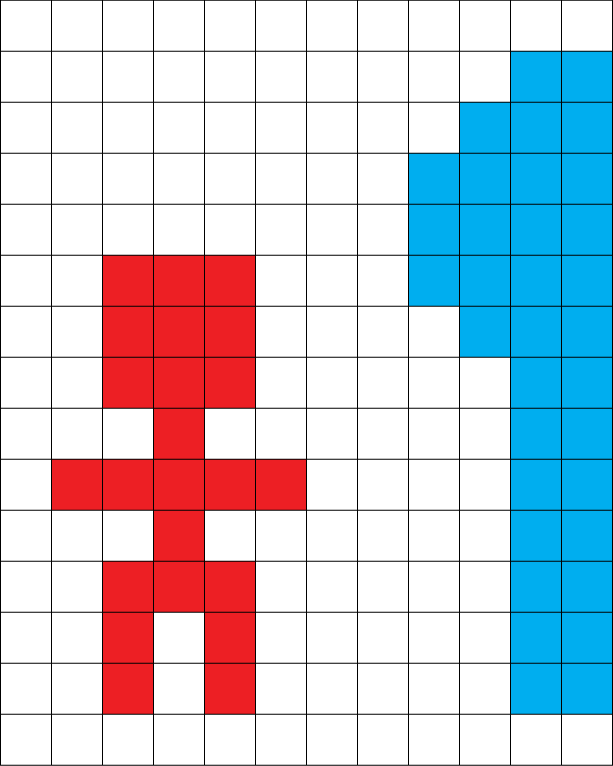}}
        \put(-17,-4){\includegraphics[width=0.118\textwidth]{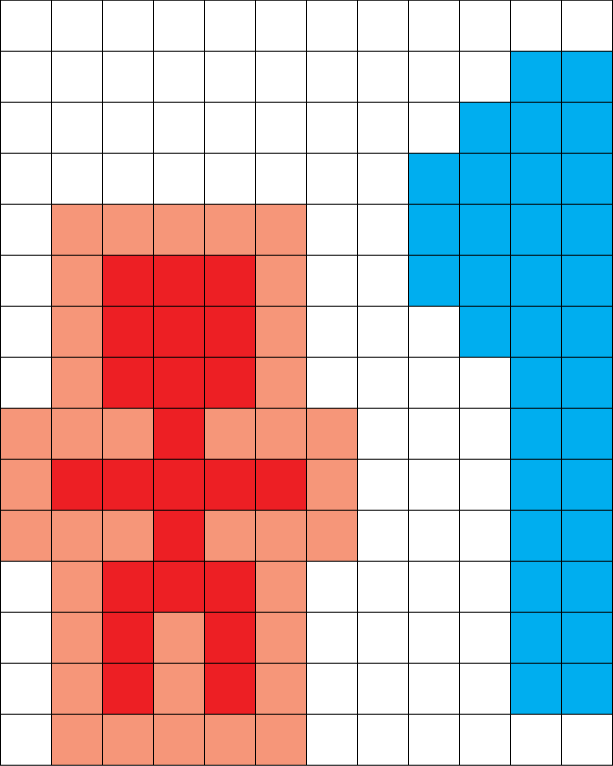}}
        \put(45,-4){\includegraphics[width=0.118\textwidth]{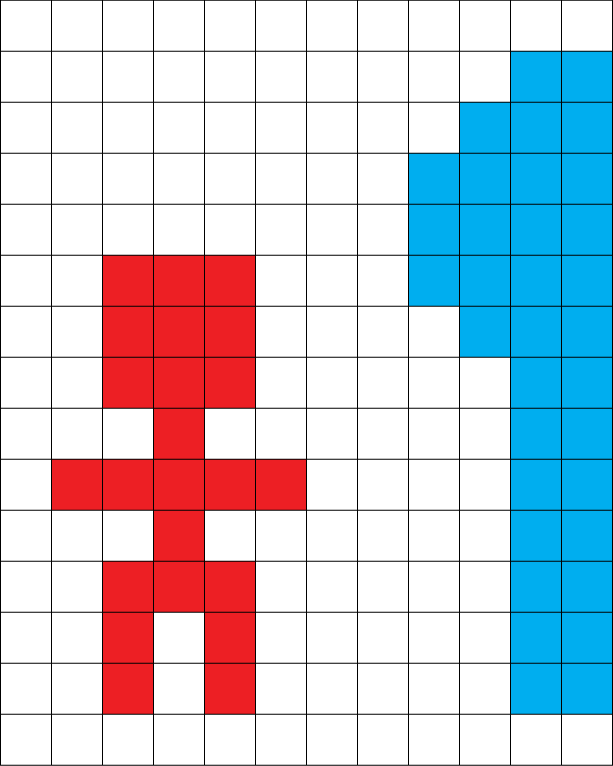}}
        
        \put(-166,380){\textcolor{black}{$g_1$}}
        \put(-140,400){\textcolor{black}{$\R{g_1}{0}$}}
        \put(-115,398){\vector(4,-3){15}}
        \put(-78,400){\textcolor{black}{$\C{g_1}{1}$}}
        \put(-53,398){\vector(4,-3){15}}
        \put(-48,200){\textcolor{black}{$[g_1]$}}
        
        \put(-205,103){\textcolor{black}{$g_2$}}
        \put(70,10){\textcolor{black}{$[g_2]$}}
        
        \put(-230,386){\textcolor{black}{$k=1$}}
        \put(-230,126){\textcolor{black}{$k=2$}}
        
        \put(-187,426){\textcolor{black}{$n=0$}}
        \put(-205,422){\makebox(58,1){\downbracefill}}
        \put(-93,426){\textcolor{black}{$n=1$}}
        \put(-142,422){\makebox(120,1){\downbracefill}}
        \put(33,426){\textcolor{black}{$n=2$}}
        \put(-17,422){\makebox(120,1){\downbracefill}}
        \put(159,426){\textcolor{black}{$n=3$}}
        \put(108,422){\makebox(120,1){\downbracefill}}
        \put(-93,339){\textcolor{black}{$n=4$}}
        \put(-142,335){\makebox(120,1){\downbracefill}}
        \put(33,339){\textcolor{black}{$n=5$}}
        \put(-17,335){\makebox(120,1){\downbracefill}}
        \put(159,339){\textcolor{black}{$n=6$}}
        \put(108,335){\makebox(120,1){\downbracefill}}
        \put(-93,252){\textcolor{black}{$n=7$}}
        \put(-142,248){\makebox(120,1){\downbracefill}}

        \multiput(-235,170)(20,0){24}{\line(5,0){5}}
        
        \put(-187,162){\textcolor{black}{$n=0$}}
        \put(-205,158){\makebox(58,1){\downbracefill}}
        \put(-93,162){\textcolor{black}{$n=1$}}
        \put(-142,158){\makebox(120,1){\downbracefill}}
        \put(33,162){\textcolor{black}{$n=2$}}
        \put(-16,158){\makebox(120,1){\downbracefill}}
        \put(159,162){\textcolor{black}{$n=3$}}
        \put(110,158){\makebox(120,1){\downbracefill}}
        \put(-93,75){\textcolor{black}{$n=4$}}
        \put(-142,71){\makebox(120,1){\downbracefill}}
        \put(33,75){\textcolor{black}{$n=5$}}
        \put(-16,71){\makebox(120,1){\downbracefill}}
    \end{picture}
    \caption{A visualization of a connected component algorithm on a set $\Omega=\{1,\ldots,12\}\times\{1,\ldots,15\}$. Here, one sees in the top figure $I=\{\text{grey pixel locations}\}$. There are two separate connected components which are identified in two steps: $k=1$ shows the process of identifying $[g_1]\coloneqq\{\text{blue pixel locations}\}$ and $k=2$ shows the process of finding the second connected component $[g_2]\coloneqq\{\text{red pixel locations}\}$.
    }
    \label{fig:20230822:VisualizationCCAlgorithm}
\end{figure*}
\subsection{How the Algorithm Works in an Example}
An example of finding connected components in (a discrete subset of) $\mathbb{R}^2$ is given in Figure~\ref{fig:20230822:VisualizationCCAlgorithm}. Here, we consider an image of a person and a tree. The set $I$ consists of the center points of the pixels of width $\rho$. Then, we execute the algorithm with $\delta=\rho\sqrt{2}$.

The first \deltaConnected\ component is identified in row 2 to 4 of Fig.~\ref{fig:20230822:VisualizationCCAlgorithm}. We initiate the algorithm by selecting $g_1=\C{g_1}{0}\in I$, depicted by the blue pixel in the first column. Then, the set $\R{g_1}{0}$ is identified in transparent blue and contains all pixels whose center point has at most distance $\delta$ to the point $g_1$. Then, set $\C{g_1}{1}$ is shown in non-transparent blue in the 2nd row, 3rd column. These steps are repeated until stabilization is reached, so until $\C{g_1}{n}=\C{g_1}{n-1}$.
The final two rows depict the identification of the second connected component in red.

%% file: Figures/VisualizationICompact.tex
\begin{tikzpicture}
    \filldraw[color=red!40, fill=red!20, very thick](-1,0) circle (1.5);
    \filldraw[color=red!40, fill=red!40, very thick](-1,0) circle (0.05) node[anchor=north] {$h_1$};
    \filldraw[color=red!40, fill=red!20, very thick](1,0) circle (1.5);
    \filldraw[color=red!40, fill=red!40, very thick](1,0) circle (0.05) node[anchor=north] {$h_2$};
    \filldraw[color=red!40, fill=red!20, very thick](3,0) circle (1.5);
    \filldraw[color=red!40, fill=red!40, very thick](3,0) circle (0.05) node[anchor=north] {$h_3$};
    \filldraw[color=red!40, fill=red!20, very thick](4.45,2) circle (1.5);
    \filldraw[color=red!40, fill=red!40, very thick](4.45,2) circle (0.05) node[anchor=north] {$h_4$};
    \filldraw[color=red!40, fill=red!20, very thick](0,6) circle (1.5);
    \filldraw[color=red!40, fill=red!40, very thick](0,6) circle (0.05) node[anchor=north] {$g_1$};
    \filldraw[color=red!40, fill=red!20, very thick](2,7) circle (1.5);
    \filldraw[color=red!40, fill=red!40, very thick](2,7) circle (0.05) node[anchor=north] {$g_2$};
    \filldraw[color=red!40, fill=red!20, very thick](4,8) circle (1.5);
    \filldraw[color=red!40, fill=red!40, very thick](4,8) circle (0.05) node[anchor=north] {$g_3$};
    \filldraw[color=red!40, fill=red!20, very thick](8.5,2.25) circle (1.5);
    \filldraw[color=red!40, fill=red!40, very thick](8.5,2.25) circle (0.05) node[anchor=north] {$f_1$};
    \filldraw[color=red!40, fill=red!20, very thick](7.75,4) circle (1.5);
    \filldraw[color=red!40, fill=red!40, very thick](7.75,4) circle (0.05) node[anchor=north] {$f_2$};
    \filldraw[color=red!40, fill=red!20, very thick](7.25,5.75) circle (1.5);
    \filldraw[color=red!40, fill=red!40, very thick](7.25,5.75) circle (0.05) node[anchor=north] {$f_3$};
    \draw[color=red, very thick](6,9) node[anchor=north] {$[g]$};
    \draw[color=red, very thick](5,0) node[anchor=north] {$[h]$};
    \draw[color=red, very thick](3.8,0) node[anchor=north] {$[h_3]$};
    \draw[color=red, very thick](5.5,2) node[anchor=north] {$[h_4]$};
    \draw[color=red, very thick](10,5) node[anchor=north] {$[f]$};
    \filldraw[blue!100, very thick, opacity = 0.2] (-2,-1) rectangle (3.5,1);
    \filldraw[blue!100, very thick, opacity = 0.2] (3.7,1) rectangle (5.2,3);
    \filldraw[blue!100, very thick, opacity = 0.2, rotate around={27:(2,7)}] (-1,6) rectangle (5,8);
    \filldraw[blue!100, very thick, opacity = 0.2, rotate around={110:(5,6)}] (-0.5,3) rectangle (4.5,5);
    \draw[black!100, very thick](-.2,1)--(-.2,4.75) (-.2,3.5)node[ right] {$D_\delta([\{g\}],[\{h\}])$};
    \draw[black!100, very thick](0.5,3.0)node[right] {$=\Delta\!\left([g],[h]\right)>\delta$};
    \draw[black!100, very thick](3.5,1)--(3.7,1) (3.5,0.7)node[anchor = west] {$D_\delta([\{h_3\}],[\{h_4\}])$};
    \draw[black!100, very thick](4.4,0.2)node[anchor = west] {$=\Delta\!\left([h_3],[h_4]\right)<\delta$};
    \draw[black!100, very thick](-1,0)--({-1.5*sqrt(1)/2-1},{1.5*sqrt(3)/2}) node[midway, left] {$\delta$};
\end{tikzpicture}

%% file: Figures/VisualizationDeltaConnectedPoints.tex
\begin{tikzpicture}
    \filldraw[color=red!40, fill=red!20, very thick](-1,0) circle (1.5);
    \filldraw[color=red!40, fill=red!40, very thick](-1,0) circle (0.05) node[anchor=north] {$h_1$};
    \filldraw[color=red!40, fill=red!20, very thick](1,0) circle (1.5);
    \filldraw[color=red!40, fill=red!40, very thick](1,0) circle (0.05) node[anchor=north] {$h_2$};
    \filldraw[color=red!40, fill=red!20, very thick](3,0) circle (1.5);
    \filldraw[color=red!40, fill=red!40, very thick](3,0) circle (0.05) node[anchor=north] {$h_3$};
    \filldraw[color=red!40, fill=red!20, very thick](4.45,2) circle (1.5);
    \filldraw[color=red!40, fill=red!40, very thick](4.45,2) circle (0.05) node[anchor=north] {$h_4$};
    \filldraw[color=red!40, fill=red!20, very thick](0,6) circle (1.5);
    \filldraw[color=red!40, fill=red!40, very thick](0,6) circle (0.05) node[anchor=north] {$g_1$};
    \filldraw[color=red!40, fill=red!20, very thick](2,7) circle (1.5);
    \filldraw[color=red!40, fill=red!40, very thick](2,7) circle (0.05) node[anchor=north] {$g_2$};
    \filldraw[color=red!40, fill=red!20, very thick](4,8) circle (1.5);
    \filldraw[color=red!40, fill=red!40, very thick](4,8) circle (0.05) node[anchor=north] {$g_3$};
    \filldraw[color=red!40, fill=red!20, very thick](8.5,2.25) circle (1.5);
    \filldraw[color=red!40, fill=red!40, very thick](8.5,2.25) circle (0.05) node[anchor=north] {$f_1$};
    \filldraw[color=red!40, fill=red!20, very thick](7.75,4) circle (1.5);
    \filldraw[color=red!40, fill=red!40, very thick](7.75,4) circle (0.05) node[anchor=north] {$f_2$};
    \filldraw[color=red!40, fill=red!20, very thick](7.25,5.75) circle (1.5);
    \filldraw[color=red!40, fill=red!40, very thick](7.25,5.75) circle (0.05) node[anchor=north] {$f_3$};
    \draw[color=red, very thick](6,9) node[anchor=north] {$[g]$};
    \draw[color=red, very thick](5,0) node[anchor=north] {$[h]$};
    \draw[color=red, very thick](10,5) node[anchor=north] {$[f]$};
    \filldraw[black!100, very thick, opacity = 0.2] (-2,-1) rectangle (3.5,1);
    \filldraw[black!100, very thick, opacity = 0.2] (3.7,1) rectangle (5.2,3);
    \filldraw[purple!100, very thick, opacity = 0.2, rotate around={27:(2,7)}] (-1,6) rectangle (5,8);
    \filldraw[blue!100, very thick, opacity = 0.2, rotate around={110:(5,6)}] (-0.5,3) rectangle (4.5,5);
\end{tikzpicture}

%% file: Review_JMIV/3MorphologicalDilations.tex
\section{Computing a Thickened Set \texorpdfstring{$A_\epsilon$}{Ae} with Hamilton-Jacobi Equations}\label{sec:MorphologicalDilations}
The above algorithm hinges on computing $\delta$-thickened versions of the sets $\C{g}{n}$ in the sense of Def.~\ref{def:20230922:epsilonThickenedSet}. In this section, we give a practical algorithm based on morphological dilations on the Lie group $G$. We carry the discussion for a generic set $A$ first. The transformations used to create a thickened set $A_\epsilon$ rely on morphological convolutions as defined next.

\begin{definition}[Morphological convolution]\ \\
Let $G$ be a Lie group. Let $f_1,f_2: G\to\mathbb{R}$ be lower semi-continuous functions, bounded from below. Then the morphological convolution $(f_1\square f_2):G\to\mathbb{R}$ is given by
\begin{align*}
    (f_1\square f_2)(g)=\inf_{h\in G} \left\{f_1(h^{-1}g)+f_2(h)\right\}.
\end{align*}
\end{definition}

To create an $\epsilon$-thickened set $A_\epsilon$, we rely on morphological dilations. 
The concept of morphological dilations is connected to the notion of viscosity solutions of Hamilton-Jacobi-Bellmann (HJB) equations. For a given $\alpha>1$ and initial condition $f\in C(G)
$, we search for a $W:G\times [0, T]\to \mathbb{R}$ that satisfies the Hamilton-Jacobi-Bellmann equation
\begin{align}
\begin{cases}
    \frac{\partial W}{\partial t}(g,t)=\frac{1}{\alpha}\left\|\nabla_\mathcal{G} W(g,t)\right\|^\alpha,& t\in (0, T]\\
    W(g,0)=f(g),&g\in G.
\end{cases}\label{eq:20230922:DilationPDE}
\end{align}
The above equation has a unique viscosity solution $W$ which is given by
\begin{align}
    W(g,t)=-(\morphkern{t}{\alpha}\square -f)(g),\label{eq:viscositySolutionW}
\end{align}
where
\begin{align}
\morphkern{t}{\alpha}(g)&\coloneqq \frac{t}{\beta}\left(\frac{d(g,e)}{t}\right)^\beta\eqqcolon\kappa_t^\alpha(d(g,e))\label{eq:20230822:MorphologicalKernelNew}
\end{align}
is called a morphological kernel with $\frac{1}{\alpha}+\frac{1}{\beta}=1,\;\alpha,\beta>1$ and where $\kappa_t^\alpha(x)=\frac{t}{\beta}(x/t)^\beta$, $x\in\mathbb{R}$ denotes the 1D morphological kernel. Equation \eqref{eq:20230822:MorphologicalKernelNew} shows that the morphology on $G$ relates to the 1D morphology on the distance map. For more details on this relation, cf. Appendix~\ref{sec:multmorphdil}.

\begin{definition}[Morphological dilation] \\
We say that $W(\cdot, t)$ is a morphological dilation of $f$ with kernel $\morphkern{t}{\alpha}$ given by \eqref{eq:viscositySolutionW} if $W$ is the viscosity solution of Eq.~\eqref{eq:20230922:DilationPDE}.
\end{definition}

For subsequent developments, it will be useful to view the viscosity solution $W(\cdot, t)$ given in \eqref{eq:viscositySolutionW} as the image of $f$ under the flow map $\varphi_t^\alpha:C(G)\to\mathbb{R}$ given by
\begin{align}
\varphi_t^\alpha(f)\coloneqq W(\cdot,t)=-(\morphkern{t}{\alpha}\square -f).\label{eq:20230922:FlowReprMorphDilationWithKernel}
\end{align}
The semigroup property of the flow implies that
\begin{align}
    \varphi_t^\alpha\circ\varphi_s^\alpha(f)=\varphi_{t+s}^\alpha(f), \quad \forall t, s\in \mathbb{R},\; \forall f\in C(G).\label{eq:semigroupProperty}
\end{align}
\begin{remark}\label{rem:FlowProperty}
    Note that the semigroup property of the flow implies that $\varphi_t^\alpha\circ^{(n-1)}\varphi_t^\alpha f=\varphi_{nt}^\alpha f$ for all $t>0$. This is due to the well-posedness of \eqref{eq:20230922:DilationPDE} in terms of viscosity solutions, but can also be seen more explicitly by looking at the corresponding kernels, see Lemma \ref{lemma:20230605CleanedUpVersion:KernelConcatenation} in Appendix \ref{sec:multmorphdil}.
\end{remark}

\begin{remark}[Relation to sub-Riemannian diffusions]\ \\
    Often, researchers rely on sub-Riemannian diffusion to determine the grouping of nearby objects and points \cite{abbasi2016geometric,PhDSanguinetti,duits2011leftInvariant,BekkersNilpotent}. Linear left-invariant diffusions are solved by a convolution with a heat kernel $h_t^\alpha$: $(h_{t}^{\alpha}* U)(g)$. It is expensive to compute this heat kernel $h_t^\alpha$ exactly, but if the Lie group $G$ is of polynomial growth, one can find an upper and lower bound for the heat kernel with Maheux' heat kernel bounds \cite[Lemma 6.6]{PhDSmets}, i.e. then there exist constants $c_1,c_2>0$ and for every $\epsilon>0$ there exists a $c_\epsilon>0$ so that
    \begin{align*}
        c_1 \eta_t e^{-\frac{d_{\cG}(g,e)^2}{4c_2t}}\leq h_t^1(g)\leq c_\epsilon \eta_t e^{-\frac{d_{\cG}(g,e)^2}{4(1+\epsilon)t}},
    \end{align*}
    with normalization constant $\eta_t$ given by $\eta_t\coloneqq \mu_\cG(B(e,\sqrt{t}))^{-1}$, where $\mu_{\cG}$ is the volume measure induced by $\cG$. Similarly, one can squeeze the approximated kernel by the real one, i.e. if $G$ is of polynomial growth, there exist constants $C\geq 1$, $D_1\in(0,1)$, $D_2>D_1$ so that for all $t>0$
    \begin{align*}
        \frac{1}{C} h_{D_1 t}^1(g)\leq h_t^{1,\text{approx}}(g)\leq C K_{D_2 t}^1(g),
    \end{align*}
    see \cite[Sec.5.2 Lemma 24]{Smets2023PDEBased}. However, dilations and erosions are \emph{exactly} solved by $(\morphkern{t}{\alpha}\square U)(g)$, where $\morphkern{t}{\alpha}$ is introduced in Eq.~\eqref{eq:20230822:MorphologicalKernelNew}, resulting in more accurate grouping results. Akin to optimal transport on $\SE{2}$ \cite{bonPai2025optimal}, one could rely for \deltaConnected\ components on Varadhan's Theorem, i.e. $d_{\cG}^2(g,e)=-\lim_{t\to 0}(4t)\log \morphkern{t}{\alpha}(g,e)$ \cite{Vadamaran1967behavior}. However, this would be less accurate \cite[Fig.5.1\&5.2]{bonPai2025optimal} than applying HJB-solvers with morphological convolutions, also when using logarithmic norm approximations.
\end{remark}

We next list some useful properties connected to the semigroup property of the HJB equation, where we use the indicator function and the support:
\begin{definition}[Indicator function]\label{def:IndicatorFunction}\\
    The indicator function $\mathbbm{1}_S:G\to\{0,1\}$ of some set $S\subset G$ is defined as
    \begin{align*}
        \mathbbm{1}_S(x)\coloneqq\begin{cases}
            1 &\text{if }x\in S,\\
            0 &\text{else}.
        \end{cases}
    \end{align*}
\end{definition}
\begin{definition}[Support of a function]\\
    The support of the function $f:G\to\mathbb{R}$ is given by
    \begin{align*}
        \text{supp}(f)\coloneqq\overline{\left\{g\in G\;|\; f(g)\neq 0\right\}}.
    \end{align*}
\end{definition}

\begin{lemma}[Properties morphological dilations]\label{lemma:propertiesMorphologicalDilations}
Let $\alpha\geq 1$. Let $\frac{1}{\alpha}+\frac{1}{\beta}=1$. For a given closed set $A\subset G$ and a function $f:G\to[0,1]$, one has
\begin{align}
&\varphi_t^\alpha(f) (g) \in [0, 1], && \forall g\in G,\label{eq:morphConvRange}\\
&A_{\varepsilon(t, \alpha)} = \text{{\normalfont supp}}\left( \varphi_t^\alpha( \mathbbm{1}_{A}) \right), && \forall t\geq 0,\qquad\qquad\label{eq:suppThickenedSetByFlow}\\
&\mathbbm{1}_{A_{t}} = \left( \varphi_t^1( \mathbbm{1}_{A}) \right), \label{eq:flow-indicative}
\end{align}
with
\begin{align}
    \varepsilon(t, \alpha) = t\sqrt[\beta]{\frac{\beta}{t}},\label{eq:SizeThickenedSetMorphConv}
\end{align}
and where $A_t$, $A_{\varepsilon(t,\alpha)}$ are $t$- and $\varepsilon(t,\alpha)$-thickened sets of A, as defined in Def.~\ref{def:20230922:epsilonThickenedSet}.
\end{lemma}

\begin{proof}
    We prove the statements one by one, starting with \eqref{eq:morphConvRange}:
    Let $f:G\to[0,1]$. Then, as the kernels are positive \eqref{eq:20230822:MorphologicalKernelNew}, we have for all $g\in G$
    \begin{align*}
        \varphi_t^\alpha(f)(g)\overset{\eqref{eq:20230922:FlowReprMorphDilationWithKernel}}{=}-\left(\morphkern{t}{\alpha}\square -f\right)(g)= \sup_{h\in G}\left\{f(h)-\morphkern{t}{\alpha}(h^{-1}g)\right\}\\\leq \|f\|_{\infty}=1,\quad
    \end{align*}
    and 
    \begin{align*}
        \varphi_t^\alpha(f)(g)\overset{\eqref{eq:20230922:FlowReprMorphDilationWithKernel}}{=}-\left(\morphkern{t}{\alpha}\square -f\right)(g)&= \sup_{h\in G}\left\{f(h)-\morphkern{t}{\alpha}(h^{-1}g)\right\}\\&\geq \sup_{h\in G}\left\{-\morphkern{t}{\alpha}(h^{-1}g)\right\}=0.
    \end{align*}
    
    We next prove \eqref{eq:suppThickenedSetByFlow}:
    We are interested in the support of $\varphi_t^\alpha(\mathbbm{1}_A)$. From \eqref{eq:morphConvRange}, we know
    \begin{align*}
        \varphi_t^\alpha(\mathbbm{1}_A)=-(\morphkern{t}{\alpha}\square-\mathbbm{1}_A)\in [0,1].
    \end{align*}
    By definition of the support and the positivity of $\varphi_t^\alpha(\mathbbm{1}_A)$, we have $\text{supp}(\varphi_t^\alpha(\mathbbm{1}_A))=\overline{\{g\in G\;|\; \varphi_t^\alpha(\mathbbm{1}_A)(g)>0\}}$, so we aim to find the set $\tilde{A}$, such that for all $g\in\tilde{A}$, we have
    \begin{align}
         0<\varphi_t^\alpha(\mathbbm{1}_A)(g)
         &=\sup_{h\in G}\left\{\mathbbm{1}_A(h)-\frac{t}{\beta}\left(\frac{d(h^{-1}g,e)}{t}\right)^\beta\right\}\leq 1,\label{eq:ProofThickenedSetA}
    \end{align}
    Note that $\overline{\tilde{A}}=\supp(\varphi_t^\alpha(\mathbbm{1}_A))$ and that the upper boundary is always satisfied as the kernel $\morphkern{t}{\alpha}$, recall Eq.~\eqref{eq:20230822:MorphologicalKernelNew}, is positive and $\mathbbm{1}_A\leq 1$. The supremum is only strictly larger than zero if the supremum is reached for $h\in A$, so $A\subset\overline{\tilde{A}}$. Then, \eqref{eq:ProofThickenedSetA} gives
    \begin{align*}
        0<1-\frac{t}{\beta}\left(\frac{d(g,A)}{t}\right)^\beta\quad\Leftrightarrow\quad d(g,A)<t\sqrt[\beta]{\frac{\beta}{t}}.
    \end{align*}
    Therefore, the set $\overline{\tilde{A}}=A_{\varepsilon(t,\alpha)}$, where $A_{\varepsilon(t,\alpha)}$ is a thickened set of $A$ as defined in Def.~\ref{def:20230922:epsilonThickenedSet}. This proves Eq.~\eqref{eq:suppThickenedSetByFlow}.

    The last statement \eqref{eq:flow-indicative} follows immediately from \eqref{eq:suppThickenedSetByFlow} by taking the limit of $\alpha\downarrow 1$, combined with 
    \begin{align*}
        \morphkern{t}{1}(g)&=\lim\limits_{\alpha\downarrow 1}\morphkern{t}{\alpha}(g)=\lim\limits_{\beta\to\infty}\frac{t}{\beta}\left(\frac{d(g,e)}{t}\right)^\beta\\
        &=\begin{cases}
            0& \text{if }d(g,e)\leq t,\\
            \infty& \text{else}.
        \end{cases}\numberthis\label{eq:morphKernelAlpha1}
    \end{align*}\qed
\end{proof}
We conclude this section by connecting one last time to the main goal in this part of the article which is to explain how to create an $\epsilon$-thickened set $A_\epsilon$. The main message is that $A_\epsilon$ can be obtained through Eqs.~\eqref{eq:suppThickenedSetByFlow} and $\eqref{eq:flow-indicative}$, which crucially rely on morphological dilations.
In practice, it is useful to rely on Eq.~\eqref{eq:flow-indicative} to determine the $\epsilon$-thickened set $A_\epsilon$ in the \deltaConnected\ component algorithm, as we will see in the next section.

%% file: Review_JMIV/4_1ConnectedComponentAlgorithm.tex
\section{Computing \texorpdfstring{$\delta$-Connected}{delta-Connected} Components using Iterative Morphological Convolutions}\label{sec:ComputingCCWithMorphConv}
In Sec.~\ref{sec:GeneralAlgorithm}, we explained the \deltaConnected\ component algorithm in a general setting. In Sec.~\ref{sec:20230822:Algorithm}, we will explain how the \deltaConnected\ component algorithm is constructed using morphological dilations, but first, we discuss the choice of the parameter $\alpha$ in the morphological dilation kernel $\morphkern{t}{\alpha}$ in Sec.~\ref{sec:choiceAlpha}. We finish with a convergence analysis for the algorithm in Sec.~\ref{sec:20230822:Analysis}.

\subsection{Fixing $\alpha=1$ in the $\delta$-Connected Component Algorithm}\label{sec:choiceAlpha}
We aim to construct a \deltaConnected\ component algorithm that relies on morphological dilations. Before we introduce the algorithm, we discuss the choice of the parameter $\alpha\geq 1$ in the morphological dilation kernel $\morphkern{t}{\alpha}$ \eqref{eq:20230822:MorphologicalKernelNew}.

We recall that identifying the \deltaConnected\ components is an iterative procedure as introduced in Sec.~\ref{sec:GeneralAlgorithm}. The identification of one single \deltaConnected\ component is mainly described by Eq.~\eqref{eq:20230822:defC_g^n+1}. In every iteration step $n$, we create a $\delta$-thickened set $C_{\delta}(g,n-1)$ and take the intersection with a reference set $I\subset G$, resulting in an updated set $C(g,n)$. This is repeated until the updated set is the same as the initial set, i.e. $C(g,n)=C(g,n-1)$.

Hence, the set $C(g,n)$ is thickened iteratively, which motivates the identification of the maximum radius of influence of one single point after $n$ morphological dilations, as will be introduced in Def.~\ref{def:MaxRofInfluence}. The maximum radius of influence follows immediately from the semigroup property in Remark~\ref{rem:FlowProperty} and the support after applying the flow to an indicator function in Equation~\eqref{eq:suppThickenedSetByFlow}.
\begin{definition}[Maximum radius of influence]\label{def:MaxRofInfluence}\ \\
    The maximum radius of the region of influence of a specific point after $n$ applications of a morphological dilation of time $t$ and kernel steepness $\alpha$ is denoted by 
    \begin{align*}
        \maxRofInfluence{n}{t}{\alpha}\coloneqq nt\left(\frac{\beta}{nt}\right)^{\frac{1}{\beta}}=\varepsilon(nt,\alpha),
    \end{align*}
    where $\frac{1}{\alpha}+\frac{1}{\beta}=1$ and $\varepsilon(t,\alpha)$ was defined in Eq.~\eqref{eq:SizeThickenedSetMorphConv}.
\end{definition}

We recall that an $\epsilon$-thickened set was produced with Eq.~\eqref{eq:suppThickenedSetByFlow} and \eqref{eq:flow-indicative}. We aim to identify \deltaConnected\ components. Therefore, the support of the set has to grow with rate $\delta$ in every iteration step. For computational efficiency, we would like to use the same kernel $\morphkern{t}{\alpha}$ in every iteration step of the algorithm. That means we would like to have for $A_{2\varepsilon(t,\alpha)}=\left(A_{\varepsilon(t,\alpha)}\right)_{\varepsilon(t,\alpha)}$ that
\begin{align*}
    A_{2\varepsilon(t,\alpha)}\!\!\!\overset{\eqref{eq:suppThickenedSetByFlow}}{=}\!\!\!\text{supp}(\varphi_t^\alpha(\varphi_t^\alpha(\mathbbm{1}_A)))\!\!\!\overset{\eqref{eq:semigroupProperty}}{=}\!\!\!\text{supp}(\varphi_{2t}^\alpha(\mathbbm{1}_A))\!\!\!\overset{\eqref{eq:suppThickenedSetByFlow}}{=}\!\!\!A_{\varepsilon(2t,\alpha)}.
\end{align*}
This is only satisfied if $\beta\to\infty\Leftrightarrow\alpha=1$ as is clear from Eq.~\eqref{eq:SizeThickenedSetMorphConv}.

If $\alpha>1$, the algorithm cannot always identify the full \deltaConnected\ component, while using the same kernel $\morphkern{t}{\alpha}$ in every set, as is visualized in Fig.~\ref{fig:20230922:Symmetry} in App.~\ref{app:symmetry}.
Therefore, in our computations of the \deltaConnected\ components, we set $\alpha=1$ and apply morphological convolutions with $\morphkern{t}{1}$, given by \eqref{eq:morphKernelAlpha1}, with $\delta\coloneqq t$.

\subsection{$\delta$-Connected Component Algorithm}\label{sec:20230822:Algorithm}
Next, we detail how we implement the general \deltaConnected\ component algorithm as discussed in Sec.~\ref{sec:GeneralAlgorithm} using the morphological dilations that were discussed in Sec.~\ref{sec:MorphologicalDilations}. More specifically, we explain how the $\delta$-thickened set $C_\delta$ in \ref{alg:findFullComponent} is determined using the properties of morphological dilations as stated in Lemma~\ref{lemma:propertiesMorphologicalDilations}.

First, we discuss how to identify one single \deltaConnected\ component $[g_i]$ for a $g_i\in I$. 
We initialize
\begin{align*}
    \C{g_i}{0}=\{g_i\}=\text{supp}(\mathbbm{1}_{\{g_i\}}).
\end{align*}
During iteration step $n$ of identifying one single \deltaConnected\ component $[g_i]$, the set $\C{g_i}{n-1}$ from the previous iteration is expanded with a radius $\delta$ using morphological dilations and Eq.~\eqref{eq:flow-indicative} in Lemma~\ref{lemma:propertiesMorphologicalDilations}.

Therefore, we define a function $U:G\times \mathbb{N}\to\{0,1\}$ which is the indicator function of the computed \deltaConnected\ component initialized in $g_i\in I$ at step $n$, i.e.
\begin{align*}
    \U{g_i}{\cdot}{n-1}\coloneqq\mathbbm{1}_{\C{g_i}{n-1}}(\cdot).
\end{align*}
This allows us to calculate the set $\C{g_i}{n}$ using morphological dilations on the function $\U{g_i}{\cdot}{n-1}$
\begin{align}
    \U{g_i}{\cdot}{n}\coloneqq \mathbbm{1}_I(\cdot)\; \varphi_\delta^1(\U{g_i}{\cdot}{n-1})\label{eq:updateUg0},
\end{align}
where
\begin{align*}
    \C{g_i}{n}&\coloneqq\text{supp}(\U{g_i}{\cdot}{n}).\numberthis\label{eq:relationCandU}
\end{align*}

These steps are repeated until $\C{g_i}{n}=\C{g_i}{n-1}$. Then, we update $I=I\,\backslash\,\C{g_i}{n}$ and pick a new $g_{i+1}\in I$ to identify the next component as explained in \ref{alg:findAllComponents} in Sec.~\ref{sec:GeneralAlgorithm}. The complete algorithm is also stated in Algorithm~\ref{alg:20230605:connectedComponents}: \ref{alg:deltaconnectedComponents} and produces \deltaConnected\ components as stated in the next proposition.

\begin{algorithm}
\caption{\texttt{find\_$\delta$CC}}
\customlabelnoprint{alg:deltaconnectedComponents}{\texttt{find\_$\delta$CC}}
\textbf{Input:} Binary map $\B:G\to\left\{0,1\right\}$, where $I\subset G$ the region of interest for the connected components, kernel $\morphkern{\delta}{\alpha}$ with $\alpha=1$, $\delta>0$\\
set $k=0$\\
set $\mathcal{I}=I$\\
while $\mathcal{I}\neq\emptyset$:\\
$\quad$ update $k=k+1$\\
$\quad$ pick $g_k\in \mathcal{I}$\\
$\quad$ initialize $\U{g_k}{\cdot}{0}=\mathbbm{1}_{\{g_k\}}$\\
$\quad$ initialize connected component sets $\C{g_k}{-1}\coloneqq\emptyset,\; \C{g_k}{0}=\text{supp}\left(\U{g_k}{\cdot}{0}\right)=\{g_k\}$.\\
$\quad$ set $n\coloneqq 0$\\
$\quad$ while $n$ such that $\C{g_k}{n}\neq \C{g_k}{n-1}$:\\
$\qquad$ do a morphological convolution calculating $\U{g_k}{\cdot}{n+1}$ as described in \eqref{eq:updateUg0}.\\
$\qquad$ create set $\C{g_k}{n+1}=\text{supp}\left(\U{g_k}{\cdot}{n+1}\right)$ (see Eq.~\eqref{eq:relationCandU})\\
$\qquad$ update $n\coloneqq n+1$\\
$\quad$ end\\
$\quad$ set $[g_k]\coloneqq \C{g_k}{n}$\\
$\quad$ update $\mathcal{I}=\mathcal{I}\,\backslash\,[g_k]$\\
end\\
$K=k$\\
\textbf{Output:} Connected components $[g_1],\ldots,[g_K]$.
\label{alg:20230605:connectedComponents}
\end{algorithm}

\begin{proposition}\label{prop:AlgCCProducesTconnectedComponents}
    Algorithm~\ref{alg:20230605:connectedComponents}: \ref{alg:deltaconnectedComponents} produces \deltaConnected\ components, i.e. $[g_1],\ldots,[g_K]$ are all \deltaConnected.
\end{proposition}
The result follows immediately from Eq.~\eqref{eq:updateUg0}, using Eq.~\eqref{eq:flow-indicative} in Lemma~\ref{lemma:propertiesMorphologicalDilations} where $t=\delta$.

Prop.~\ref{prop:AlgCCProducesTconnectedComponents} shows that \ref{alg:deltaconnectedComponents} finds \deltaConnected\ components. It remains to be shown that these are all \deltaConnected\ components and that the algorithm is indeed a stable algorithm which we will address next.

\subsection{Convergence Analysis}\label{sec:20230822:Analysis}
In this section, we prove that Algorithm~\ref{alg:20230605:connectedComponents}: \ref{alg:deltaconnectedComponents} computes the \deltaConnected\ components of any compact set $I$ on any Lie group $G$ in a finite number of calculation steps for $\alpha=1$ in Theorem~\ref{th:20230605:ConvergenceConnectedComponentAlgorithm}.

For this, we need to have an explicit non-recursive expression of $\U{g_0}{\cdot}{n}$ because this describes the state of the \deltaConnected\ component algorithm in the $n$-th iteration step. In Proposition~\ref{th:20230605:Un}, we find a non-recursive expression for $\U{g_i}{\cdot}{n}$, where
Definition~\ref{def:MaxRofInfluence} and Lemma~\ref{lemma:propertiesMorphologicalDilations}
are key ingredients in proving the convergence of our \deltaConnected\ component algorithm (Theorem~\ref{th:20230605:ConvergenceConnectedComponentAlgorithm}).

\begin{proposition}\label{th:20230605:Un}
    Let $\alpha = 1$, $t=\delta>0$, $g\in G$, $n\in\mathbb{N}$. Let $g_0\in G$ be the reference point. Set $\delta=\maxRofInfluence{1}{t}{1}$. Then, we have, for $\U{g_0}{g}{n}$ given by Eq. \eqref{eq:updateUg0}, that $0\leq \U{g_0}{g}{n}\leq 1$ where
    \begin{align*}
        \U{g_0}{g}{n}\!=\!\!\left.\begin{cases}
            1 &\text{if }g\simdelta g_0 \wedge m_\delta(g,g_0)\leq n\\
            0&\text{else}.
        \end{cases}\right\}\!,
    \end{align*}
    where $m_\delta(g,g_0)$ was defined in Def.~\ref{def:20230605:Connectedness} and Remark~\ref{remark:MdependsOngAndh}.
\end{proposition}
For the details of the proof of Proposition~\ref{th:20230605:Un}, see Appendix~\ref{app:20230922:ProofPropUn}. Intuitively, the condition $m_\delta(g,g_0)\leq n$ ensures that the $n$-steps of the connected component algorithm are sufficient to reach the point $g$ starting from $g_0$.

Knowing the algorithm's state after $n$ iterations, allows us to confidently say that the algorithm can identify all \deltaConnected\ components in a finite number of steps, as proved in the next theorem.
\begin{theorem}[Convergence of 
Alg.~\ref{alg:20230605:connectedComponents}:~ \ref{alg:deltaconnectedComponents}]\label{th:20230605:ConvergenceConnectedComponentAlgorithm}
%
    Let $\delta>0$, $\alpha=1$, and assume that $I$ has $K$ connected components $[g_1],\dots, [g_K]$ such that $I=[g_1]\cup\cdots\cup[g_K]\subset G$ and $[g_i]\cap[g_j]=\emptyset$ for all $i\neq j\in\{1,\ldots,K\}$.
    Then Algorithm~\ref{alg:20230605:connectedComponents}: \ref{alg:deltaconnectedComponents} correctly finds all the \deltaConnected\ components in a finite number of steps.
\end{theorem}
\begin{proof}
    Consider $\alpha=1$. Let $i\in\left\{1,\ldots,K\right\}$ be given. Note that in this case, one has the dilated volume $\U{g_i}{\cdot}{n}: G\to\left\{0,1\right\}$ and the connected components
    \begin{align*}
        \C{g_i}{n}=\text{supp}\left(\U{g_i}{\cdot}{n}\right)=\left\{g\in G\;\left| \; \U{g_i}{g}{n}\neq 0\right.\right\}.
    \end{align*}
    In the proof of Prop.~\ref{th:20230605:Un}, we defined the set $\tildeB{n}{\delta}{g_i}$ by
    \begin{align}
        \tildeB{n}{\delta}{g_i}\coloneqq \left\{g\in I\;\left|\;m_\delta(g_i,g)\leq n\right.\right\}.\label{eq:setTildeB}
    \end{align}
    Then, by Proposition~\ref{th:20230605:Un} where we set $g_i$ instead of $g_0$, one has
    \begin{align*}
        &\U{g_i}{g}{n}=\mathbbm{1}_{[g_i]}(g)\cdot \mathbbm{1}_{
        \tildeB{n}{\delta}{g_i}}(g)\\
        &\C{g_i}{n}=\tildeB{n}{\delta}{g_i}.
    \end{align*}
    From the compactness of $I$ and the fact that $I=\bigcup_i [g_i]$, $[g_i]\cap [g_j]=\emptyset$, one can deduce that $[g_i]$ is closed. Therefore, also $[g_i]$ is compact.
    Since $[g_i]$ is compact, its covering number $n_\delta([g_i])$, as introduced in \eqref{eq:20230822:CoveringNumber}, is finite and thus $\U{g_i}{\cdot}{n}$ converges in at most $n=n_\delta([g_i])+1<\infty$ steps (see Corollary~\ref{cor:upperBoundMaxItCCAlg}), i.e. $\U{g_i}{\cdot}{n_\delta([g_i])+1}=\U{g_i}{\cdot}{n_\delta([g_i])+2}$ and thus $\C{g_i}{n_\delta([g_i])+1}=\C{g_i}{n_\delta([g_i])+2}$.
    \qed
\end{proof}
As a direct consequence of Theorem~\ref{th:20230605:ConvergenceConnectedComponentAlgorithm}, we can give an upper bound for the maximum number of iterations that \ref{alg:deltaconnectedComponents} needs to identify all \deltaConnected\ components.
\begin{corollary}[Upper bound 
number of iterations]\label{cor:upperBoundMaxItCCAlg}
    Let $I=[g_1]\cup\ldots\cup [g_K]\subset G$ compact, where the \deltaConnected\ components $[g_i]\cap [g_j]=\emptyset$ for all $i\neq j\in\{1,\ldots,K\}$. The \deltaConnected\ component algorithm given in Algorithm~\ref{alg:20230605:connectedComponents}: \ref{alg:deltaconnectedComponents} finishes in at most $n=n_\delta(I)+K$ steps.
\end{corollary}
\begin{proof}
    We start by noticing that $n_\delta(I)=\sum\limits_{i=1}^K n_\delta([g_i])$ because $d([g_i],[g_j])>\delta$ for $i\neq j$.

    To identify \deltaConnected\ component $[g_i]$, we need at most 1 step to get from the starting point $g_i\in I$ to a point $\tilde{g}_i\in C$, where $C$ is the minimal $\delta$-covering of $I$. Then, by Lemma~\ref{lemma:20230605:maxNrdeltaConnectedComponents}, we need at most $n_\delta([g_i])$ steps to identify the full \deltaConnected\ component $[g_i]$.

    Hence, identifying all \deltaConnected\ components in $I=\bigcup\limits_{i=1}^K [g_i]$ is done in at most $\sum\limits_{i=1}^K 1+n_\delta([g_i])=K+n_\delta(I)$ steps.\qed
\end{proof}

\begin{remark}[Complexity of the Algorithm] \\
    The presented method has a complexity of $\cO(n N_o K_o N_x N_y K_s^2)$ in $\SE{2}$ \cite[Table 3.2]{PhDFranken}, where $n\coloneqq n_\delta+K$ denotes the maximum number of iterations of the algorithm, $N_x$, $N_y$, and $N_o$ the dimensions of the orientation score in the $x$-, $y-$ and $\theta$-direction, and $K_s$, and $K_o$ the spatial dimensions and the number of sampled orientations of the $\SE{2}$-kernel, respectively. Note that in \cite{PhDFranken}, they consider linear convolutions instead of morphological convolutions. However, the complexity directly carries
    over from the linear to the morphological setting by replacing the semifield $(\bR,\cdot,+)$ by the semifield $(\bR,+,\inf)$.
\begin{figure*}[h]
    \centering
    \begin{subfigure}[t]{0.2\textwidth}
         \centering
        \includegraphics[width=\textwidth]{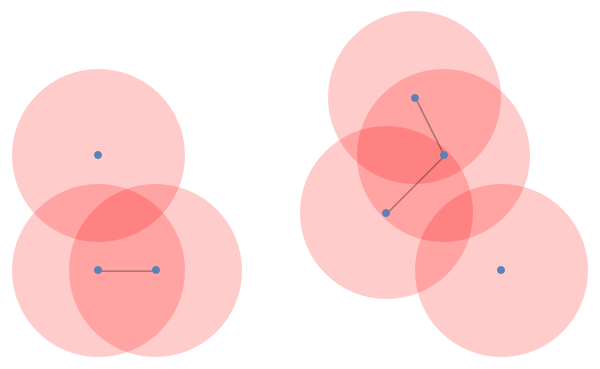}
        \caption{Set $I_0$ including balls of radius $\delta=1.5$.}
        \label{fig:20230605:SetPersistenceCC}
    \end{subfigure}~
    \begin{subfigure}[t]{0.25\textwidth}
         \centering
        \scalebox{0.5}{\input {Figures/PersistenceCC/R2/PersistenceDiagramCC}}
        \caption{Persistence Diagram of $I_0$.}
        \label{fig:20230605:SetPersistenceDiagram}
    \end{subfigure}~
    \begin{subfigure}[t]{0.25\textwidth}
         \centering

        \scalebox{0.5}{\input {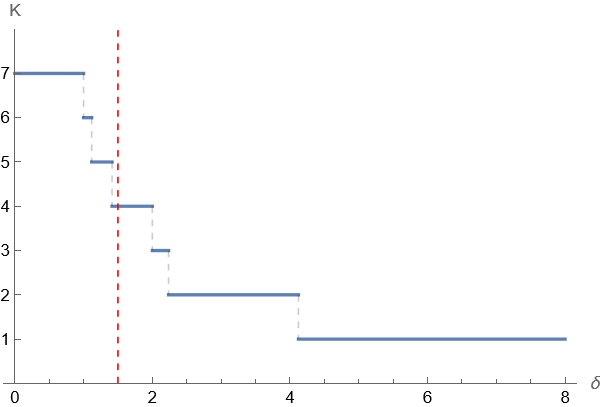}}
        \caption{Amount of \deltaConnected\ components in $I_0$ for various values of $\delta$.}
        \label{fig:20230605:SetPersistenceNCC}
    \end{subfigure}~
    \begin{subfigure}[t]{0.25\textwidth}
         \centering

        \scalebox{0.5}{\input {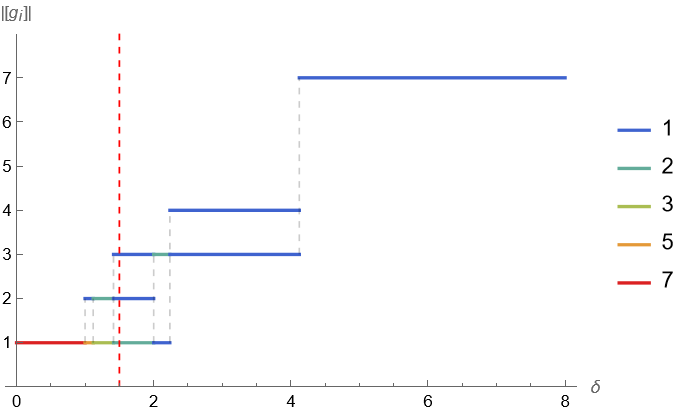}}
        \caption{Number of elements per connected component in $I_0$ for various values of $\delta$.}
        \label{fig:20230605:SetPersistenceNElements}
    \end{subfigure}
    \begin{subfigure}[t]{0.2\textwidth}
         \centering
        \includegraphics[width=\textwidth]{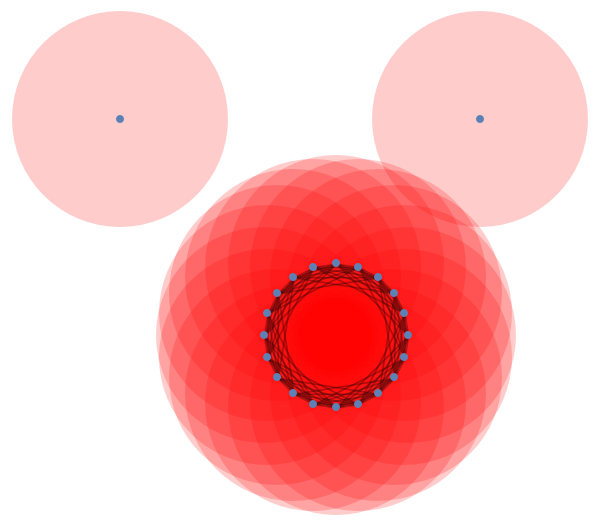}
        \caption{Set $I_1$ including balls of radius $\delta=1.5$.}
        \label{fig:20230605:Set2PersistenceCC}
    \end{subfigure}~
    \begin{subfigure}[t]{0.25\textwidth}
         \centering
        \scalebox{0.5}{\input {Figures/PersistenceCC/R2/PersistenceDiagramCC2}}
        \caption{Persistence Diagram of $I_1$.}
        \label{fig:20230605:Set2PersistenceDiagram}
    \end{subfigure}~
    \begin{subfigure}[t]{0.25\textwidth}
         \centering
        \scalebox{0.5}{\input {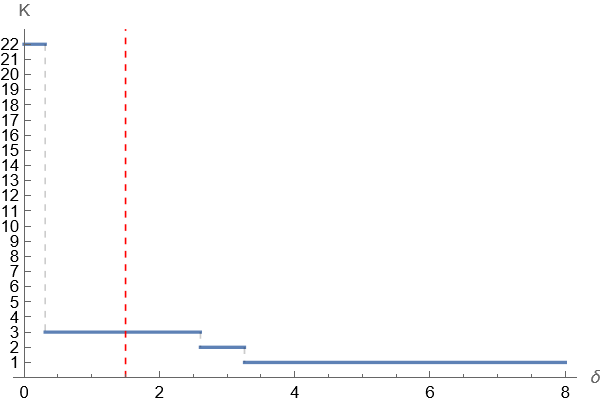}}
        \caption{Amount of \deltaConnected\ components in $I_1$ for various values of $\delta$.}
        \label{fig:20230605:Set2PersistenceNCC}
    \end{subfigure}~
    \begin{subfigure}[t]{0.25\textwidth}
         \centering
        \scalebox{0.5}{\input {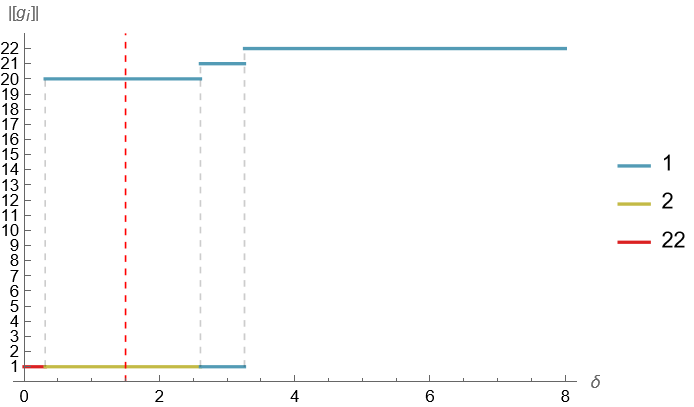}}
        \caption{Number of elements per connected component in $I_1$ for various values of $\delta$.}
        \label{fig:20230605:Set2PersistenceNElements}
    \end{subfigure}
    \caption{Persistency plots for different sets $I\subset\mathbb{R}^2$. The longest lines (most values of $\delta$) denote the most persistent connected components. Depending on the application, the optimal choice of $\delta$ is where the \deltaConnected\ components are the most stable/persistent (i.e. long lines in Fig.~\ref{fig:20230605:SetPersistenceNCC},\ref{fig:20230605:SetPersistenceNElements},\ref{fig:20230605:Set2PersistenceNCC},\ref{fig:20230605:Set2PersistenceNElements}) while aiming for the smallest possible choice of $\delta$ (to be as distinctive as possible).}
    \label{fig:PersistenceConnectedComponent}
\end{figure*}

    In a more general setting, consider a Lie group $G$ of dimension $\dim(G)$. Then, the complexity of the morphological convolution is $\cO(\prod_{i=1}^{\dim(G)} N_i K_i)$, where the image data and the kernel have dimensions \mbox{$N_1\times\ldots\times N_{\dim(G)}$} and $K_1\times\ldots\times K_{\dim(G)}$, respectively. At most $n\coloneqq n_\delta+K$ morphological convolutions are calculated, leading to a complexity of $\cO(n\prod_{i=1}^{\dim(G)} N_i K_i)$.
\end{remark}

%% file: Figures/PersistenceCC/R2/PersistenceDiagramCC.tex

\begin{tikzpicture}
\begin{axis}[
    ytick ={0,1,2,3,4,6},
    yticklabels = {0,1,2,3,4,$\infty$},
    yticklabel style = {align = center, font = \large},
    xtick ={0,1,2,3,4,6},
    xticklabels = {0,1,2,3,4,$\infty$},
    xticklabel style = {align = center, font = \large},
    xlabel style = {font = \Large},
    ylabel style = {font = \Large},
    xlabel = {Birth}, 
    ylabel = {Death}, 
    boxplot/variable width min target = 0,
    ymin = 0, ymax = 6, xmin = 0, xmax =6,
    axis lines*=left,
    legend style={at={(1.1,0.5)},
    anchor=west,legend columns=1, draw = none, cells = {anchor = west}},
    ]

\addplot[scatter, scatter/use mapped color={draw=none,fill=black, fill opacity = 1, draw opacity = 0},only marks]
table[row sep=newline, x = birth, y=death] {Figures/PersistenceCC/R2/PersistenceDataCC.txt};

\addplot[mark = none, black, dashed] coordinates {(0,0) (7,7)};

\end{axis}
\end{tikzpicture}

%% file: Figures/PersistenceCC/R2/TotalNumberOfCC.tex
\pgfplotstableread[]{
birth death N
0 1. 7
1. 1.11803 6
1.11803 1.41421 5
1.41421 2. 4
2. 2.23607 3
2.23607 4.12311 2
4.12311 8 1
}\birthdeath

\begin{tikzpicture}
\begin{axis}[
    xlabel = {$\delta$}, 
    xtick ={0,1,1.5,2,3,4,5,6,7},
    xticklabels ={0,1,\footnotesize\textcolor{red}{$\delta^*$},2,3,4,5,6,7},
    xticklabel style = {align = center, font = \large},
    axis x line=middle,
    xlabel style={font = \Large, anchor=west},
    ylabel = {K}, 
    axis y line=middle,
    ylabel style = {font = \Large, anchor = south},
    ytick ={1,2,3,4,5,6,7},
    yticklabel style = {align = center, font = \large},
    boxplot/variable width min target = 0,
    ymin = 0, ymax = 7.5, xmin = 0, xmax =8,
    axis lines*=left,
    legend style={at={(1.1,0.5)},
    anchor=west,legend columns=1, draw = none, cells = {anchor = west}},
    ]

\addplot[mark = none, blue] coordinates {(0,7) (1,7)};
\addplot[mark = none, blue] coordinates {(1,6) (1.11803,6)};
\addplot[mark = none, blue] coordinates {(1.11803,5) (1.41421,5)};
\addplot[mark = none, blue] coordinates {(1.41421,4) (2,4)};
\addplot[mark = none, blue] coordinates {(2,3) (2.23607,3)};
\addplot[mark = none, blue] coordinates {(2.23607,2) (4.12311,2)};
\addplot[mark = none, blue] coordinates {(4.12311,1) (8,1)};

\addplot[mark = none, black, dashed] coordinates {(1,7) (1,6)};
\addplot[mark = none, black, dashed] coordinates {(1.11803,6) (1.11803,5)};
\addplot[mark = none, black, dashed] coordinates {(1.41421,5) (1.41421,4)};
\addplot[mark = none, black, dashed] coordinates {(2,4) (2,3)};
\addplot[mark = none, black, dashed] coordinates {(2.23607,3) (2.23607,2)};
\addplot[mark = none, black, dashed] coordinates {(4.12311,2) (4.12311,1)};

\addplot[mark = none, red, dashed] coordinates {(1.5,0) (1.5,8)};
    
\end{axis}
\end{tikzpicture}

%% file: Figures/PersistenceCC/R2/ElementsPerCC_colors.tex
\pgfplotstableread[]{
birth death N
0 1. 7
1. 1.11803 6
1.11803 1.41421 5
1.41421 2. 4
2. 2.23607 3
2.23607 4.12311 2
4.12311 8 1
}\birthdeath

\begin{tikzpicture}
\begin{axis}[
    xlabel = {$\delta$}, 
    xtick ={0,1,1.5,2,3,4,5,6,7},
    xticklabels ={0,1,\footnotesize\textcolor{red}{$\delta^*$},2,3,4,5,6,7},
    xticklabel style = {align = center, font = \large},
    axis x line=middle,
    xlabel style={anchor=west,font = \Large},
    ylabel = {size $|[g_i]|$}, 
    axis y line=middle,
    ylabel style = {anchor = south, font = \Large},
    ytick ={1,2,3,4,5,6,7},
    yticklabel style = {align = center, font = \large},
    boxplot/variable width min target = 0,
    ymin = 0, ymax = 7.5, xmin = 0, xmax =8,
    axis lines*=left,
    legend style={at={(1.1,0.5)},
    anchor=west,legend columns=1, draw = none, cells = {anchor = west}},
    ]

\addplot[mark = none, blue] coordinates {(2,1) (2.23607,1)};
\addlegendentry{\Large 1}
\addplot[mark = none, teal] coordinates {(1.41421,1) (2,1)};
\addlegendentry{\Large 2}
\addplot[mark = none, green] coordinates {(1.11803,1) (1.41421,1)};
\addlegendentry{\Large 3}
\addplot[mark = none, orange] coordinates {(1,1) (1.11803,1)};
\addlegendentry{\Large 5}
\addplot[mark = none, magenta] coordinates {(0,1) (1,1)};
\addlegendentry{\Large 7}

\addplot[mark = none, blue] coordinates {(1,2) (1.11803,2)};
\addplot[mark = none, teal] coordinates {(1.11803,2) (1.41421,2)};
\addplot[mark = none, blue] coordinates {(1.41421,2) (2,2)};

\addplot[mark = none, blue] coordinates {(1.41421,3) (2,3)};
\addplot[mark = none, teal] coordinates {(2,3) (2.23607,3)};
\addplot[mark = none, blue] coordinates {(2.23607,3) (4.12311,3)};

\addplot[mark = none, blue] coordinates {(2.23607,4) (4.12311,4)};

\addplot[mark = none, blue] coordinates {(4.12311,7) (8,7)};

\addplot[mark = none, gray, dashed] coordinates {(1,1) (1,2)};
\addplot[mark = none, gray, dashed] coordinates {(1.11803,1) (1.11803,2)};
\addplot[mark = none, gray, dashed] coordinates {(1.41421,1) (1.41421,3)};
\addplot[mark = none, gray, dashed] coordinates {(2,1) (2,3)};
\addplot[mark = none, gray, dashed] coordinates {(2.23607,1) (2.23607,4)};
\addplot[mark = none, gray, dashed] coordinates {(4.12311,3) (4.12311,7)};

\addplot[mark = none, red, dashed] coordinates {(1.5,0) (1.5,8)};
    
\end{axis}
\end{tikzpicture}

%% file: Figures/PersistenceCC/R2/PersistenceDiagramCC2.tex

\begin{tikzpicture}
\begin{axis}[
    ytick ={0,1,2,3,4,6},
    yticklabels = {0,1,2,3,4,$\infty$},
    yticklabel style = {align = center, font = \large},
    xtick ={0,1,2,3,4,6},
    xticklabels = {0,1,2,3,4,$\infty$},
    xticklabel style = {align = center, font = \large},
    xlabel style = {font = \Large},
    ylabel style = {font = \Large},
    xlabel = {Birth}, 
    ylabel = {Death}, 
    boxplot/variable width min target = 0,
    ymin = 0, ymax = 6, xmin = 0, xmax =6,
    axis lines*=left,
    legend style={at={(1.1,0.5)},
    anchor=west,legend columns=1, draw = none, cells = {anchor = west}},
    ]

\addplot[scatter, scatter/use mapped color={draw=none,fill=black, fill opacity = 1, draw opacity = 0},only marks]
table[row sep=newline, x = birth, y=death] {Figures/PersistenceCC/R2/PersistenceDataCC.txt};

\addplot[mark = none, black, dashed] coordinates {(0,0) (7,7)};

\end{axis}
\end{tikzpicture}

%% file: Figures/PersistenceCC/R2/TotalNumberOfCC2.tex
\pgfplotstableread[]{
birth death N
0 0.312869 22
0.312869 2.60668 3
2.60668 3.25871 2
3.25871 8 1
}\birthdeath

\begin{tikzpicture}
\begin{axis}[
    xlabel = {$\delta$}, 
    xtick ={0,1,1.5,2,3,4,5,6,7},
    xticklabels ={0,1,\footnotesize\textcolor{red}{$\delta^*$},2,3,4,5,6,7},
    xticklabel style = {align = center, font = \large},
    axis x line=middle,
    xlabel style = {font = \Large,anchor=west},
    ylabel = {K}, 
    axis y line=middle,
    ylabel style = {font = \Large,anchor = south},
    ytick ={1,2,3,4,5,6,7,8,9,10,11,12,13,14,15,16,17,18,19,20,21,22},
    yticklabels = {1,,3,,,,,,,,,,,,,,,,,,,22},
    yticklabel style = {align = center, font = \large},
    boxplot/variable width min target = 0,
    ymin = 0, ymax = 23, xmin = 0, xmax =8,
    axis lines*=left,
    legend style={at={(1.1,0.5)},
    anchor=west,legend columns=1, draw = none, cells = {anchor = west}},
    ]

\addplot[mark = none, blue] coordinates {(0,22) (0.312869,22)};
\addplot[mark = none, blue] coordinates {(0.312869,3) (2.60668,3)};
\addplot[mark = none, blue] coordinates {(2.60668,2) (3.25871,2)};
\addplot[mark = none, blue] coordinates {(3.25871,1) (8,1)};

\addplot[mark = none, black, dashed] coordinates {(0.312869,22) (0.312869,3)};
\addplot[mark = none, black, dashed] coordinates {(2.60668,3) (2.60668,2)};
\addplot[mark = none, black, dashed] coordinates {(3.25871,2) (3.25871,1)};

\addplot[mark = none, red, dashed] coordinates {(1.5,0) (1.5,23)};
    
\end{axis}
\end{tikzpicture}

%% file: Figures/PersistenceCC/R2/ElementsPerCC2_colors.tex
\pgfplotstableread[]{
birth death N
0 0.312869 22
0.312869 2.60668 3
2.60668 3.25871 2
3.25871 8 1
}\birthdeath

\begin{tikzpicture}
\begin{axis}[
    xlabel = {$\delta$}, 
    xtick ={0,1,1.5,2,3,4,5,6,7},
    xticklabels ={0,1,\footnotesize\textcolor{red}{$\delta^*$},2,3,4,5,6,7},
    xticklabel style = {align = center, font = \large},
    axis x line=middle,
    xlabel style={font = \Large, anchor=west},
    ylabel = {size $|[g_i]|$}, 
    axis y line=middle,
    ylabel style = {font = \Large, anchor = south},
    ytick ={1,2,3,4,5,6,7,8,9,10,11,12,13,14,15,16,17,18,19,20,21,22},
    yticklabels = {1,,,,,,,,,,,,,,,,,,,20,,22},
    yticklabel style = {align = center, font = \large},
    boxplot/variable width min target = 0,
    ymin = 0, ymax = 23.5, xmin = 0, xmax =8,
    axis lines*=left,
    legend style={at={(1.1,0.5)},
    anchor=west,legend columns=1, draw = none, cells = {anchor = west}},
    ]

\addplot[mark = none, blue] coordinates {(2.60668,1) (3.25871,1)};
\addlegendentry{\Large 1}
\addplot[mark = none, green] coordinates {(0.312869,1) (2.60668,1)};
\addlegendentry{\Large 2}
\addplot[mark = none, magenta] coordinates {(0,1) (0.312869,1)};
\addlegendentry{\Large 22}

\addplot[mark = none, blue] coordinates {(0.312869,20) (2.60668,20)};

\addplot[mark = none, blue] coordinates {(2.60668,21) (3.25871,21)};

\addplot[mark = none, blue] coordinates {(3.25871,22) (8,22)};

\addplot[mark = none, gray, dashed] coordinates {(0.312869,1) (0.312869,20)};
\addplot[mark = none, gray, dashed] coordinates {(2.60668,1) (2.60668,21)};
\addplot[mark = none, gray, dashed] coordinates {(3.25871,1) (3.25871,22)};

\addplot[mark = none, red, dashed] coordinates {(1.5,0) (1.5,25)};
    
\end{axis}
\end{tikzpicture}

%% file: Review_JMIV/4_2ChoiceOfDelta.tex
\section{Choosing \texorpdfstring{$\delta$}{delta} with Persistence Diagrams}\label{sec:20230922:Persistence}
In the previous section, we have discussed how to find \deltaConnected\ components using morphological convolutions. We chose a fixed $\delta$ at the start of the algorithm, and use \ref{alg:deltaconnectedComponents} to identify all \deltaConnected\ components.
Consequently, the algorithm's output heavily depends on the choice of $\delta$ (recall Fig.~\ref{fig:20230605:VisualizationConnectedComponents}).

To gain some insight into the behavior of the \deltaConnected\ components, it is common in Topological Data Analysis to create so-called `persistence diagrams' to study the stability of the components \cite{Edelsbrunner2002topological,Zomorodian2005computing,Chazal2013persistenceBased,Skraba2010persistenceBased}. Inspired by this approach, we visualize the persistence of components in three different ways. Firstly, we create the classical persistence diagrams from topological data analysis. Secondly, we create a figure where we plot the total number of \deltaConnected\ components for different values of $\delta$. When the distance between two components is smaller than the considered $\delta$, they merge and hence the number of connected components decreases. Lastly, we visualize, for sets of points, how many points are grouped in each \deltaConnected\ component. These plots (for instance in Fig.~\ref{fig:PersistenceConnectedComponent}) allow us to verify the persistence of the \deltaConnected\ components. Depending on the application, we want them to be as persistent as possible, e.g. when looking at vascular images, or we want to choose $\delta$ small enough to differentiate between the smallest possible \deltaConnected\ components. Note that using a very small value of $\delta$ results in a high number of \deltaConnected\ components due to noise being assigned its own \deltaConnected\ component.

In Fig.~\ref{fig:PersistenceConnectedComponent}, we show the persistence diagrams for two examples in $\mathbb{R}^2$. The first concerns a set of points (7 in total) that belong to our set $I_0$, depicted in Fig.~\ref{fig:20230605:SetPersistenceCC}. Then, we determine the \deltaConnected\ components for a range of values of $\delta$. We plot the number of connected components per value of $\delta$ (cf. Fig.~\ref{fig:20230605:SetPersistenceNCC}) and the number of vertices that are part of every \deltaConnected\ component (cf. Fig.~\ref{fig:20230605:SetPersistenceNElements}). One can see from the plots that the most stable connected components (at least one connection, but not everything is connected) occur when $\delta\in (2.2, 4.1)$.

We do the same for the second example, a different set $I_1$, where 20 points are placed in a circle, and with two outliers (see Fig.~\ref{fig:20230605:Set2PersistenceCC}). Then, the persistence graphs in Fig.~\ref{fig:20230605:Set2PersistenceNCC} and Fig.~\ref{fig:20230605:Set2PersistenceNElements} show a stabilization for $\delta\in(0.3,2.6)$.

\begin{figure}[ht!]
    \centering
    \begin{subfigure}[t]{0.2\textwidth}
         \centering
        \includegraphics[width=0.7\textwidth]{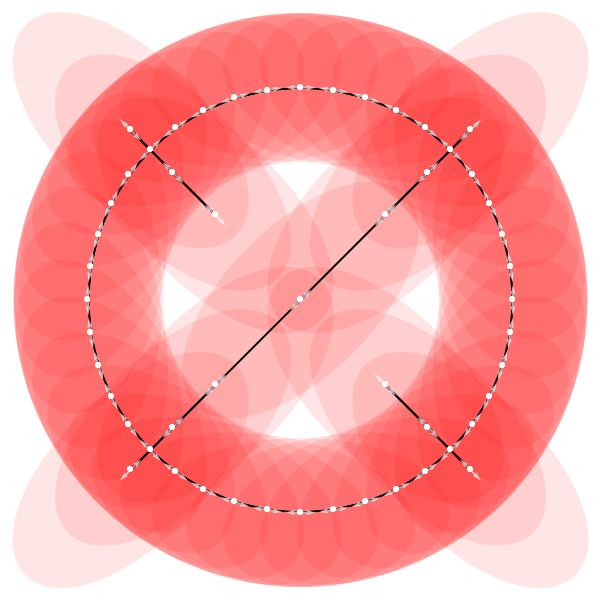}
        \caption{Set $I_0\subset\SE{2}$ including balls of radius $\delta^*=0.7$.}
        \label{fig:20230605:SetSE2PersistenceCC}
    \end{subfigure}~
    \begin{subfigure}[t]{0.2\textwidth}
         \centering
        \includegraphics[width=0.7\textwidth]{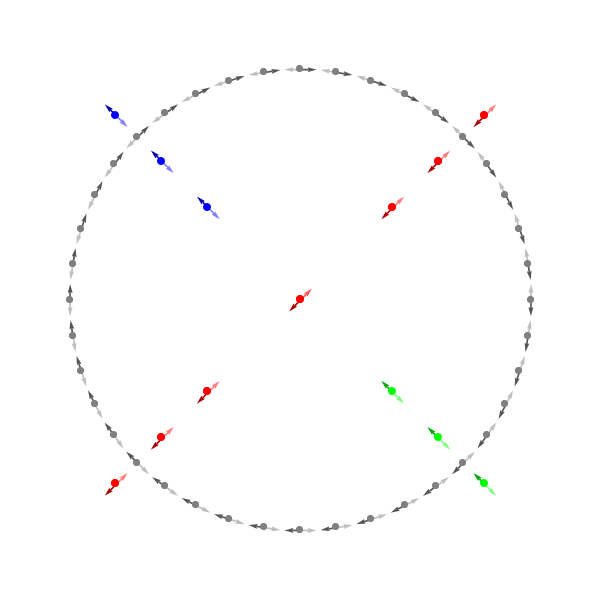}
        \caption{Different \deltaConnected\ components identified with $\delta^*=0.7$.}
        \label{fig:20230605:SetSE2CC}
    \end{subfigure}
    \begin{subfigure}[t]{0.45\textwidth}
         \centering
        \scalebox{0.7}{\input {Figures/PersistenceCC/SE2/PersistenceDiagram}}
        \caption{Persistence Diagram.}
        \label{fig:20230605:SetSE2PersistenceDiagram}
    \end{subfigure}
    \begin{subfigure}[t]{0.45\textwidth}
         \centering
        \scalebox{0.6}{\input {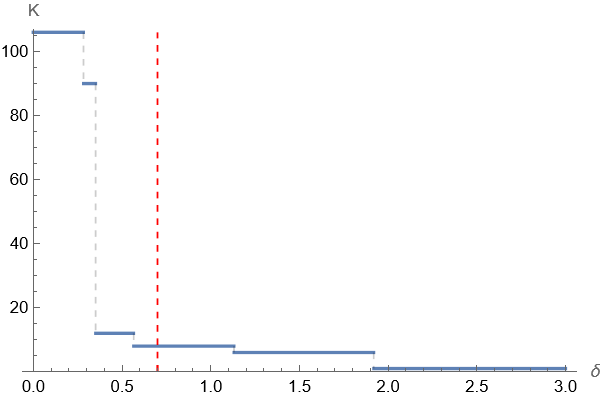}}
        \caption{Amount of \deltaConnected\ components in $I_0$ for various values of $\delta$.}
        \label{fig:20230605:SetSE2PersistenceNCC}
    \end{subfigure}
    \begin{subfigure}[t]{0.45\textwidth}
         \centering
        \scalebox{0.6}{\input {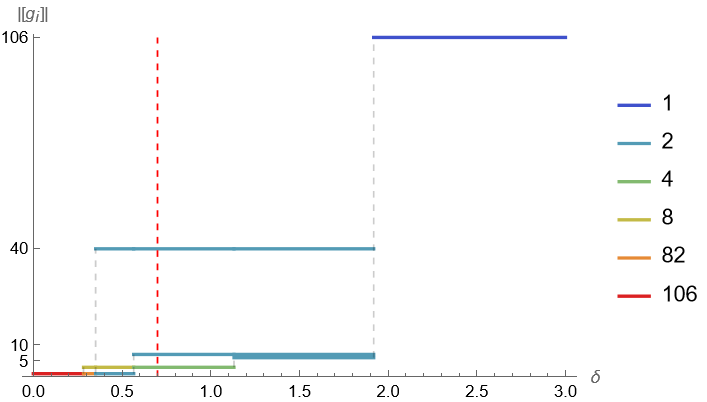}}
        \caption{Number of elements per connected component in $I_0$ for various values of $\delta$.}
        \label{fig:20230605:SetSE2PersistenceNElements}
    \end{subfigure}
    \caption{Persistency plots for different sets $I\subset G=\SE{2}$ using metric tensor field parameters $(w_1,w_2,w_3)=(1,2,2)$. The longest lines (most values of $\delta$) denote the most persistent connected components. Depending on the application, the optimal choice of $\delta$ is where the \deltaConnected\ components are the most stable/persistent (i.e. long lines in Fig~\ref{fig:20230605:SetSE2PersistenceNCC},\ref{fig:20230605:SetSE2PersistenceNElements}) while aiming for the smallest possible choice of $\delta$ (to be as distinctive as possible).}
    \label{fig:PersistenceConnectedComponentonG}
\end{figure}

In Fig.~\ref{fig:PersistenceConnectedComponentonG}, we do the same, but now on the Lie group $G=SE(2)$. Again, we see persistence of the \deltaConnected\ components for $\delta\in(0.57,1.13)$ and $\delta\in(1.13,1.92)$ with 6 and 4 \deltaConnected\ components respectively.
Therefore, choosing $\delta$ in these intervals results in the most persistent \deltaConnected\ components. 

%% file: Figures/PersistenceCC/SE2/PersistenceDiagram.tex

\begin{tikzpicture}
\begin{axis}[
    ytick ={0,1,2,3,4.5},
    yticklabels = {0,1,2,3,$\infty$},
    yticklabel style = {align = center},
    xtick ={0,1,2,3,4.5},
    xticklabels = {0,1,2,3,$\infty$},
    xticklabel style = {align = center},
    xlabel = {Birth}, 
    xlabel style = {font = \large},
    ylabel = {Death}, 
    ylabel style = {font = \large},
    boxplot/variable width min target = 0,
    ymin = 0, ymax = 4.5, xmin = 0, xmax =4.5,
    axis lines*=left,
    legend style={at={(1.1,0.5)},
    anchor=west,legend columns=1, draw = none, cells = {anchor = west}},
    ]

\addplot[scatter, scatter/use mapped color={draw=none,fill=black, fill opacity = 1, draw opacity = 0},only marks]
table[row sep=newline, x = birth, y=death] {Figures/PersistenceCC/SE2/PersistenceData.txt};

\addplot[mark = none, black, dashed] coordinates {(0,0) (4,4)};

\end{axis}
\end{tikzpicture}

%% file: Figures/PersistenceCC/SE2/TotalNumberOfCCSE2r07.tex

\begin{tikzpicture}
\begin{axis}[
    xlabel = {$\delta$}, 
    xtick ={0,1,1.5,2,3,4,5,6,7},
    xticklabels ={0,1,\footnotesize\textcolor{red}{$\delta^*$},2,3,4,5,6,7},
    xticklabel style = {align = center, font = \large},
    axis x line=middle,
    xlabel style = {font = \Large,anchor=west},
    ylabel = {K}, 
    axis y line=middle,
    ylabel style = {font = \Large,anchor = south},
    ytick ={1,10,100},
    yticklabel style = {align = center, font = \large},
    boxplot/variable width min target = 0,
    ymin = 0, ymax = 110, xmin = 0, xmax = 3.5,
    axis lines*=left,
    legend style={at={(1.1,0.5)},
    anchor=west,legend columns=1, draw = none, cells = {anchor = west}},
    ]

\addplot[mark = none, blue] coordinates {(0,106) (0.282843,106)};
\addplot[mark = none, blue] coordinates {(0.282843,90) (0.351241,90)};
\addplot[mark = none, blue] coordinates {(0.351241,12) (0.565685,12)};
\addplot[mark = none, blue] coordinates {(0.565685,8) (1.13137,8)};
\addplot[mark = none, blue] coordinates {(1.13137,6) (1.91902,6)};
\addplot[mark = none, blue] coordinates {(1.91902,1) (8,1)};

\addplot[mark = none, black, dashed] coordinates {(0.282843,106) (0.282843,90)};
\addplot[mark = none, black, dashed] coordinates {(0.351241,90) (0.351241,12)};
\addplot[mark = none, black, dashed] coordinates {(0.565685,12) (0.565685,8)};
\addplot[mark = none, black, dashed] coordinates {(1.13137,8) (1.13137,6)};
\addplot[mark = none, black, dashed] coordinates {(1.91902,6) (1.91902,1)};

\addplot[mark = none, red, dashed] coordinates {(1.5,0) (1.5,120)};
    
\end{axis}
\end{tikzpicture}

%% file: Figures/PersistenceCC/SE2/ElementsPerCCSE2r07_colors.tex

\begin{tikzpicture}
\begin{axis}[
    xlabel = {$\delta$}, 
    xtick ={0,1,1.5,2,3,4,5,6,7},
    xticklabels ={0,1,\footnotesize\textcolor{red}{$\delta^*$},2,3,4,5,6,7},
    xticklabel style = {align = center, font = \large},
    axis x line=middle,
    xlabel style={font = \Large, anchor=west},
    ylabel = {size $|[g_i]|$}, 
    axis y line=middle,
    ylabel style = {font = \Large, anchor = south},
    ytick ={10,100},
    yticklabel style = {align = center, font = \large},
    boxplot/variable width min target = 0,
    ymin = 0, ymax = 110, xmin = 0, xmax =3.5,
    axis lines*=left,
    legend style={at={(1.1,0.5)},
    anchor=west,legend columns=1, draw = none, cells = {anchor = west}},
    ]

\addplot[mark = none, blue] coordinates {(1.91902,106) (8,106)};
\addlegendentry{1}
\addplot[mark = none, teal] coordinates {(0.351241,40) (1.91902,40)};
\addlegendentry{2}
\addplot[mark = none, green] coordinates {(0.565685,3) (1.13137,3)};
\addlegendentry{4}
\addplot[mark = none, orange] coordinates {(0.282843,3) (0.565685,3)};
\addlegendentry{8}
\addplot[mark = none, magenta] coordinates {(0.282843,1) (0.351241,1)};
\addlegendentry{82}
\addplot[mark = none, violet] coordinates {(0,1) (0.282843,1)};
\addlegendentry{106}

\addplot[mark = none, teal] coordinates {(0.351241,1) (0.565685,1)};
\addplot[mark = none, teal] coordinates {(0.565685,7) (1.91902,7)};
\addplot[mark = none, teal] coordinates {(1.13137,6) (1.91902,6)};

\addplot[mark = none, gray, dashed] coordinates {(0.282843,1) (0.282843,3)};
\addplot[mark = none, gray, dashed] coordinates {(0.351241,1) (0.351241,40)};
\addplot[mark = none, gray, dashed] coordinates {(0.565685,1) (0.565685,7)};
\addplot[mark = none, gray, dashed] coordinates {(1.13137,3) (1.13137,6)};
\addplot[mark = none, gray, dashed] coordinates {(1.91902,6) (1.91902,106)};

\addplot[mark = none, red, dashed] coordinates {(1.5,0) (1.5,130)};
    
\end{axis}
\end{tikzpicture}

%% file: Review_JMIV/4_3AffinityMatrices.tex
\section{Affinity Matrices between \texorpdfstring{$\delta$-Connected}{delta-Connected} Components}\label{sec:20230922:Affinity}
Now that we have described a way to choose the parameter $\delta$ we can calculate the \deltaConnected\ components by Algorithm~\ref{alg:20230605:connectedComponents}: \ref{alg:deltaconnectedComponents}. Then we would like to quantify how well-aligned the \deltaConnected\ components are as a whole. For this 
we introduce and analyze affinities between connected components in this section. Intuitively, `affinity' measures the proximity and alignment of \deltaConnected\ components. 

We define the affinity matrix $A=(a_{ij})$ by
\begin{align*}
    a_{ij}=\sup&\left\{\left(\frac{1}{\mu([g_i])}\int_{[g_i]}\varphi_t^\alpha\left(W_{[g_j]}^{(0)}\right)(h)^p\mathrm{d}h\right)^{1/p},\right.\\
    & \left. \left(\frac{1}{\mu([g_j])}\int_{[g_j]}\varphi_t^\alpha\left(W_{[g_i]}^{(0)}\right)(h)^p\mathrm{d}h\right)^{1/p}\right\}, \numberthis\label{eq:affinityMatrix}
\end{align*}
with measure normalisation $\mu([g_i]):= \int_{[g_i]} {\rm d}g$, and $\alpha>1,\; t>0, p\geq 1$ fixed, and where $W_{[g_i]}^{(0)}$ is given by
\begin{align}
    W_{[g_i]}^{(0)}(g)=\frac{D(g)\mathbbm{1}_{[g_i]}(g)}{\sup\limits_{h\in [g_i]}\left\{D(h)\right\}},\quad g\in G,\label{eq:def:affinityW0}
\end{align}
where $D:G\to[0,1]$ is a (non-zero) data term that can be freely chosen. In the experimental section, we chose the data term $D$ to be the absolute value
of the orientation score, defined in Eq.~\eqref{eq:orientationscoreWphif}, rescaled to be between 0 and~1.
Let us also recall that the flow operator $\varphi_{t}^{\alpha}$ was given by \eqref{eq:20230922:FlowReprMorphDilationWithKernel} and implemented by the morphological convolution in \eqref{eq:viscositySolutionW}. We say $a_{ij}$ denotes the affinity between \deltaConnected\ components $[g_i]$ and $[g_j]$.

In the computation of the affinity matrices
we fix the parameter $t$, based on the choice of $\alpha$ and the compact set $I$, such that the dilation of any \deltaConnected\ component $[g_i]$ has non-zero values for all $h\in I$. How this can be achieved is explained in the next proposition.
\begin{proposition}[Non-zero affinity between \deltaConnected\ components]\label{20230922:prop:IdentifyWhenAffinityIsNonZero}
    Let $\alpha>1$ be given and $\frac{1}{\alpha}+\frac{1}{\beta}=1$. Let $I=\bigcup\limits_{k=1}^K [g_k]\subset G$ compact be given, $[g_i]\neq[g_j]$ for all $i\neq j$. Then
    \begin{align*}
        \varphi_t^\alpha\left(W_{[g_i]}^{(0)}\right)(g)\neq 0\text{ for all }g\in I,\; [g_i]\in \tildeIdelta
    \end{align*}
    at least when $t>\left(\sup\limits_{q_1,q_2\in I}d(q_1,q_2)\right)^\alpha\beta^{1-\alpha}$.
\end{proposition}
\begin{proof}
    We need to show that one dilation step with morphological kernel $k_{t}^{\alpha}$ on $W_{[g_i]}^{(0)}(g)$ results in all non-zero values for all $g\in I$. That means that we start with calculating:
    \begin{align*}
        0\leq\varphi_t^\alpha\left(W_{[g_i]}^{(0)}\right)(g)&=-\left(\morphkern{t}{\alpha}\square -W_{[g_i]}^{(0)}\right)(g)\\&=\sup_{h\in G}\left\{W_{[g_i]}^{(0)}(h)-\frac{t}{\beta}\left(\frac{d(g,h)}{t}\right)^\beta\right\}\leq 1,
    \end{align*}
    where we find the inequalities by Lemma~\ref{lemma:propertiesMorphologicalDilations}. Then, we calculate the value of $t>0$ such that
    \begin{align*}
        0<\sup_{h\in I}\left\{W_{[g_i]}^{(0)}(h)-\frac{t}{\beta}\left(\frac{d(g,h)}{t}\right)^\beta\right\}\leq 1.
    \end{align*}
    By definition of $W_{[g_i]}^{(0)}$, for any element in the connected component $q\in [g_i]$ where the data does not vanish ($D(q)\neq 0$), we have
    \begin{align*}
        W_{[g_i]}^{(0)}(q)&-\frac{t}{\beta}\left(\frac{d(g,q)}{t}\right)^\beta\\&\leq \sup_{h\in I}\left\{W_{[g_i]}^{(0)}(h)-\frac{t}{\beta}\left(\frac{d(g,h)}{t}\right)^\beta\right\}\leq 1.
    \end{align*}
    We are sure that all affinities are nonzero when 
    \begin{align*}
         W_{[g_i]}^{(0)}(q)&-\frac{t}{\beta}\left(\frac{d(g,q)}{t}\right)^\beta>0, \text{ for all }g\in I.
    \end{align*}
    By definition of $W_{[g_i]}^{(0)}$ in \eqref{eq:def:affinityW0}, 
    there exists a $q\in [g_i]$ such that $W_{[g_i]}^{(0)}(q)=1$, so that for all $g\in I$
    \begin{align*}
        &W_{[g_i]}^{(0)}(q)-\frac{t}{\beta}\left(\frac{d(g,q)}{t}\right)^\beta=1-\frac{t}{\beta}\left(\frac{d(g,q)}{t}\right)^\beta > 0\\
        \Leftrightarrow\; &t>\left(\frac{d(g,q)^\beta}{\beta}\right)^{\frac{1}{\beta-1}}\!=\left(d(g,q)\beta^{-1/\beta}\right)^\alpha=\left(d(g,q)\right)^\alpha\beta^{1-\alpha}
    \end{align*}
    for all $g\in I$, and where $\frac{1}{\alpha}+\frac{1}{\beta}=1$ with $\alpha,\beta\geq 1$ and the last estimate holds because the distance $d(g,q)$ is positive. Then, since the inequality has to hold for all $q\in I$ and should be independent of the origin, we have
    \begin{align*}
        t>\left(\sup_{q_1,q_2\in I}d(q_1,q_2)\right)^\alpha\beta^{1-\alpha}
        \geq\left(\sup_{\substack{q_1\in I,\\q_2\in [g_i]}}d(q_1,q_2)\right)^\alpha\beta^{1-\alpha}. \textrm{ \hfill $\Box$}
    \end{align*}
\end{proof}

Algorithm~\ref{alg:20230605:affinityMatrices}:~\ref{alg:affinityMatrices} describes the algorithm to find the affinity matrices between the \deltaConnected\ components of any set $I\subset G$.

\begin{algorithm}
\caption{\texttt{find\_affinity}}
\customlabelnoprint{alg:affinityMatrices}{\texttt{find\_affinity}}
\textbf{Input:} \deltaConnected\ components $[g_1],\ldots,[g_K]$ determined by \ref{alg:deltaconnectedComponents}, $I=\bigcup\limits_{k=1}^K [g_k]\subset G$ region of interest for the connected components, kernel $\morphkern{t}{\alpha}$ with $\alpha>1$, $t>0$, data $D:G\to[0,1]$, parameter $p\geq 1$.\\
set $t>\left(\sup\limits_{q_1,q_2\in I} d(q_1,q_2)\right)^\alpha \beta^{1-\alpha}$ fixed\\
set $k=1$\\
for $k\leq K$\\
$\quad$ initialize $W_{[g_k]}^{(0)}(g)\coloneqq \frac{D(g)\cdot\mathbbm{1}_{[g_k]}(g)}{\sup\limits_{h\in [g_k]}\{D(h)\cdot\mathbbm{1}_{[g_k]}(h)\}}$\\
$\quad$ update $W_{[g_k]}^{(1)}(g)=\varphi_t^\alpha\left(W_{[g_k]}^{(0)}\right)(g)=-\left(\morphkern{t}{\alpha}\square - W_{[g_k]}^{(0)}\right)(g)$\\
$\quad$ for $1\leq l\leq K$\\
$\qquad$ set $\tilde{a}_{kl}
\coloneqq 
\frac{\left(\int_{[g_l]}\varphi_t^\alpha\left(W_{[g_k]}^{(0)}\right)(h)^p\mathrm{d}h\right)^{1/p}}{\left(\int_{[g_l]}\mathrm{d}h\right)^{1/p}}
$\\
$\quad$ update $k\coloneqq k+1$\\
end\\

for $i,j=1,\ldots,K$\\
$\quad $ set $a_{ij}=\sup\{\tilde{a}_{ij},\tilde{a}_{ji}\}$\\
end\\
\textbf{Output:} Affinity Matrix $A=(a_{ij})$ for all $i,j=1,\ldots,K$.
\label{alg:20230605:affinityMatrices}
\end{algorithm}

If the data $D(\cdot)$ is constant and non-zero in the connected components, we can give an upper and lower bound for the affinity of $[g_i]$ on $[g_j]$.
\begin{proposition}\label{prop:AffinityConstantData}
    Let the data $D:G\to[0,1]$ be such that $D(g)=c_{[g]}\neq 0$ for all $g\in I$, where $c_{[g]}$ is constant in the same \deltaConnected\ component. Let $\alpha>1$ be given and $\frac{1}{\alpha}+\frac{1}{\beta}=1$. Let $t>\sup\limits_{q_1,q_2\in I} d(q_1,q_2)^\alpha \beta^{1-\alpha}$ fixed. \\
    Then the affinities on and off-diagonal satisfy respectively
    \begin{align*}
        a_{ii}=1\quad\text{and}\quad 0<1-\frac{t}{\beta}\left(\frac{\Diameter}{t}\right)^\beta\leq a_{ij}<1-\frac{t}{\beta}\left(\frac{\delta}{t}\right)^\beta<1
    \end{align*}
    for all $[g_i]\neq[g_j]\subset I$, with diameter  $\Diameter=\sup\limits_{q_1,q_2\in I}d(q_1,q_2)$.
\end{proposition}
\begin{proof}
    First, we note that since the data term $D(g)=c_{[g]}\neq 0$ for all $g\in I$, with $c_{[g]}$ constant in the same \deltaConnected\ component, we have 
    \begin{align*}
        W_{[g_i]}^{(0)}(g)=\mathbbm{1}_{[g_i]}(g).
    \end{align*}
    Then, the morphological dilation of the initialization yields
    \begin{align}
        W_{[g_i]}^{(1)}(g)\coloneqq\varphi_{\alpha}^t\left(W_{[g_i]}^{(0)}\right)(g)
        =1-\inf_{h\in [g_i]}\frac{t}{\beta}\left(\frac{d(g,h)}{t}\right)^\beta\label{eq:W1gi-prop}
    \end{align}
    where we used the relation $\frac{1}{\alpha}+\frac{1}{\beta}=1$, with $\alpha>1$. The distance satisfies two conditions:
    \begin{enumerate}
        \item If $g\in [g_i]$, then $\inf\limits_{h\in[g_i]} d(g,h)=0$.\label{prop:itemize:gingi}
        \item If $g\notin [g_i]$, then $\delta<\inf\limits_{h\in[g_i]} d(g,h)\leq \Diameter\coloneqq \sup\limits_{g,h\in I}d(g,h)$. The upper bound follows from $[g_i]\subset I$. The lower bound follows from $g$ not being part of the \deltaConnected\ component $[g_i]$. Recall that by Definition \ref{def:20230922:deltaConnectedComponents}, we know that the distance between two \deltaConnected\ components is larger than $\delta$.\label{prop:itemize:gnotingi}
    \end{enumerate}
    Using this information on the distances, we can further simplify the expression in \eqref{eq:W1gi-prop} to
    \begin{align}
        W_{[g_i]}^{(1)}(g)=1 \text{ if }g\in [g_i], \label{eq:simple}
    \end{align}
    and
    \begin{align}
        1-\frac{t}{\beta}\left(\frac{\Diameter}{t}\right)^\beta\leq 
        W_{[g_i]}^{(1)}(g)
        <1-\frac{t}{\beta}\left(\frac{\delta}{t}\right)^\beta \text{ if }g\notin [g_i],\label{eq:prop:upperandlowerboundW1}
    \end{align}
    using the distance estimates in items \ref{prop:itemize:gingi} and \ref{prop:itemize:gnotingi}. Then, the affinity \eqref{eq:affinityMatrix} equals $a_{ij}=\sup \{\tilde{a}_{ij},\tilde{a}_{ji}\}$
   with
    \begin{align*}
        \tilde{a}_{ij}
        &\coloneqq\left(\frac{1}{\mu([g_j])}\int_{[g_j]}\left(
        W_{[g_i]}^{(1)}(h)
        \right)^p\mathrm{d}h\right)^{1/p}.
    \end{align*}
    Clearly if $i=j$ then $a_{ij}=a_{ii}=1$ by (\ref{eq:simple}).
 For $i\neq j$,  Eq.\!~\eqref{eq:prop:upperandlowerboundW1} provides an upper bound
    \begin{align*}
         \tilde{a}_{ij}
         &<\left(\frac{1}{\mu([g_j])}\int_{[g_j]}\left(1-\frac{t}{\beta}\left(\frac{\delta}{t}\right)^\beta\right)^p\mathrm{d}h\right)^{1/p}\\
         &=1-\frac{t}{\beta}\left(\frac{\delta}{t}\right)^\beta<1,\numberthis\label{eq:prop:upperbound}
    \end{align*}
    and a lower bound
    \begin{align*}
         \tilde{a}_{ij}
         &\geq\left(\frac{1}{\mu([g_j])}\int_{[g_j]}\left(1-\frac{t}{\beta}\left(\frac{\Diameter}{t}\right)^\beta\right)^p\mathrm{d}h\right)^{1/p}\\
         &=1-\frac{t}{\beta}\left(\frac{\Diameter}{t}\right)^\beta>0,\numberthis\label{eq:prop:lowerbound}
    \end{align*}
    where $1-\frac{t}{\beta}\left(\frac{\Diameter}{t}\right)^\beta>0$ follows from the choice of $t$ as was explained and proven in Prop.~\ref{20230922:prop:IdentifyWhenAffinityIsNonZero}. 
Thereby one has
    \begin{align*}
        0<1-\frac{t}{\beta}\left(\frac{\Diameter}{t}\right)^\beta\leq
        a_{ij}=\max\{\tilde{a}_{ij},\tilde{a}_{ji}\}
        &<1-\frac{t}{\beta}\left(\frac{\delta}{t}\right)^\beta<1,
    \end{align*}
  and the result follows. 
    \qed
\end{proof}

\begin{remark}
    Note that an alternative to the affinity matrices would be optimal transport. One can develop a model on the Lie group $G$ that finds the optimal way to transform one vessel into another. The transformations corresponding to the lowest (Wasserstein-)distance \cite{bonPai2025optimal} are the most likely to be connected.
\end{remark}

%% file: Review_JMIV/5Experiments.tex
\section{Experiments}\label{sec:20230922:Experiments}
In the previous section, we introduced several concepts and algorithms, such as a formal algorithm to identify connected components (\ref{alg:deltaconnectedComponents}) and an algorithm to calculate affinity matrices (\ref{alg:affinityMatrices}). Here, we will show the results of these algorithms applied to several images of the STAR-dataset \cite{Zhang2016,AbbasiSureshjani}. All experiments are performed on the Lie group $G=\SE{2}$, and the Mathematica notebooks are available via \cite{Berg2024connectedNotebooks}.

We start with discussing the experimental set-up in Sec.~\ref{sec:experimentalSetup} and our approach to identifying the \deltaConnected\ components and the corresponding results in Sec.~\ref{sec:Exp:CC}. Then, we discuss the results of the affinity matrices on some of the \deltaConnected\ component experiments in Sec.~\ref{sec:Exp:AM}.

\subsection{Experimental Set-Up}\label{sec:experimentalSetup}
In the experiments, we will identify the \deltaConnected\ components (flowchart in Fig.~\ref{fig:flowchartExperimentsCC}) in retinal images from the STAR-dataset \cite{Zhang2016,AbbasiSureshjani} on the Lie group $G=\SE{2}$. The retinal images in this dataset are standard 2D-images. Therefore, we first explain how we prepare the images for processing in the Lie group $\SE{2}$.


\begin{figure*}[ht!]
    \centering
        \scalebox{1.0}{\input {Figures/FlowchartCC}}
    \caption{Flowchart of the executed steps in the connected component experiments shown in Fig.~\ref{fig:CCIm5}-\ref{fig:CCIm24}.}
    \label{fig:flowchartExperimentsCC}
\end{figure*}

\begin{figure*}[ht!]
    \centering
    \scalebox{1.0}{\input {Figures/FlowchartAM}}
    \caption{Flowchart of the executed steps in the affinity matrices experiments shown in Fig.~\ref{fig:CCIm21}-\ref{fig:CCIm24}.}
    \label{fig:flowchartExperimentsAM}
\end{figure*}

We consider an input image $f:\mathbb{R}^2\to\mathbb{R}$. We lift it to the space of positions and orientations $\SE{2}$ (by creating an orientation score) to disentangle crossing structures. The orientation score $W_\phi f$ is calculated by a convolution with a rotating anisotropic wavelet $\phi$ and is given by Eq.~\eqref{eq:orientationscoreWphif}. For $\phi$, we use real-valued cake wavelets (see \cite{InvertibleOrientationScores,InvertibleOrientationScores2}). We use 32 orientations for all experiments.

From the orientation score $W_\phi f$, we calculate a crossing-preserving vesselness $\mathcal{V}^{\SE{2}}(W_\phi f)$, which is done as described in \cite[Appendix~D]{Berg2024geodesic}. In all experiments, we use parameter settings $\xi=1$, $\zeta=1$, $s=\{1.5,2,2.5,3\}$, $\sigma_{s,Ext}=0$ and $\sigma_{a,Ext}=0$.

From this crossing-preserving vesselness, we compute a cost function $\mathcal{C}=1/(1+\lambda \mathcal{V}^{\SE{2}}(W_\phi f)^p)$ with parameters $\lambda>0$ and $p>0$. Then, we binarize the cost function using Otsu's method to identify the binarization threshold value. Because the width of the structures varies within the same image, we choose to identify the centerline of the binarized structure, using the approach proposed in \cite{Lee1994skeleton}. Since this is done in an isotropic way, we compensate for that by slightly dilating the structures in the $\mathcal{A}_1$-direction; the principal direction of the local structure (see Eq.~\eqref{eq:MetricTensorFrameA} for the explicit formula of $\mathcal{A}_1$). The support of the resulting binarized function on $\SE{2}$ together with the lifted bifurcation information (for all orientations at bifurcation position equal to 1) is then used for the reference set $I$ in the \deltaConnected\ component algorithm, more specifically, we identified $\mathbbm{1}_I$.

Then, we use the connected component algorithm as introduced in \ref{alg:deltaconnectedComponents} to identify all connected components in the lifted image. To visualize the output of the \deltaConnected\ component algorithm, we define the function $f^{CC}:\SE{2}\to\mathbb{N}$ by
\begin{align*}
    f^{CC}(x,y,\theta)=\begin{cases}
        i & \text{if }(x,y,\theta)\in [g_i]\\
        0 & \text{else.}
    \end{cases}
\end{align*}
The \deltaConnected\ components are visualized by projecting them back onto $\mathbb{R}^2$ (taking per location the maximum over all orientations), i.e.
\begin{align*}
    f_{\text{out}}^{CC}(x,y)=\max_{\theta\in[-\pi,\pi)}f^{CC}(x,y,\theta).
\end{align*}

Moreover, we also calculate the affinity matrices (flowchart in Fig.~\ref{fig:flowchartExperimentsAM}). The algorithm used to compute them, uses the \deltaConnected\ components as input. Therefore, one first follows the steps for computing the \deltaConnected\ components described in the previous paragraph. The projection step to visualize the results can be skipped as the identification of the affinity matrices also happens in $\SE{2}$. Instead, we perform the affinity matrices algorithm as given in \ref{alg:affinityMatrices} to identify all affinity matrices. Here, we use the orientation score data $W_\phi f$ as the data term $D=|W_\phi f|/\sup |W_\phi f|$. To visualize the results of the output of the affinity matrices algorithm, we determine a threshold value of $T$ and group the \deltaConnected\ components having an affinity higher than this threshold value $T$, i.e.
\begin{align*}
    \tilde{C}_{i}=\bigcup_{j\in \mathbb{I}_i} [g_j],
\end{align*}
with $\mathbb{I}_i$ the set containing
the vertices that belong to the same connected component of the adjacency graph \mbox{$A>T$}.
Consequently, multiple \deltaConnected\ components will be grouped in the new visualization
\begin{align*}
    f^{AM}(x,y,\theta)=\begin{cases}
        i& \text{if }(x,y,\theta)\in \tilde{C}_i\\
        0& \text{else.}
    \end{cases}
\end{align*}
To visualize the results, we again project the affinity matrices results back onto $\mathbb{R}^2$ by taking per location the maximum over all orientations
\begin{align*}
    f_{\text{out}}^{AM}(x,y)=\max_{\theta\in[-\pi,\pi)}f^{AM}(x,y,\theta)
\end{align*}

The additional parameter settings used to process all images can be found in Table~\ref{tab:parameters}. The Mathematica notebooks are publicly available via \cite{Berg2024connectedNotebooks}.

\begin{table*}[ht!]
    \centering
    \begin{tabular}{c c|
    cc|
    cc|
    lll|
    cccc}
        Figure & Image & 
        \multicolumn{2}{|c|}{Cost} & 
        \multicolumn{2}{|c|}{Dilation} & 
        \multicolumn{3}{|c|}{\deltaConnected\ Components} &  
        \multicolumn{4}{|c}{Affinity Matrices}\\
        
          & &
         $\lambda$ & $p$ &
         $(w_1,w_2,w_3)$ & $\alpha$ &
         $(\tilde{w}_1,\tilde{w}_2,\tilde{w}_3)$, $w_i=k\tilde{w}_i$ & $\alpha$ & $\delta$ &
         $(w_1,w_2,w_3)$ & $\alpha$ & $T$ & $p$\\\hline
         
         Fig.~\ref{fig:CCIm5} & STAR48 & 
         $100$ & $3$ &
         $(0.2,1.5,50)$ & $1.3$ &
         $(0.06,0.7,2)$ & $1$ & $1$ & 
         -&-&-&-\\
         
        Fig.~\ref{fig:CCIm9} & STAR13 &
         $\phantom{0}50$ & $\phantom{.5}3.5$ & 
         $(0.1,1.5,50)$ & $1.3$ & 
         $(0.1\phantom{0},0.7,4)$ & $1$ & $1$ & 
         -&-&-&-\\
         
         Fig.~\ref{fig:CCIm21} & STAR34 & 
         $200$ & $3$ & 
         $(0.2,1.5,50)$ & $1.3$ &
         $(0.2\phantom{0},1\phantom{.7},2)$ & $1$ & $1$ & 
         $(0.1,4,4\phantom{.1})$ & $2$ & $0.99957$ & $2$ \\
         
         Fig.~\ref{fig:CCIm23} & STAR37 & 
         $100$ & $3$ &
         $(0.2,1.5,50)$ & $1.3$ &
         $(0.08,0.7,4)$ & $1$ & $1$ &
         $(0.5,2,0.5)$ & $2$ & $0.9985\phantom{9}$ & $2$ \\
         
         Fig.~\ref{fig:CCIm24} & STAR38 & 
         $100$ & $3$ & 
         $(0.2,1.5,50)$ & $1.3$ & 
         $(0.2\phantom{0},1\phantom{.7},2)$ & $1$ & $1$ &
         $(0.5,2,0.1)$ & $2$ & $0.993\phantom{99}$ & $2$ \\
    \end{tabular}
    \caption{Parameter settings for preprocessing of data and calculation of \deltaConnected\ components and affinity matrices. The parameter $k$ in the \deltaConnected\ components is given by $k=(13/3)^{6/13}$.}
    \label{tab:parameters}
\end{table*}

\begin{remark}
    The output of the \deltaConnected\ component algorithm primarily depends on the metric tensor weights, which directly affect the left-invariant Riemannian distance used to compute point-to-point distances. By selecting weight parameters based on image resolution, one can ensure that the resulting \deltaConnected\ components remain consistent across different discretizations—up to rounding errors.
\end{remark}

\subsection{$\delta$-Connected Components}\label{sec:Exp:CC}
\begin{figure}[ht!]
    \centering
    \begin{subfigure}[t]{0.2\textwidth}
         \centering
        \includegraphics[width=\textwidth]{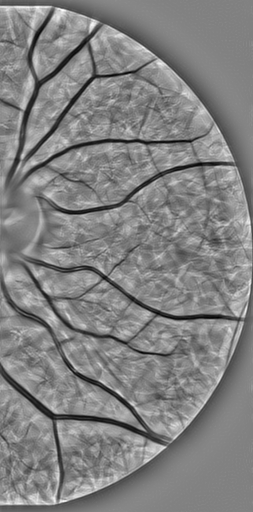}
        \caption{Underlying Image (STAR48).}
    \end{subfigure}~
    \begin{subfigure}[t]{0.2\textwidth}
         \centering
        \includegraphics[width=\textwidth]{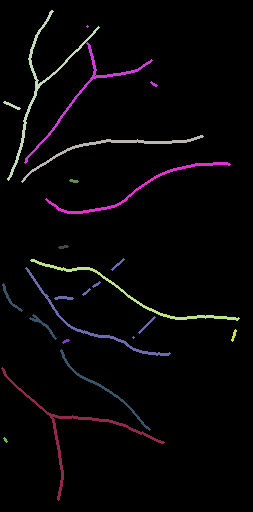}
        \caption{Connected components calculated by the flowchart in Fig.~\ref{fig:flowchartExperimentsCC} and parameter settings in Table~\ref{tab:parameters}.}
    \end{subfigure}
    \caption{The output of the \deltaConnected\ component algorithm executed on a thinning of the calculated vesselness of the given image. The parameter settings can be found in Table~\ref{tab:parameters}.}
    \label{fig:CCIm5}
\end{figure}

\begin{figure}[ht!]
    \centering
    \begin{subfigure}[t]{0.2\textwidth}
         \centering
        \includegraphics[width=\textwidth]{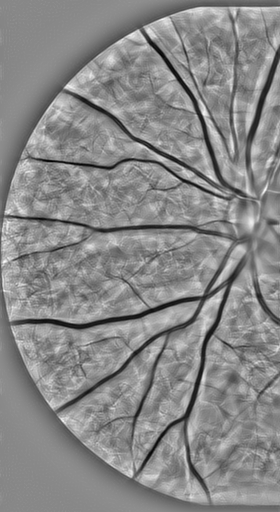}
        \caption{Underlying Image (STAR13).}
    \end{subfigure}~
    \begin{subfigure}[t]{0.2\textwidth}
         \centering
        \includegraphics[width=\textwidth]{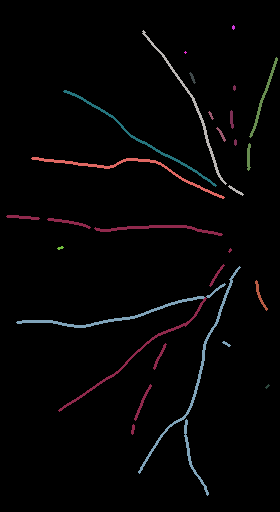}        
        \caption{Connected components calculated by the flowchart in Fig.~\ref{fig:flowchartExperimentsCC} and parameter settings in Table~\ref{tab:parameters}.}
    \end{subfigure}
    \caption{The output of the \deltaConnected\ component algorithm executed on a thinning of the calculated vesselness of the given image. The parameter settings can be found in Table~\ref{tab:parameters}.}
    \label{fig:CCIm9}
\end{figure}

\begin{figure}[ht!]
    \centering
    \begin{subfigure}[t]{0.48\textwidth}
         \centering
        \includegraphics[width=\textwidth]{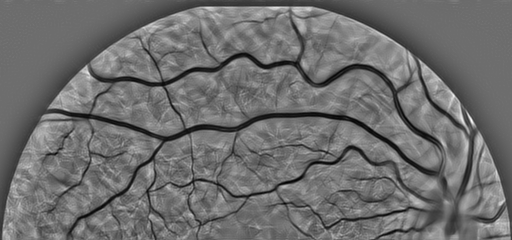}
        \caption{Underlying Image (STAR34).}
    \end{subfigure}
    \begin{subfigure}[t]{0.48\textwidth}
         \centering
        \includegraphics[width=\textwidth]{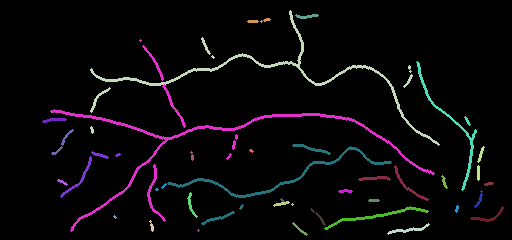}
        \caption{Connected components calculated by the flowchart in Fig.~\ref{fig:flowchartExperimentsCC} and parameter settings in Table~\ref{tab:parameters}.}
        \label{fig:Im21_CC}
    \end{subfigure}
    \begin{subfigure}[t]{0.48\textwidth}
         \centering
        \begin{tikzpicture}[scale=1.1]
            \node[] (pic) at (0,0) {\includegraphics[width=\textwidth]{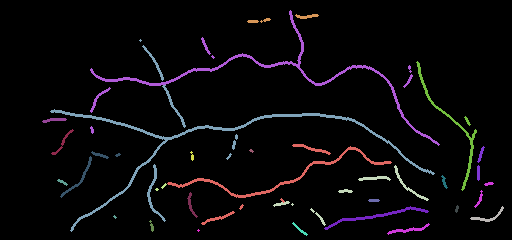}};
            \draw[white,thick] (0,1.4) rectangle (1,1.7);
            \draw[white,thick] (0.3,-1.2) rectangle (1.1,-1.5);
            \draw[white,thick] (1.3,-.8) rectangle (2.3,-1.2);
        \end{tikzpicture}
        \caption{Affinity Matrices calculated on the connected components in Fig.~\ref{fig:Im21_CC}. Newly grouped \deltaConnected\ components are indicated by a white box.}
        \label{fig:Im21_AM}
    \end{subfigure}
    \caption{The output of the \deltaConnected\ component algorithm and affinity matrices algorithm executed on a thinning of the calculated vesselness of the given image. The parameter settings can be found in Table~\ref{tab:parameters}.}
    \label{fig:CCIm21}
\end{figure}

\begin{figure}[ht!]
    \centering
    \begin{subfigure}[t]{0.48\textwidth}
         \centering
        \includegraphics[width=\textwidth]{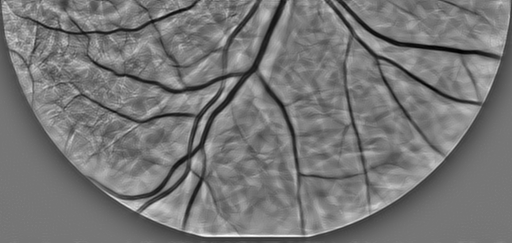}
        \caption{Underlying Image (STAR37).}
    \end{subfigure}
    \begin{subfigure}[t]{0.48\textwidth}
         \centering
        \includegraphics[width=\textwidth]{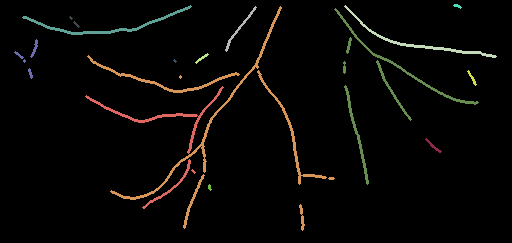}
        \caption{Connected components calculated by the flowchart in Fig.~\ref{fig:flowchartExperimentsCC} and parameter settings in Table~\ref{tab:parameters}.}
        \label{fig:Im23_CC}
    \end{subfigure}
    \begin{subfigure}[t]{0.48\textwidth}
         \centering
        \begin{tikzpicture}[scale=1.1]
            \node[] (pic) at (0,0) {\includegraphics[width=\textwidth]{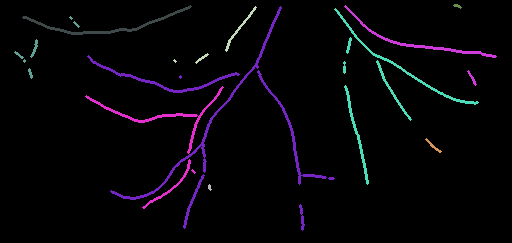}};
            \draw[white,thick] (-1.4,0.8) rectangle (-0.2,1.3);
            \draw[white,thick] (3.0,0.5) rectangle (3.4,1.2);
        \end{tikzpicture}
        \caption{Affinity Matrices calculated on the connected components in Fig.~\ref{fig:Im23_CC}. Newly grouped \deltaConnected\ components are indicated by a white box.}
        \label{fig:Im23_AM}
    \end{subfigure}
    \caption{The output of the \deltaConnected\ component algorithm and affinity matrices algorithm executed on a thinning of the calculated vesselness of the given image. The parameter settings can be found in Table~\ref{tab:parameters}.}
    \label{fig:CCIm23}
\end{figure}

\begin{figure}[ht!]
    \centering
    \begin{subfigure}[t]{0.48\textwidth}
         \centering
        \includegraphics[width=\textwidth]{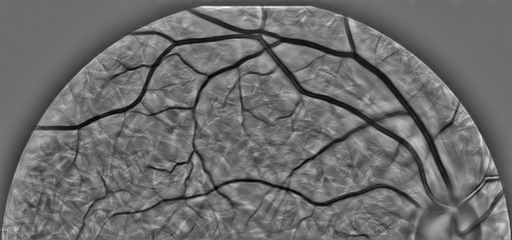}
        \caption{Underlying Image (STAR38).}
    \end{subfigure}
    \begin{subfigure}[t]{0.48\textwidth}
         \centering
        \includegraphics[width=\textwidth]{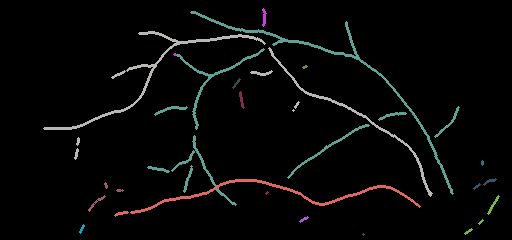}
        \caption{Connected components calculated by the flowchart in Fig.~\ref{fig:flowchartExperimentsCC} and parameter settings in Table~\ref{tab:parameters}.}
        \label{fig:Im24_CC}
    \end{subfigure}
    \begin{subfigure}[t]{0.48\textwidth}
         \centering
        \begin{tikzpicture}[scale=1.1]
            \node[] (pic) at (0,0) {\includegraphics[width=\textwidth]{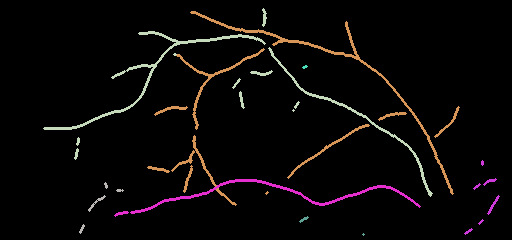}};
            \draw[white,thick] (-2.8,-1.8) rectangle (-2.2,-1.1);
            \draw[white,thick] (0,-1.3) rectangle (0.7,-0.8);
            \draw[white,thick] (3.2,-1.5) rectangle (3.8,-0.6);
            \draw[white,thick] (0,1.0) rectangle (0.3,1.7);
            \draw[white,thick] (-0.4,0.8) rectangle (0.2,0.3);
        \end{tikzpicture}
        \caption{Affinity Matrices calculated on the connected components in Fig.~\ref{fig:Im24_CC}. Newly grouped \deltaConnected\ components are indicated by a white box.}
        \label{fig:Im24_AM}
    \end{subfigure}
    \caption{The output of the \deltaConnected\ component algorithm and affinity matrices algorithm executed on a thinning of the calculated vesselness of the given image. The parameter settings can be found in Table~\ref{tab:parameters}.}
    \label{fig:CCIm24}
\end{figure}

We identify the \deltaConnected\ components using the method described in Sec.~\ref{sec:experimentalSetup} and visualized in the flowchart in Fig.~\ref{fig:flowchartExperimentsCC}. We identify the components for five different images from the STAR-dataset \cite{Zhang2016,AbbasiSureshjani}, and show the results in Fig.~\ref{fig:CCIm5}-\ref{fig:CCIm24}.

We start with the retinal image STAR48 in Fig.~\ref{fig:CCIm5}. Due to the binarization, we see that not all vascular structures are nicely connected. However, the \deltaConnected\ component algorithm can compensate for interrupted vascular structures by choosing the right threshold value for $\delta$, and the correct distance parameters $w_1$, $w_2$ and $w_3$. We see that the algorithm has correctly grouped most of the segments that belong to the same vessel in the underlying image. 

The \deltaConnected\ components of the retinal image STAR13 are visualized in Fig.~\ref{fig:CCIm9}. The vascular structure in this image does not contain a lot of bifurcations or crossing structures. The algorithm is good at identifying components that correspond to the structures in the underlying image. We also see that the choice of the threshold in the binarization (found with Otsu's method) has a big influence on the input in the \deltaConnected\ component algorithm; many vascular structures are interrupted, and therefore the algorithm needs to compensate for that as well. The results are good as long as the gaps are smaller than a certain threshold value. If the gaps are bigger, the \deltaConnected\ component algorithm is not able to identify them as belonging to the same vascular structure.

In Fig.~\ref{fig:CCIm21}, we applied the algorithm to STAR34. The results are shown in Fig.~\ref{fig:Im21_CC}. The vessels in the image are more sinuous. We see that this causes some challenges if the vascular structure gets interrupted in the binarization. This is because we need to choose the metric parameters in the \deltaConnected\ component algorithm. We chose these parameters, cf. column 7\&8 in Table~\ref{tab:parameters}, such that forward movement is allowed, but sideways movement and changing orientation are not, to avoid crossing structures being connected. However, one does need to change orientation and forward movement at interrupted tortuous structures. Consequently, the \deltaConnected\ component algorithm has trouble connecting the vessel segments that are interrupted at highly tortuous parts.

In Fig.~\ref{fig:CCIm23}, we performed the algorithm on the retinal image STAR37. The algorithm groups the vessel segments correctly using the chosen parameters. However, it does not identify full vascular trees, but only parts of it. This is due to the relatively large spatial gaps between vessel parts, often also changing orientation.

Lastly, we look at STAR38 in Fig.~\ref{fig:Im24_CC}. The \deltaConnected\ component algorithm can identify the large vascular structures correctly. Some small vessels are not correctly connected to the main vessel at the bifurcations, but all vessel segments are correctly connected. 

We conclude that this \deltaConnected\ component algorithm allows us to identify parts of vascular trees. Additionally, it can differentiate between different structures when they are crossing. However, when the vascular structure is interrupted at a very tortuous part, the algorithm can have difficulties connecting the right parts, depending on the chosen metric parameters.

\subsubsection{Runtime and Accuracy}
We compare two different models: a) our \deltaConnected\ component algorithm in $\SE2$, and b) the classical connected components on $\mathbb{R}^2$ where voxels are connected if they share at least a corner. We use the STAR dataset to compare both methods.

We report the calculation times of both methods. For b), we show the calculation times both for the standard Mathematica implementation and our own (non-optimized) method to identify \deltaConnected\ components applied to $\mathbb{R}^2$, 
using only 1 orientation layer.

Lastly, we report the accuracy of the identified components compared to the ground truth for all images in the STAR-dataset. We do this using two different measures $E_{split}$ and $E_{merge}$, where $E_{split}$ measures how many components one vascular tree of the ground truth is divided on average, and $E_{merge}$ measures how many vascular trees are covered by a single component on average, i.e.,
\begin{align}
    E_{split}=\frac{\sum\limits_{i=1}^{N}\left(\sum\limits_{k=1}^K |[g_k]\cap T_i|\right)}{N};\\
    E_{merge}=\frac{\sum\limits_{k=1}^{K}\left(\sum\limits_{i=1}^N |[g_k]\cap T_i|\right)}{K},
\end{align}
where $T_i$ represents one of the $N$ vascular trees, and $[g_k]$ one of the identified $K$ \deltaConnected\ components.

The calculation times of the \deltaConnected\ component algorithm are the longest (cf. Fig.~\ref{fig:calculationTimesCC}). This is as we expect, as we have not optimized the algorithm. The calculation times for (non-)optimized $\mathbb{R}^2$ components suggest that one can significantly improve calculation times by relying on more sophisticated methods originally designed for the classical connected component algorithm.

In the accuracy plots in Figs.~\ref{fig:accuracyEsplit}~and~\ref{fig:accuracyEmerge}, each point represents a different image. The points indicated by an `o' show improved results using the \deltaConnected\ component algorithm, whereas those indicated by an `x' performed better using the classical connected component algorithm. We see that in a significant majority of the cases, the \deltaConnected\ component algorithm outperforms the classical connected components.
\begin{figure}
    \centering
    \begin{subfigure}[t]{0.48\textwidth}
    \centering
    \begin{tikzpicture}
    \begin{axis}[
        xlabel = {model},
        ylabel= {calculation time},
        xtick ={1,2,3},
        xticklabel style = {align = center, font = \small},
        xticklabels = {our method\\ in SE(2), our method\\ in $\mathbb{R}^2$, standard method\\ in $\mathbb{R}^2$},
        ytick={0.001,0.01,0.1,1,10,100,1000},
        scale only axis=true,
        ymode = log,
        scale=0.7,
        scatter/classes={%
            se2={mark=x,draw=black,opacity=0.5},
            r2={mark=x,draw=black,opacity=0.5},
            fastr2={mark=x,draw=black,opacity=0.5}}]
        \addplot[scatter,only marks,%
            scatter src=explicit symbolic]%
        table[meta=label] {
        a times label
        1 146.48 se2
        1 223.314 se2
        1 269.526 se2
        1 215.637 se2
        1 192.838 se2
        1 168.37 se2
        1 150.881 se2
        1 131.838 se2
        1 238.618 se2
        1 232.363 se2
        1 227.525 se2
        1 168.174 se2
        1 154.005 se2
        1 253.033 se2
        1 177.614 se2
        1 162.878 se2
        1 185.244 se2
        1 271.338 se2
        1 152.186 se2
        2 30.6339 r2 
        2 38.9753 r2 
        2 54.9913 r2 
        2 48.1003 r2 
        2 40.2378 r2 
        2 31.8566 r2 
        2 37.542 r2 
        2 29.7748 r2 
        2 45.5991 r2 
        2 32.4323 r2 
        2 39.7954 r2 
        2 33.1301 r2 
        2 35.5333 r2 
        2 45.9364 r2 
        2 69.1813 r2 
        2 29.4604 r2 
        2 57.8239 r2 
        2 25.0486 r2 
        2 36.7287 r2
        3 0.0180056 fastr2
        3 0.0056523 fastr2
        3 0.0046587 fastr2
        3 0.0045451 fastr2
        3 0.0053311 fastr2
        3 0.0051041 fastr2
        3 0.005467 fastr2
        3 0.0055715 fastr2
        3 0.0056881 fastr2
        3 0.0051134 fastr2
        3 0.0051419 fastr2
        3 0.0058608 fastr2
        3 0.0042747 fastr2
        3 0.004481 fastr2
        3 0.005089 fastr2
        3 0.0050566 fastr2
        3 0.0052091 fastr2
        3 0.0053437 fastr2
        3 0.0051033 fastr2
            };
    \end{axis}
    \end{tikzpicture}
    \caption{Calculation times for different methods to identify connected components.}
    \label{fig:calculationTimesCC}
    \end{subfigure}
    \begin{subfigure}[t]{0.45\textwidth}
    \centering
    \begin{tikzpicture}
    \begin{axis}[
        xlabel = {$\SE2$},
        ylabel= {$\mathbb{R}^2$},
        xmin = 0, xmax = 30,
        ymin = 0, ymax = 30,
        scale only axis=true,
        axis equal image,
        scale=0.7,
        scatter/classes={%
            a={mark=o,draw=black},
            b={mark=x,draw=black}}]
        \addplot[scatter,only marks,%
            scatter src=explicit symbolic]%
        table[row sep=newline, x=SE2,y=R2,meta=label] {Review_JMIV/Accuracies/accuracy_CC_A.txt};
        \path[name path = A] (axis cs:0,0) -- (axis cs:35,35);
        \path[name path = B] (axis cs:35,0) -- (axis cs:35,35);
        \path[name path = C] (axis cs:0,0) -- (axis cs:0,35);
        \addplot+ [red, fill opacity=0.2] fill between [of=A and B];
        \addplot+ [green, fill opacity=0.2] fill between [of=A and C];
        \draw[dotted] (axis cs:0,0)--(axis cs:35,35);
    \end{axis}
    \end{tikzpicture}
    \caption{\# found components covering one real one: $E_{split}$.}
    \label{fig:accuracyEsplit}
    \end{subfigure}
    \begin{subfigure}[t]{0.45\textwidth}
    \centering
    \begin{tikzpicture}
    \begin{axis}[
        xlabel = {$\SE2$},
        ylabel= {$\mathbb{R}^2$},
        xmin = 0.9, xmax = 1.3,
        ymin = 0.9, ymax = 1.3,
        scale only axis=true,
        axis equal image,
        scale=0.7,
        scatter/classes={%
            a={mark=o,draw=black},
            b={mark=x,draw=black}}]
        \addplot[scatter,only marks,%
            scatter src=explicit symbolic]%
        table[row sep=newline, x=SE2,y=R2,meta=label] {Review_JMIV/Accuracies/accuracy_CC_B.txt};
        \path[name path = A] (axis cs:0,0) -- (axis cs:3,3);
        \path[name path = B] (axis cs:3,0) -- (axis cs:3,3);
        \path[name path = C] (axis cs:0,0) -- (axis cs:0,3);
        \addplot+ [red, fill opacity=0.2] fill between [of=A and B];
        \addplot+ [green, fill opacity=0.2] fill between [of=A and C];
        \draw[dotted] (axis cs:0,0)--(axis cs:3,3);
    \end{axis}
    \end{tikzpicture}
    \caption{\# real components covered by one found: $E_{merge}$.}
    \label{fig:accuracyEmerge}
    \end{subfigure}
    \caption{Visualizations of the runtime and accuracy measures $E_{split}$ and $E_{merge}$. In the green upper-triangle the \deltaConnected\ component outperforms the classical algorithm (indicated by an `o'), whereas in the red lower-triangle the classical algorithm is more accurate than the presented \deltaConnected\ component algorithm (points indicated by an `x').}
    \label{fig:runtimeAccuracy}
\end{figure}

\subsection{Affinity Matrices}\label{sec:Exp:AM}
In the previous section, we discussed the experimental results of the \deltaConnected\ component algorithm. We demonstrated that the algorithm can identify parts of vessels belonging to the same vascular structure, but often does not identify full vascular trees. We will use the affinity matrices to group different \deltaConnected\ components that are most likely to belong to the same vascular structure based on their local alignment and proximity. We have performed the affinity matrices algorithm on the \deltaConnected\ component results of Fig.~\ref{fig:CCIm21}-\ref{fig:CCIm24}. We use the orientation score data $|W_\phi f|$ as the data term in the affinity matrices algorithm given in \ref{alg:affinityMatrices}.

First, we calculated the affinity matrices on the \deltaConnected\ component output in Fig.~\ref{fig:Im21_AM}. Thresholding on the affinities allows us to group different components, which results in slightly more complete vascular trees. We cannot pick parameter settings such that all parts that belong to the same vessel are connected without connecting them to other structures. However, the output in Fig.~\ref{fig:Im21_AM} is a more complete vascular tree classification than the \deltaConnected\ component output in Fig.~\ref{fig:Im21_CC}. We have marked the correctly grouped vascular structures with a white box.

In Fig.~\ref{fig:Im23_AM}, we have also applied the affinity matrices algorithm to the output of the \deltaConnected\ component algorithm in Fig.~\ref{fig:Im23_CC}. The thresholded output of the affinity matrices has correctly grouped different \deltaConnected\ components that belong to the same vascular structure. The newly grouped structures have been indicated with a white box. We note that the output is closer to the actual underlying vascular structure. 
Still, some cases have not been grouped correctly (due to the imperfect data term $D$).

Last, we have calculated the affinity matrices for the retinal image of STAR38 (cf. Fig~\ref{fig:CCIm24}). The \deltaConnected\ component output in Fig.~\ref{fig:Im24_CC} has been used as input for \ref{alg:affinityMatrices}. The output of this algorithm has been thresholded and the \deltaConnected\ components having an affinity higher than the threshold are grouped. The output is visualized in Fig.~\ref{fig:Im24_AM}. Here, full vascular trees are identified completely and correctly. The white boxes indicate the changes compared to the \deltaConnected\ component output.

We found that thresholding the affinity matrices improves the results; \deltaConnected\ components are grouped such that more complete vascular trees are identified. It is important to note that the choice of the metric parameters of Table~\ref{tab:parameters} and threshold value $T$ are important to the output.

%% file: Figures/FlowchartCC.tex
\begin{tikzpicture}[xscale=0.8,yscale=1.0]
    \node[draw, align = center] (input) at (0,1) {Input image $f$};
    \node[draw, align = center] (os) at (0,-1) {Orientation scores\\$W_\phi f$};
    \node[draw, align = center] (V) at (4,-1) {vesselness \\ $\mathcal{V}^{\SE{2}}(W_\phi f)$};
    \node[draw, align = center] (bin) at (7,-1) {binarize};
    \node[draw, align = center] (thinDil) at (10,-1) {thinning \& \\dilation in \\ $\mathcal{A}_1$-direction};
    \node[draw, align = center] (bif) at (13.5,-1) {set \\bifurcations\\ to 1};
    \node[draw, align = center] (CC) at (18,-1) {Execute Algorithm~\ref{alg:20230605:connectedComponents}:\\\ref{alg:deltaconnectedComponents}};
    \node[draw, align = center] (output) at (18,1) {\deltaConnected\ component\\ output};

    \draw[-latex] (input) to node[right] {lifting to $SE(2)$} (os);
    \draw[-latex] (os) to (V);
    \draw[-latex] (V) to (bin);
    \draw[-latex] (bin) to (thinDil);
    \draw[-latex] (thinDil) to (bif);
    \draw[-latex] (bif) to (CC);
    \draw[-latex] (CC) to node[left] {projecting to $\mathbb{R}^2$} (output);
    
    \draw[dashed, -latex] (input) to (output);
\end{tikzpicture}

%% file: Figures/FlowchartAM.tex
\begin{tikzpicture}[xscale=0.8,yscale=1.0]
    \node[draw, align = center] (input) at (0,1) {Input image $f$};
    \node[draw, align = center] (os) at (0,-1) {Orientation scores\\$W_\phi f$};
    \node[draw, align = center] (CC) at (6,-1) {\deltaConnected\ component \\ algorithm via flowchart\\ in Fig.~\ref{fig:flowchartExperimentsCC}};
    \node[draw, align = center] (AM) at (12,-1) {Execute Algorithm~\ref{alg:20230605:affinityMatrices}:\\\ref{alg:affinityMatrices}};
    \node[draw, align = center] (threshold) at (18,-1) {Threshold on output \\ and group \\\deltaConnected\ components};
    \node[draw, align = center] (output) at (18,1) {Affinity matrices output};

    \draw[-latex] (input) to node[right] {lifting to $SE(2)$} (os);
    \draw[-latex] (os) to (CC);
    \draw[-latex] (CC) to (AM);
    \draw[-latex] (AM) to (threshold);
    \draw[-latex] (threshold) to node[left] {projecting to $\mathbb{R}^2$} (output);
    
    \draw[dashed, -latex] (input) to (output);
\end{tikzpicture}

%% file: Review_JMIV/6ConclusionFutureWork.tex
\section{Conclusion and Future Work}\label{sec:20230922:DiscussionConclusion}
In this article, we have introduced a way to identify so-called \deltaConnected\ components on a Lie group $G$. These components consist of sets (of points) with a maximum distance $\delta$ from each other. First, we introduced the general idea behind the algorithm and some theoretical background on morphological dilations. Then, we connected this theory to the general algorithm, resulting in the \deltaConnected\ component algorithm, stated in \ref{alg:deltaconnectedComponents}.

We studied the convergence of the \deltaConnected\ component algorithm in Theorem~\ref{th:20230605:ConvergenceConnectedComponentAlgorithm}. We proved that the algorithm always finishes in a finite number of iteration steps. Subsequently, we discussed the choice of the parameter $\delta$. We suggested using persistence diagrams to choose the optimal value for $\delta$ and illustrated it with a few examples.

Once the \deltaConnected\ components are introduced and calculated, we aim to determine a hierarchy between the different components. Therefore, we propose to use specific affinity matrices. They describe a way to group components based on their proximity and local alignment. To account for subtleties in the data, we include a data term in the initialization of the affinity matrices. The full algorithm can be found in \ref{alg:affinityMatrices}.

To show the performance of the \deltaConnected\ component algorithm and the affinity matrices algorithm, we have tested both algorithms on several 2D images of the retina. All experiments show that the \deltaConnected\ component algorithm can distinguish different structures at crossings. Additionally, the \deltaConnected\ components group well-aligned structures, resulting in more complete vessels in the output. However, it cannot always identify full vascular trees. Therefore, we calculated the affinity matrices for several results of the \deltaConnected\ component algorithm. We see that this leads to more complete vascular trees, where different \deltaConnected\ components are grouped that belong to the same vascular tree.

The algorithms are challenged by gaps at highly bending parts in the vascular structure. The algorithm cannot connect these structures without connecting different vascular trees at crossings. Large spatial gaps also form challenges: choosing $\delta$ too high results in connecting vessel segments belonging to different vascular trees.

For future work, it would be interesting to train the metric tensor weights $(w_1,w_2,w_3)$ via a PDE-G-CNN to improve results further, and to use optimal transport to determine a hierarchical structure on different \deltaConnected\ components. In this article, we have done this with the affinity matrices. However, the application is also very suitable for creating a model relying on optimal transport in the Lie group $G$, determining which vascular structures are close and well-aligned, using existing models like \cite{bonPai2025optimal}, or by creating new variants that are more suitable for the application at hand.
Additionally, aiming for faster GPU implementations via TaiChi (as done for geodesic tracking in \cite{berg_sherry2024SO3}) by fully using the parallelization options would significantly improve the runtime.

%% file: Appendices.tex
\appendix
\section{Logarithmic Norm Approximation}\label{app:20230922:LogarithmicNormApproximation}
Calculating the exact distances in a Lie group $G$ is expensive. Therefore, we approximate the norm with a so-called `logarithmic norm approximation'. This section elaborates on how this is done, as explained in \cite{Smets2023PDEBased,Bellaard2023analysis}.

Let us consider any Lie group $G$, equipped with a left-invariant metric tensor field $\mathcal{G}$, Lie group elements $g\in G$ and identity element $e\in G$. We consider an exponential curve $\gamma:[0,1]\to G$, connecting $e$ to $g$ where $\|\dot{\gamma}(0)\|$ is minimal, i.e. $\gamma(0)=e$, $\gamma(1)=g$ with the property $\gamma(s+t)=\gamma(s)\cdot\gamma(t)$.
We want to be able to calculate distances between different elements in the Lie group, e.g. the distance between $e$ and $g$, denoted by $d_{\mathcal{G}}(e,g)$. This is quite an expensive operation. Therefore, we approximate the exact distance with the distance of a logarithm, as explained in this section:
\begin{subequations}
    \begin{align}
        d_{\mathcal{G}}(e,g)&\leq \text{Length}_{\mathcal{G}}(\gamma)\\
        &=\int_0^1\left\|\dot{\gamma}(t)\right\|_{\mathcal{G}}\mathrm{d}t\\
        &=\int_0^1\left\|\left(L_{\dot{\gamma}(t)}\right)_*\dot{\gamma}(0)\right\|_{\mathcal{G}}\mathrm{d}t\\
        &=\int_0^1 \left\|\dot{\gamma}(0)\right\|_{\mathcal{G}}\mathrm{d}t\label{eq:20230922:gammadotteqgammadot0}\\
        &=\left\|\dot{\gamma}(0)\right\|_{\mathcal{G}}\label{eq:20230922:gammadot0teqgammadot0}\\
        &=\left\|\log g\right\|_{\mathcal{G}}\label{eq:20230922:gammadot0eqlogg},
    \end{align}
\end{subequations}
where the equality in \eqref{eq:20230922:gammadot0eqlogg} holds because by definition $\log g\coloneqq\dot{\gamma}(0)\in T_e(G)$, \eqref{eq:20230922:gammadot0teqgammadot0} holds because the metric tensor field is left-invariant, and where \eqref{eq:20230922:gammadotteqgammadot0} holds because 
\begin{align*}
    \dot{\gamma}(s)&=\left.\frac{\mathrm{d}}{\mathrm{d}t}\gamma(s+t)\right|_{t=0}
    =\left.\frac{\mathrm{d}}{\mathrm{d}t}\left(\gamma(s)\cdot\gamma(t)\right)\right|_{t=0}
    \\&
    =\left.\left(\frac{\mathrm{d}}{\mathrm{d}t}L_{\gamma(s)}\gamma(t)\right)\right|_{t=0}
    =\left(L_{\gamma(s)}\right)_*\dot{\gamma}(0),
\end{align*}
with $L_g:G\to G$ the left action defined by $L_g h= g\cdot h$ for all $h\in G$, and where the last equality follows directly by the definition of the push-forward.
Hence, we have found an upper bound for the distance from $e\in G$ to a point $g\in G$ in any Lie group $G$.

\begin{remark}
    By restricting the space of curves over which we optimize in the Riemannian distance $d_{\mathcal{G}}$ (cf. Eq. \eqref{eq:distance}) to the exponential curves (which have frozen coefficients in the left-invariant frame), we obtain $\rho_{\mathcal{G}}$. Essentially, the Riemannian exponential and Riemannian logarithm are then approximated by the Lie group exponential and the Lie group logarithm. Note that $d_{\mathcal{G}}(g,e)\leq \rho_{\mathcal{G}}(g)$ because we restrict ourselves to exponential curves in the optimization.
\end{remark}

In this article, we consider two Lie groups; the special Euclidean group $G=SE(2)$ and the special orthogonal group $G=SO(3)$. For both, we will calculate $\log g$ for $g\in G$. 
Let $R_G(g)$ denote the matrix representation of the element $g\in G$, such that $R_G(g_1)R_G(g_2)=R_G(g_1 g_2)$. Then, the associated Lie algebra is denoted by $\text{span}(A_i)$, where
\begin{align}
    A_i=\left.\frac{\partial R_G(g)}{\partial g^i}\right|_e,\label{eq:LieAlgebraGenerator}
\end{align}
with $e\in G$ the identity element and where $\left\{\mathcal{A}_i\right\}$ denotes the left-invariant frame, with dual $\left\{\omega^i\right\}$ such that $\langle\omega^i,\mathcal{A}_j\rangle=\delta_j^i$. Then, \begin{align*}
    \sum_{i=1}^n g^i A_i\leftrightarrow g\in G.
\end{align*}

\subsection{Logarthmic Norm Approximation in \SE{2}}
We start with calculating $\log g$ for the special Euclidean group $G=SE(2)$, where $g=(x,y,\theta)\in G$. The generating matrix $R_{SE(2)}$ is given by
\begin{align*}
    R_{SE(2)}(x,y,\theta)=\begin{pmatrix}
        \cos\theta&-\sin\theta&x\\
        \sin\theta&\cos\theta&y\\
        0&0&1
    \end{pmatrix}.
\end{align*}
Then, the associated Lie algebra $se(2)=\text{span}(A_1,A_2,A_3)$ is spanned by
\begin{align*}
    A_1&=\begin{pmatrix}
        0&0&1\\
        0&0&0\\
        0&0&0
    \end{pmatrix},&A_2&=\begin{pmatrix}
        0&0&0\\
        0&0&1\\
        0&0&0
    \end{pmatrix},&A_3&=\begin{pmatrix}
        0&-1&0\\
        1&0&0\\
        0&0&0
    \end{pmatrix},
\end{align*}
calculated by \eqref{eq:LieAlgebraGenerator} using identity element $e=(0,0,0)$. Hence, 
\begin{align*}
    &\SE{2}\ni(x,y,\theta)\leftrightarrow\\
    &\qquad\begin{pmatrix}
        \cos\theta&-\sin\theta&x\\
        \sin\theta&\cos\theta&y\\
        0&0&1
    \end{pmatrix}=\exp\left(x A_1\right)\exp\left(y A_2\right)\exp\left(\theta A_3\right)\\
    &\qquad\qquad=\exp\left(c^1 A_1+c^2 A_2+c^3 A_3\right)\\
    &\qquad\qquad=\begin{pmatrix}
        \cos c^3&-\sin c^3&\left(c^1\cos \tfrac{c^3}{2}-c^2\sin \tfrac{c^3}{2}\right)\text{sinc }\tfrac{c^3}{2}\\
        \sin c^3&\cos c^3&\left(c^1\sin \tfrac{c^3}{2}+c^2\cos\tfrac{c^3}{2}\right)\text{sinc }\tfrac{c^3}{2}\\
        0&0&1\\
    \end{pmatrix}.
\end{align*}

Hence, one has the following relations:
\begin{align}
    \begin{cases}
        x=\left(c^1\cos \tfrac{c^3}{2}-c^2\sin \tfrac{c^3}{2}\right)\text{sinc }\tfrac{c^3}{2}\\
        y=\left(c^1\sin \tfrac{c^3}{2}+c^2\cos\tfrac{c^3}{2}\right)\text{sinc }\tfrac{c^3}{2}\\
        \theta=c^3\mod 2\pi
    \end{cases}\hspace{-10px},\quad\begin{cases}
        c^1=\frac{x\cos\frac{\theta}{2}+y\sin\frac{\theta}{2}}{\text{sinc }\frac{\theta}{2}}\\
        c^2=\frac{-x\sin\frac{\theta}{2}+y\cos\frac{\theta}{2}}{\text{sinc }\frac{\theta}{2}}\\
        c^3=\theta,
    \end{cases}\label{eq:logarithmicNormApproxSE2}
\end{align}
where $c^1\left.\partial_x\right|_e+c^2\left.\partial_y\right|_e+c^3\left.\partial_\theta\right|_e\in T_e(SE(2))$ represent the logarithmic coordinates.

\subsection{Logarthmic Norm Approximation in $SO(3)$}
We follow the same approach as above for the special orthogonal group $G=SO(3)$, with group elements $(\alpha,\beta,\varphi)\in SO(3)$. The generating matrix $R_{SO(3)}$ is given by
\begin{align*}
    R_{SO(3)}&(\alpha,\beta,\varphi)\\&=\begin{pmatrix}
        \cosalpha \cosbeta & \sinbeta \cosvarphi-\sinalpha \cosbeta \sinvarphi & -\sinalpha \cosbeta \cosvarphi-\sinbeta \sinvarphi \\
        -\cosalpha \sinbeta & \sinalpha \sinbeta \sinvarphi+\cosbeta \cosvarphi & \sinalpha \sinbeta \cosvarphi-\cosbeta \sinvarphi \\
        \sinalpha & \cosalpha \sinvarphi & \cosalpha \cosvarphi
    \end{pmatrix},
\end{align*}
where $\cosalpha\!=\cos\alpha$, $\cosbeta\!=\cos\beta$, $\cosvarphi\!=\cos\varphi$, $\sinalpha\!=\sin\alpha$, $\sinbeta\!=\sin\beta$ and $\sinvarphi\!=\sin\varphi$.
The associated Lie algebra $so(3)$ is spanned by
\begin{align*}
    A_1&=\begin{pmatrix}
        0&0&0\\
        0&0&-1\\
        0&1&0
    \end{pmatrix},&A_2&=\begin{pmatrix}
        0&0&1\\
        0&0&0\\
        -1&0&0
    \end{pmatrix},&A_3&=\begin{pmatrix}
        0&-1&0\\
        1&0&0\\
        0&0&0
    \end{pmatrix},
\end{align*}
calculated by \eqref{eq:LieAlgebraGenerator} using identity element $e=(0,0,0)$. Hence, 
\begin{align*}
    SO(3)\ni (\alpha,\beta,\varphi)&\leftrightarrow R_{SO(3)}(\alpha,\beta,\varphi)\\
    &=
    \exp(-\beta A_3) \exp(-\alpha A_2) \exp(\varphi A_1)    \\
    &=\exp(c^1 A_1+c^2 A_2+c^3 A_3).\numberthis\label{eq:20230922:SO3conditionLogNormApprox}
\end{align*}
We need to find the relation between $(\alpha,\beta,\varphi)$ and $(c^1,c^2,c^3)$ such that Eq.~\eqref{eq:20230922:SO3conditionLogNormApprox} holds. To this end we note that 
%
%
we can interpret the rotation described by $R_{SO(3)}$ as a counterclockwise rotation around an axis $\mathbf{a}$ with an angle $\phi$. Both can be identified from the matrix representation of the group elements of $SO(3)$ given by $R_{SO(3)}$. The rotation axis, notated by the vector $\mathbf{a}$, is the eigenvector corresponding to the eigenvalue 1. 
If the rotation matrix $R_{SO(3)}$ is not symmetric, the eigenvector $\mathbf{a}$ and the corresponding rotation angle $\phi$ are given by
\begin{align}
    \mathbf{a}=\frac{1}{2\sin\phi} \begin{pmatrix}
    R_{32}-R_{23}\\ R_{13}-R_{31}\\ R_{21}-R_{12}
\end{pmatrix}, \quad 
    \phi=\arccos\left(\frac{Tr(R_{SO(3)})-1}{2}\right).\label{eq:20230922:rotationAngle}
\end{align}
The expression for the rotation angle $\phi$ is readily verified, as for a rotation matrix $\tilde{R}$ rotating a vector with an angle $\tilde\phi$ around the $z$-axis. Then, using the conjugation invariance of the trace operator and the fact that there exists a matrix $P$ such that $P^{-1}\tilde{R}P=R_{SO(3)}$, Eq.~\eqref{eq:20230922:rotationAngle} holds for the rotation with angle $\phi$ around the axis $\mathbf{a}$.
Then, one has, by Rodrigues' rotation formula and the well-known Lie group exp formula in $SO(3)$ that
\begin{align*}
\exp&(\phi (a^{1}A_1+ a^2 A_2 + a^3 A_3))(\mathbf{x})=
    \exp(\phi\cdot \mathbf{a})(\mathbf{x}) \\&=\mathbf{x}\cos\phi+\left(\mathbf{a}\times\mathbf{x}\right)\sin\phi+\mathbf{a}\left(\mathbf{a}\cdot\mathbf{x}\right)\left(1-\cos\phi\right)=R_{SO(3)}(\mathbf{x}).
\end{align*}
Therefore, we can now conclude that the logarithmic coordinates in $SO(3)$ are given by 
\begin{align*}
   &(c^1,c^2,c^3)^T= \phi\cdot \mathbf{a}=\frac{1}{2\;\text{sinc }\phi}\begin{pmatrix}
    R_{32}-R_{23}\\ R_{13}-R_{31}\\ R_{21}-R_{12}
\end{pmatrix} \\
\Leftrightarrow \;&
\Omega_\mathbf{a}= 
a^1 A_1+ a^2 A_2 + a^3 A_3= \frac{1}{2 \, \textrm{sinc}\, \phi} (R-R^T),\numberthis\label{eq:logarithmicNormApproxSO3}
\end{align*}
where $\Omega_{\mathbf{a}} \in so(3)=T_{e}(SO(3))$ is the antisymmetric matrix associated to
$\Omega_{\mathbf{a}}(\mathbf{x})=\mathbf{a} \times \mathbf{x}$.
The formula above has an (innocent) removable singularity at $\phi=\pm \pi$ where the matrix $R$ becomes symmetric (with 1,-1,-1 on the diagonal after diagonalisation). At $\phi=0$ we have $\textrm{sinc}(0)=1$.  

\section{Reflectional Symmetries}\label{app:ReflectionalSymmetries}
Note that by \eqref{eq:logapprox}, we have left-invariance $d_\mathcal{G}(g,h)=d_\mathcal{G}(h^{-1}g,e)\approx\|\log h^{-1}g\|_{\mathcal{G}}$. Furthermore, this distance's (reflectional) symmetries carry over to the \deltaConnected\ component algorithm as we will show next.
\begin{lemma}\label{lemma:normInvarianceReflections}
    Set $q=h^{-1}g$. Then $\log q=\sum_{i=1}^n c^i(q)\left.\mathcal{A}_i\right|_e\in T_e(G)$. Now its norm $\|\log q\|_{\mathcal{G}}$ is invariant under $c^i\mapsto \pm c^i$, $i=1,\ldots,n$, which gives rise to $2^n$ symmetries.
\end{lemma}
\begin{proof}
    We recall from App.~\ref{app:20230922:LogarithmicNormApproximation} that $\|\log q\|_{\mathcal{G}}=\sqrt{\sum_{i=1}^n w_i (c^i)^2}$. Then, its invariance under $c^i\mapsto \pm c^i$ follows immediately. \qed
\end{proof}

\begin{figure}[ht!]
    \centering
    \begin{subfigure}[b]{0.45\linewidth}
        \includegraphics[width=\linewidth]{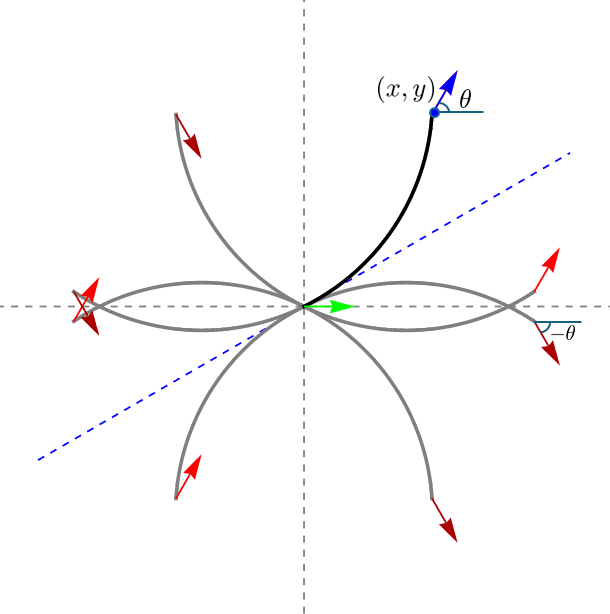}
        \caption{Reflectional symmetries with $c^2$ unconstrained (general Riemannian case \cite{Bellaard2023analysis}).}
        \label{fig:ReflectionalSymmetriesRiemann}
    \end{subfigure}\hfill
    \begin{subfigure}[b]{0.45\linewidth}
        \includegraphics[width=\linewidth]{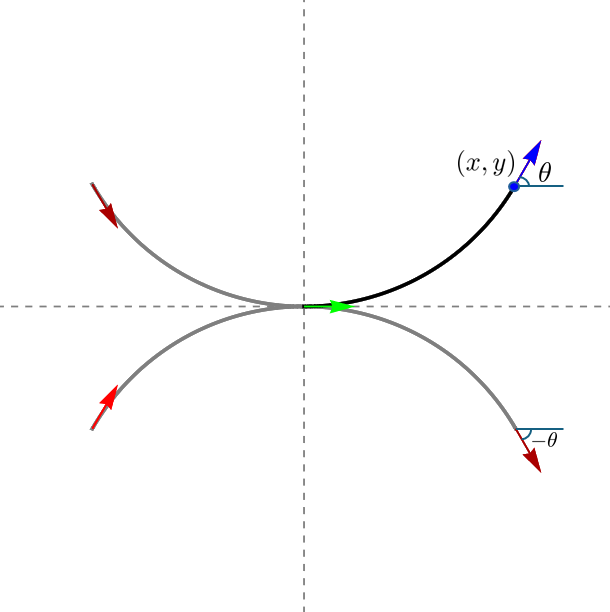}
        \caption{Reflectional symmetries with $c^2=0$ (sub-Riemannian case \cite{moiseev2010maxwell}.}
    \end{subfigure}    
    \caption{Reflectional symmetries of the blue point $(x,y,\theta)\equiv (x,y,R_\theta)\in SE(2)$, reflected in the three reflection axes. Depending on the axis, the orientation is updated to $-\theta$ or stays the same. The blue reflection axis is only for the endpoints of the geodesics.}
    \label{fig:ReflectionalSymmetries}
\end{figure}

See Fig.~\ref{fig:ReflectionalSymmetries} for what this means on SE(2) where $n=\dim SE(2)=3$. For further details, see \cite{Bellaard2023analysis,moiseev2010maxwell}.

\begin{corollary}\label{cor:ReflectionalSymmetries}
    If the connected component algorithm is applied to a point cloud $\{q_1,\ldots,q_n\}\subset SE(2)$, then it produces the same reflected connected components as when it is applied to $\{\tilde{\varepsilon}^i(q_1),\ldots, \tilde{\varepsilon}^i(q_n)\}$, where $\tilde{\varepsilon}^i(q)=\exp(\varepsilon^i \log(q))$, $i=0,\ldots,2^n-1$, denotes the symmetry given by the reflection of the logarithmic coordinates as denoted in \cite[Table 5]{Bellaard2023analysis}.
\end{corollary}
\begin{proof}
    The statement follows immediately from Lemma~\ref{lemma:normInvarianceReflections} applied to $\|\log\tilde{\varepsilon}^i(g)$:
    \begin{align*}
        \|\log\tilde{\varepsilon}^i(q)\|_\mathcal{G}=\|\varepsilon^i (\log q)\|_{\mathcal{G}}\overset{\text{Lemma~}\ref{lemma:normInvarianceReflections}}{=} \|\log q\|_{\mathcal{G}}.
    \end{align*}
    Thus, the distances between elements are the same, and consequently, the \deltaConnected\ components are also the same. \qed
\end{proof} 

The precise formulas of these reflections can be found in \cite{moiseev2010maxwell,Bellaard2023analysis}. For a quick intuition, see Fig.~\ref{fig:ReflectionalSymmetries}. Essentially, they are found by line reflections in the half-angle axis, the $x$- and $y$-axes, possibly negating the angle. Intuitively, they are depicted in Fig.~\ref{fig:ReflectionalSymmetriesRiemann}.

\section{Proofs of different Lemma's and Propositions}\label{app:Proofs}
\input{AppendixFiles}

\subsection{Proof of Proposition~\ref{th:20230605:Un}}\label{app:20230922:ProofPropUn}
\begin{proof}
    First of all, note that $m_\delta(g,g_0)\leq n$ implies that $d(g,g_0)\leq n\delta$ by applying the triangular inequality to the sequence of intermediate points in Def.~\ref{def:20230605:Connectedness}.
    
    Then, note that the inequalities follow directly from Lemma~\ref{lemma:propertiesMorphologicalDilations}, Eq.~\eqref{eq:morphConvRange}. We will prove the rest of this statement by induction w.r.t. $n$. Without loss of generality in view of left-invariance, we prove the statement for $g_0=e$. This is analogous to doing a roto-translation with $g_0^{-1}$ of the given binary map $\B$, and afterwards pushing forward the output $\U{e}{n}{\cdot}$ to $\U{g_0}{n}{\cdot}$. First, we show that the statement is true for the first step of the algorithm. We recall Eq.~\eqref{eq:morphKernelAlpha1}, which stated
    \begin{align*}
        \morphkern{t}{1}(g)\coloneqq\begin{cases}
            0 & \text{if }d(g,e)<t\\
            \infty & \text{else}.
        \end{cases}
    \end{align*}
    
    We can easily express the first step of the connected component algorithm $n=1$ by calculating $\U{e}{1}{\cdot}$:
    \begin{align*}
        \U{e}{1}{g}&=\left(-\left(\morphkern{\delta}{1}\square -\U{e}{0}{\cdot}\right)(g)\right)\B(g)\\
        &=\begin{cases}
            \sup\limits_{h\in G} \left\{\U{e}{0}{g}-\morphkern{\delta}{1}(h^{-1}g)\right\} & \text{if } g\in I,\\
            0 & \text{else};
        \end{cases}\\
        &=\begin{cases}
            1
            & \text{if }g\in I \text{ and }d(e,g)<\delta
            \\
            0 & \text{else},
        \end{cases}
    \end{align*}    
    where $g\in I$ and $d(e,g)<\delta$ is equivalent to $g\simdelta e$ and $m_\delta(e,g)\leq 1$.
    Note that $g\simdelta e$ automatically implies $g\in I$.
    
    Next, we assume that the statement is true for a given $n\in\mathbb{N}$, meaning that one has \begin{align*}
        \U{e}{n}{g}&=\begin{cases}
            1 &\text{if }g\simdelta e \wedge m_\delta(e,g)\leq n\\
            0&\text{else}.
        \end{cases}\\
        &=\mathbbm{1}_{\setIendelta{e}{n}{\delta}}(g),
    \end{align*}
    where $$\setIendelta{e}{n}{\delta}=[e]\cap \tildeB{n}{\delta}{e},$$ where $[e]$ denotes the set of points that are \deltaConnected\ to $e$ and where $\tildeB{n}{\delta}{e}$ is defined by
    \begin{align}
        \tildeB{n}{\delta}{e}\coloneqq \left\{g\in I\;\left|\; m_\delta(e,g)\leq n\right.\right\}.
    \end{align} 
    Then, we aim to show that the statement is also true for iteration $n+1$. That means that one investigates whether one can write $\U{e}{n+1}{\cdot}$ in the same form as $\U{e}{n}{\cdot}$. Let $g\in I$, then
    \begin{align*}
        \U{e}{n+1}{g}&=
        -\left(-\morphkern{\delta}{1}\square -\U{e}{n}{\cdot}\right)(g)
        \\
        &=
        -\left(-\morphkern{\delta}{1}\square -\mathbbm{1}_{\setIendelta{e}{n}{\delta}}\right)(g)
        \\
        &=
        \sup_{h\in G}\left\{\mathbbm{1}_{\setIendelta{e}{n}{\delta}}(h)-\morphkern{\delta}{1}(h^{-1}g)\right\}
        .
        \numberthis\label{eq:20230822CleanedUpVersion:Un+1indicatorsets}
    \end{align*}
    Then, in order to reach the supremum in \eqref{eq:20230822CleanedUpVersion:Un+1indicatorsets}, we can simplify the expression to
    \begin{align*}
        \U{e}{n+1}{g}&=\begin{cases}
            1&\text{if }g\in \setIendelta{e}{n}{\delta}\text{ or if }\exists h\in \setIendelta{e}{n}{\delta}\text{ s.t. }d(h,g)<\delta\\
            0&\text{else}.
        \end{cases}
    \end{align*}
    We will reformulate the condition for $\U{e}{n+1}{g}$ to equal 1. First, it is important to note that $\setIendelta{e}{n}{\delta}\subseteq \setIendelta{e}{n+1}{\delta}$. Let us assume $g\notin \setIendelta{e}{n}{\delta}$, but there exists a $h\in \setIendelta{e}{n}{\delta}$ such that $d(h,g)<\delta$. Then, the first condition gives us
    \begin{align}
        m_\delta(e,g)>n.\label{eq:Mdeltageqn}
    \end{align}
    Simultaneously, there exists a $h\in \setIendelta{e}{n}{\delta}$ such that $d(h,g)<\delta$, so
    \begin{align}
        m_\delta(e,h)\leq n \quad \text{and} \quad d(h,g)<\delta, \quad \text{so} \quad m_\delta(e,g)\leq n+1.\label{eq:Mdeltaleqn+1}
    \end{align}
    Combining \eqref{eq:Mdeltageqn} and \eqref{eq:Mdeltaleqn+1} tells us that $m_\delta(e,g)=n+1$. Additionally, the value of $\U{e}{n+1}{g}$ is equal to 1 when $g\in \setIendelta{e}{n}{\delta}$, i.e., $m_\delta(e,g)\leq n$. Hence, $\U{e}{n+1}{g}$ can be reformulated to
    \begin{align*}
        \U{e}{n+1}{g}&=\begin{cases}
            1&\text{if }g\in \setIendelta{e}{n+1}{\delta}\text{, i.e. } m_\delta(e,g)\leq n+1 \text{ and }g\simdelta e\\
            0&\text{else}.
        \end{cases}
    \end{align*}
    This is the same expression as was assumed to be true for some $n\in\mathbb{N}$ but now for $n+1\in\mathbb{N}$, and hence, the statement is proven.

    \qed
\end{proof}

\section{For \texorpdfstring{$\alpha>1$}{a>1} the equivalence relation breaks }\label{app:symmetry}
For an example where  the equivalence relation \texorpdfstring{$\simdelta$}{ } 
breaks down for $\alpha>1$; see  Fig.~\ref{fig:20230922:Symmetry}. It motivates why in connected component algorithms we choose $\alpha \downarrow 1$ (or $\alpha=1$) and why we used the evolutions for $\alpha>1$ only for the affinity measure experiments (where $\alpha>1$ shrinks the wavefront propagation and softens the max-pooling). 
Some theoretical results on viscosity solutions for HJB equations \cite{Bellaard2023analysis,azagra2005nonsmooth,diop2021extension,fathi2007weakkam}, do require $\alpha>1$. In the algorithms for the experiments (e.g.~\ref{alg:deltaconnectedComponents}) one can either take $\alpha=1$ or $\alpha \downarrow 1$ which boils down to the same when taking the pointwise limit in the morphological convolutions (\ref{eq:viscositySolutionW}) that solve the HJB PDE system (\ref{eq:20230922:DilationPDE}) almost everywhere. 
\begin{figure}
    \centering
    \begin{subfigure}[b]{0.45\textwidth}
         \centering
        \scalebox{.8}{\input {Figures/SymmetryOfConnectedComponents/SymmetryInConnectedComponentsSetting}}
        \caption{Example setting of one connected component $[h]=\bigcup\limits_{i=1}^3 h_i$ with $d(h_1,h_2)=d(h_2,h_3)=\delta$ and $d(h_1,h_3)=2\delta$.}
        \label{fig:20230922:SettingSymmetry}
    \end{subfigure}\hfill
    \begin{subfigure}[b]{0.45\textwidth}
         \centering
        \scalebox{.82}{\input {Figures/SymmetryOfConnectedComponents/SymmetryInConnectedComponentsStarth2}}
        \begin{picture}(0,0)(-55,-70)
            \framebox{\shortstack[l]{According to \\
            \ref{alg:deltaconnectedComponents}:\\$h_1\simdelta h_2\simdelta h_3$}}
        \end{picture}
        \caption{Identification of \deltaConnected\ component with the algorithm starting at $g_0=h_2$, identifying the whole component in 1 step, independent of the value of $\alpha$.}
        \label{fig:20230922:Symmetryh2}
    \end{subfigure}\hfill
    \begin{subfigure}[b]{0.45\textwidth}
         \centering
        \scalebox{.85}{\input {Figures/SymmetryOfConnectedComponents/SymmetryInConnectedComponentsStarth1}}
        \begin{picture}(0,0)(-58,-90)
            \framebox{\shortstack[l]{According to \\
            \ref{alg:deltaconnectedComponents}:\\$h_1\simdelta h_2$\\ $h_2\simdelta h_3$\\$h_1\!\!\not\simdelta\! h_3$}}
        \end{picture}
        \caption{Identification of \deltaConnected\ component with $\alpha=2$ and the algorithm starting at $g_0=h_1$ converging after 1 step, but not identifying $h_3$ as part of the \deltaConnected\ component $[h_1]$.}
        \label{fig:20230922:Symmetryh1}
    \end{subfigure}
    \caption{The \deltaConnected\ component $[h]$ (in red in Fig.~\ref{fig:20230922:SettingSymmetry}) is not always completely identified with \ref{alg:deltaconnectedComponents} when not choosing $\alpha=1$. In the example, we have used $\alpha=2$ to show that the algorithm returns different \deltaConnected\ components for different starting points that may differ from the full \deltaConnected\ component. The example uses starting points $h_2$ and $h_1$ for results in Fig.~\ref{fig:20230922:Symmetryh2}~and ~\ref{fig:20230922:Symmetryh1} respectively.
    }
    \label{fig:20230922:Symmetry}
\end{figure}

%% file: AppendixFiles.tex
\subsection{Multiple consecutive morphological dilations}\label{sec:multmorphdil}
We will show what happens when multiple consecutive morphological dilations are executed in Lemma~\ref{lemma:20230605CleanedUpVersion:KernelConcatenation}. In the statement and proof, we use the morphological delta which we define first.
\begin{definition}[Morphological delta]
    We define the morphological delta $\delta_{e}^M:G\to\mathbb{R}\cup\{\infty\}$ as
    \begin{align}
        \delta_{e}^M(g)=\begin{cases}
            0 &\text{if }g=e\\
            \infty &\text{else.}
        \end{cases}
    \end{align}
\end{definition}
\begin{lemma}\label{lemma:20230605CleanedUpVersion:KernelConcatenation}
    Let $1\leq\alpha\leq \infty$. Consider the morphological kernel $\morphkern{t}{\alpha}$. Let $g\in G$ be given. When defining\footnote{These pointwise limits arise from \eqref{eq:20230822:MorphologicalKernelNew}.} $k_0^\alpha=\delta_e^M$ for $\alpha<\infty$, and $k_0^\infty=d(\cdot,e)$, we have:
    \begin{align*}
        \forall t,s\geq 0: \;\left(\morphkern{t}{\alpha}\square \morphkern{s}{\alpha}\right)(g)=\morphkern{t+s}{\alpha}(g)=(\kappa_t^\alpha\square_{\mathbb{R}}\kappa_s^\alpha)(d(g,e)),
    \end{align*}
    where $\kappa_t^\alpha:\mathbb{R}^+\to\mathbb{R}^+$, is defined as $\kappa_t^\alpha(d(g,e))=\morphkern{t}{\alpha}(g)$.
\end{lemma}
\begin{proof}
    For $t=0$, $\alpha<\infty$, one has $\morphkern{t}{\alpha} = \delta_e^M$, and $\delta_e^M\square f = f$, so let us examine $1\leq\alpha<\infty$, $t,s>0$. Let $\kappa_t^\alpha(r)\coloneqq \kappa_t^\alpha(|r|)$ for all $r\in\mathbb{R}$.
    Then, one can write
    \begin{subequations}
        \begin{align}
            \left(\morphkern{t}{\alpha}\square \morphkern{s}{\alpha}\right)(g)&=\inf_{h\in G} \morphkern{t}{\alpha}(h^{-1}g)+\morphkern{s}{\alpha}(h)\\
            &=\inf_{h\in G} \kappa_t^{\alpha}(d(h^{-1}g,e))+\kappa_s^{\alpha}(d(h,e))\\
            &=\inf_{h\in \gamma_{e,g}^{\min}} \kappa_t^{\alpha}(d(g,e)-d(h,e))+\kappa_s^{\alpha}(d(h,e))\label{eq:Lemma31equality3}\\
            &=\inf_{0\leq r\leq d(g,e)} \kappa_t^{\alpha}(d(g,e)-r)+\kappa_s^{\alpha}(r)\\
            &=\inf_{r\in \mathbb{R}} \kappa_t^{\alpha}(d(g,e)-r)+\kappa_s^{\alpha}(r)\label{eq:Lemma31equality5}\\
            &=(\kappa_t^\alpha\square_{\mathbb{R}}\kappa_s^\alpha)(d(g,e))\label{eq:Lemma31equality6}\\
            &=\kappa_{t+s}^{\alpha}(d(g,e))=\morphkern{t+s}{\alpha}(g),\label{eq:Lemma31equality7}
        \end{align}
    \end{subequations}
    where the third equality \eqref{eq:Lemma31equality3} holds because one has $\geq$ due to the triangle inequality and $\kappa_t^{\alpha}$ being monotonous. The $\leq$ holds because $\gamma_{e,g}^{\min}\subset G$, where $\gamma_{e,g}^{\min}$ denotes the minimizing geodesic connecting $e$ and $g$. To show why the fifth equality \eqref{eq:Lemma31equality5} holds, we define $G(r)\coloneqq\kappa_t^{\alpha}(d(g,e)-r)+\kappa_s^{\alpha}(r)=\frac{t}{\beta}\left(\frac{d-r}{t}\right)^\beta+\frac{s}{\beta}\left(\frac{r}{s}\right)^\beta$, where we omit the dependence on the elements $g$ and $e$ in the notation for the distance $d=d(g,e)$. Then 
    \begin{align*}
        G'(r)=-\left(\frac{d-r}{t}\right)^{\beta-1}+\left(\frac{r}{s}\right)^{\beta-1}=0
        &\;\Leftrightarrow\; r=\frac{s}{t+s}d.
    \end{align*}
    Since $s,t>0$, $r$ attains a global minimimum between $0$ and distance $d=d(g,e)$. The sixth equality \eqref{eq:Lemma31equality6} is by definition of the morphological convolutions on $\mathbb{R}$. The seventh equality \eqref{eq:Lemma31equality7} follows by well-known results for morphological PDEs on $\mathbb{R}$ \cite{SchmidtWeickert2016Morphological}.

    Note that the case $\alpha=\infty$, $t=0$, $s>0$ is the same as $\alpha=\infty$, $t,s>0$ since $k_0^\infty(g)=k_t^\infty(g)$ for all $t\geq 0$, and where
    \begin{align*}
        (\morphkern{t}{\infty}\square\morphkern{s}{\infty})(g)=\inf_{h\in G} d(h^{-1}g,e)+d(h,e)=d(g,e).\text{\qed}
    \end{align*}
\end{proof}

%% file: Figures/SymmetryOfConnectedComponents/SymmetryInConnectedComponentsSetting.tex
\begin{tikzpicture}
    \draw[<->]({-1+0.06},0)--({1-0.06},0) node at (0,0)[anchor=south] {$\delta$};
    \draw[<->]({1+0.06},0)--({3-0.06},0) node at (2,0)[anchor=south] {$\delta$};
    \filldraw[color=red, fill=red, very thick](-1,0) circle (0.05) node[anchor=north] {$h_1$};
    \filldraw[color=red, fill=red, very thick](1,0) circle (0.05) node[anchor=north] {$h_2$};
    \filldraw[color=red, fill=red, very thick](3,0) circle (0.05) node[anchor=north] {$h_3$};
\end{tikzpicture}

%% file: Figures/SymmetryOfConnectedComponents/SymmetryInConnectedComponentsStarth2.tex
\begin{tikzpicture}
    \filldraw[color=red!40, fill=red!20, very thick](1,0) circle (2.06);
    \draw[<->]({-1+0.06},0)--({1-0.06},0) node at (0,0)[anchor=south] {$\delta$};
    \draw[<->]({1+0.06},0)--({3-0.06},0) node at (2,0)[anchor=south] {$\delta$};
    \filldraw[color=red, fill=red, very thick](-1,0) circle (0.05) node[anchor=north] {$h_1$};
    \filldraw[color=red, fill=red, very thick](1,0) circle (0.05) node[anchor=north] {$h_2$};
    \filldraw[color=red, fill=red, very thick](3,0) circle (0.05) node[anchor=north] {$h_3$};
    \draw[color=black] node at ({-1-2*sqrt(2)},{1}) [anchor=west] {$\alpha=2$};
    \draw[color=black] node at ({-1-2*sqrt(2)},{.5}) [anchor=west] {$g_0=h_2$};
    \draw[color=red!0] node at ({3+2*sqrt(2)+1},0)[anchor=south] {};
\end{tikzpicture}

%% file: Figures/SymmetryOfConnectedComponents/SymmetryInConnectedComponentsStarth1.tex
\begin{tikzpicture}
    \draw[color=red!0] node at ({3+2*sqrt(2)+1},0)[anchor=south] {};
    \filldraw[color=red!40, fill=red!20, very thick, dashed] (-1,0) circle ({2*sqrt(2)});
    \filldraw[color=red!40, fill=red!20, very thick](-1,0) circle (2.06);
    \filldraw[color=red!40, fill=red!20, very thick](1,0) circle ({2*sqrt(2)-2});
    \filldraw[color=red!40, fill=red!20, very thick, dashed, fill opacity = 0](-1,0) circle (2.06);
    \draw[<->]({-1+0.06},0)--({1-0.06},0) node at (0,0)[anchor=south] {$\delta$};
    \draw[<->]({-1-0.05},0.05)--({-1+2*sqrt(2)*cos(120)},{2*sqrt(2)*sin(120)}) node at ({-1+sqrt(2)*cos(120)},{sqrt(2)*sin(120)})[anchor=east] {$\delta\sqrt{2}$};
    \draw[<->]({1+0.06},0)--({3-0.06},0) node at (2,0)[anchor=south] {$\delta$};
    \filldraw[color=red, fill=red, very thick](-1,0) circle (0.05) node[anchor=north] {$h_1$};
    \filldraw[color=red, fill=red, very thick](1,0) circle (0.05) node[anchor=north] {$h_2$};
    \filldraw[color=red, fill=red, very thick](3,0) circle (0.05) node[anchor=north] {$h_3$};
    \draw[color=black] node at ({-1-2*sqrt(2)},{2*sqrt(2)-0.2}) [anchor=west] {$\alpha=2$};
    \draw[color=black] node at ({-1-2*sqrt(2)},{2*sqrt(2)-0.2-0.5}) [anchor=west] {$g_0=h_1$};
\end{tikzpicture}

%% file: Review_JMIV/AppendixSAM.tex
\section[Comparison of \deltaConnected\ Component Method to Existing Approaches]{Comparison of \deltaConnected\ Component Method to Existing Segmentation and Clustering Methods}\label{app:comparisonSAMToMATo}
\subsection{Segment Anything Model 2}\label{app:SAM}
The Segment Anything Model 2 trained on the SA-V dataset \cite{ravi2024sam2}, is an improvement of the Segment Anything Model \cite{kirillov2023segany}. The updated model can segment both videos and images.

\begin{figure*}[h!]
    \centering
    \begin{subfigure}[t]{0.32\linewidth}
        \centering
        \includegraphics[width=0.94\textwidth]{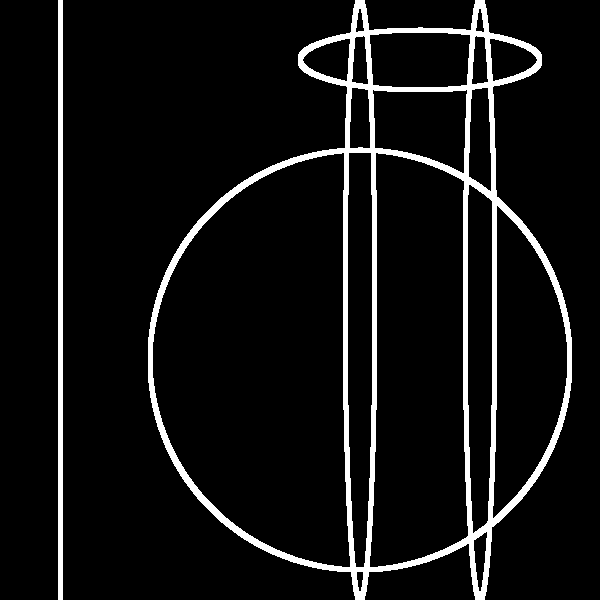}
        \caption{Reference image, cf.~Fig.~\ref{fig:20230922:comparisonR2vsSE2:chip}.}
        \label{fig:Appendix:Circles:reference}
    \end{subfigure}\hfill
    \begin{subfigure}[t]{0.32\linewidth}
        \centering
        \includegraphics[width=0.94\textwidth]{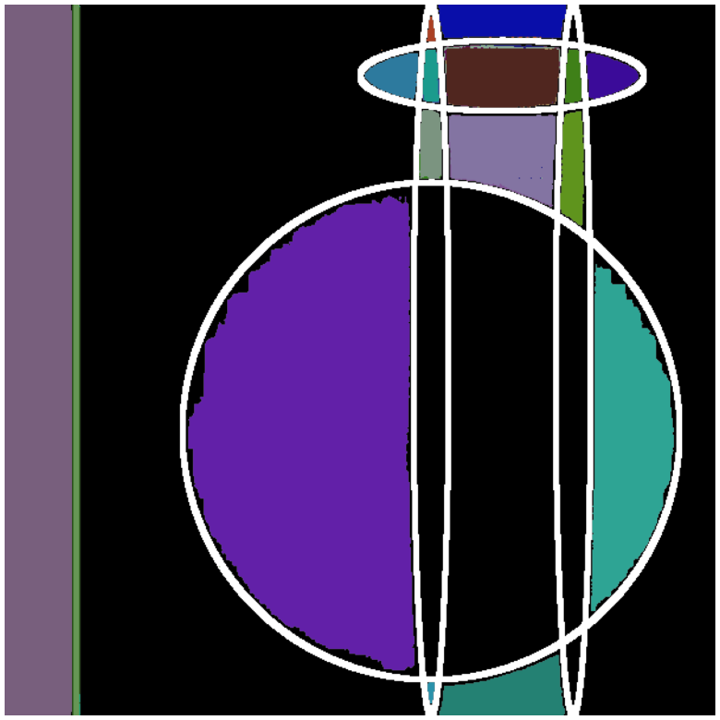}
        \caption{Automatic segmentation.}
        \label{fig:Appendix:Circles:autosegmentation}
    \end{subfigure}\hfill
    \begin{subfigure}[t]{0.32\linewidth}
        \centering
        \includegraphics[width=0.94\textwidth]{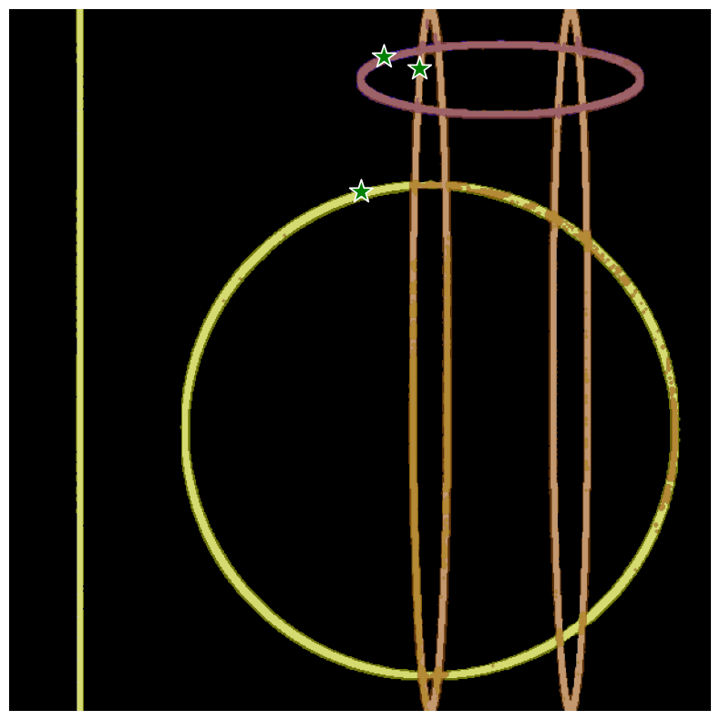}
        \caption{Segmentation.
        }
        \label{fig:Appendix:Circles:segmentation}
    \end{subfigure}        
    \caption[The output of the SAM 2 applied to an image containing lines and ovals.]{The output of the SAM 2 \cite{ravi2024sam2} applied to the reference image using to two different versions of the model; the auto masklet generator in Fig.~\ref{fig:Appendix:Circles:autosegmentation} and interactive segmentation that requires the (de)selection of pixels by the user in Fig.~\ref{fig:Appendix:Circles:segmentation}. The three stars were manually selected by the user (see Remark~\ref{remark:app:prediction}). See Table~\ref{tab:SAM:predictionScore} for the prediction scores of the segments in Fig.~\ref{fig:Appendix:Circles:segmentation}.}
    \label{fig:Appendix:Circles}
\end{figure*}

The model can be used in two different ways: 1) an automatic segmentation method that does not require user input, and 2) an interactive segmentation method that requires the user to select (or deselect) pixels that need to be included (or excluded) from the segmentation. Fig.~\ref{fig:Appendix:Circles} shows the output of the two different methods applied to the reference image in Fig.~\ref{fig:Appendix:Circles:reference}. The automatic method automatically samples its starting point. This results in segments of the background (black) and of the foreground (white). The density of the sample points can be adjusted. An increased number of sample points slows the model down, but increases the chance of correctly segmenting smaller features, such as thin blood vessels.
For example, in Fig.~\ref{fig:Appendix:Circles:autosegmentation}, the model missed the segment of the circles and ellipses. With the interactive method, the user selects the pixels that they want to be included in the segmentation.

Fig.~\ref{fig:Appendix:Circles:segmentation} shows the output of the three segmentations corresponding to the three selected pixels (visualized with a green star). Note that the line on the left is included in the segmentation of the large circle, the vertically oriented ellipses are grouped together, and there is overlap between the segmentations.

\begin{figure*} 
    \centering
    \begin{subfigure}[t]{0.32\linewidth}
        \centering
        \includegraphics[width=0.94\textwidth]{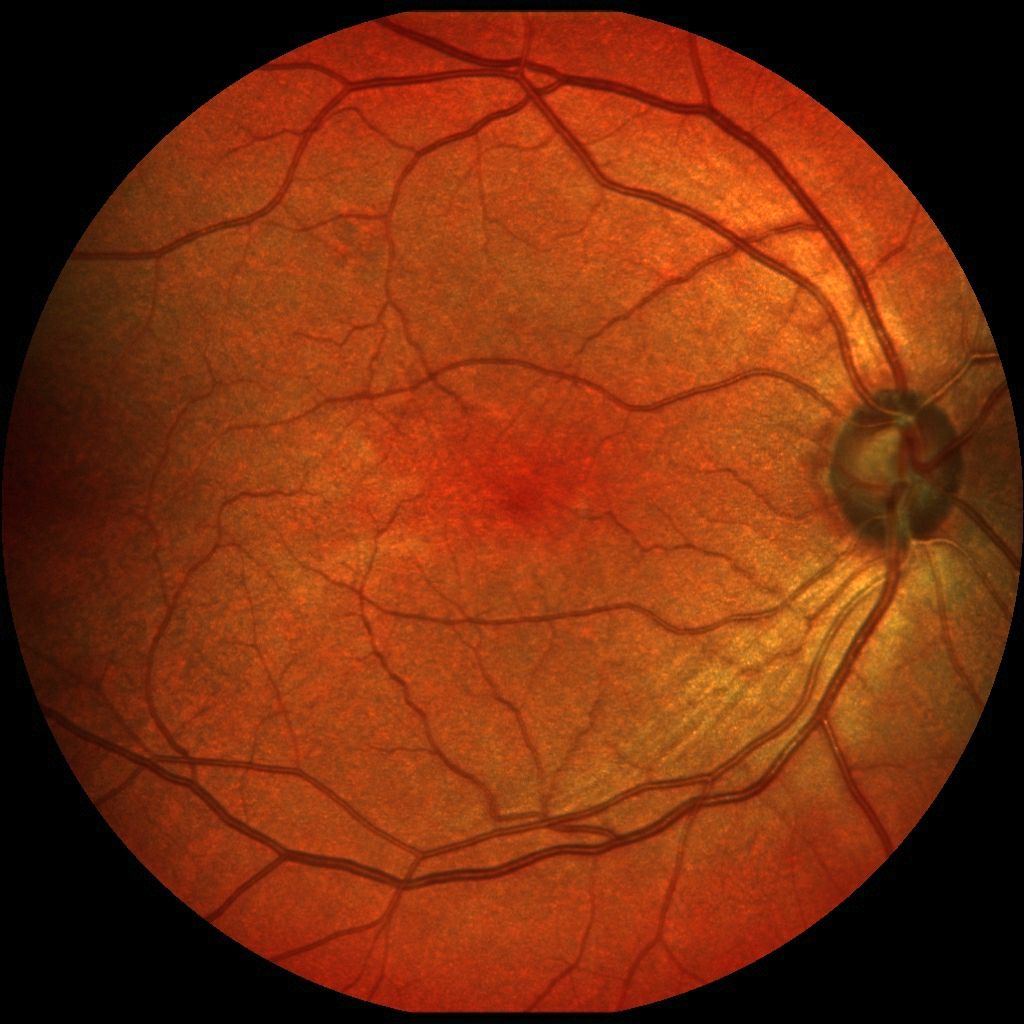}
        \caption{Reference image.}
    \end{subfigure}\hfill
    \begin{subfigure}[t]{0.32\linewidth}
        \centering
        \includegraphics[width=0.94\textwidth]{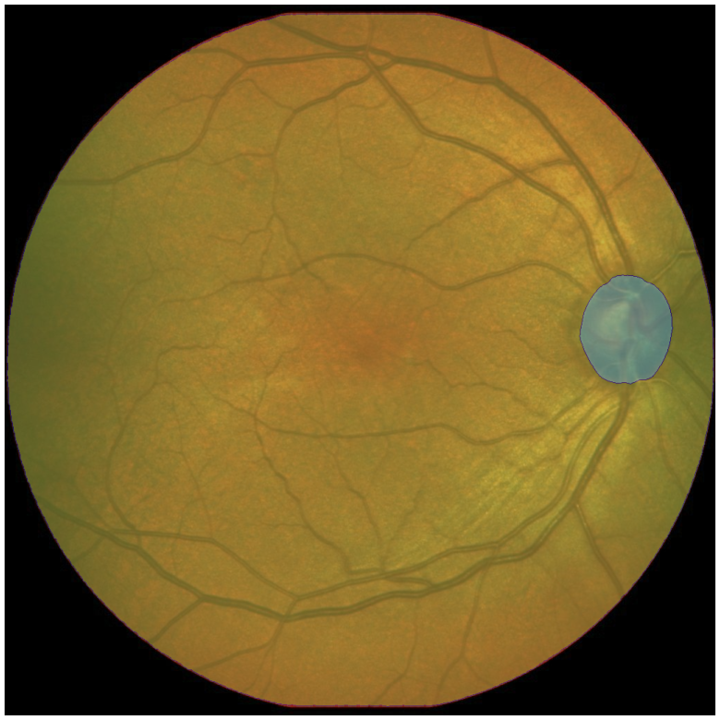}
        \caption{Automatic segmentation.}
        \label{fig:Appendix:STAR38:autosegmentation}
    \end{subfigure}\hfill
    \begin{subfigure}[t]{0.32\linewidth}
        \centering
        \includegraphics[width=0.94\textwidth]{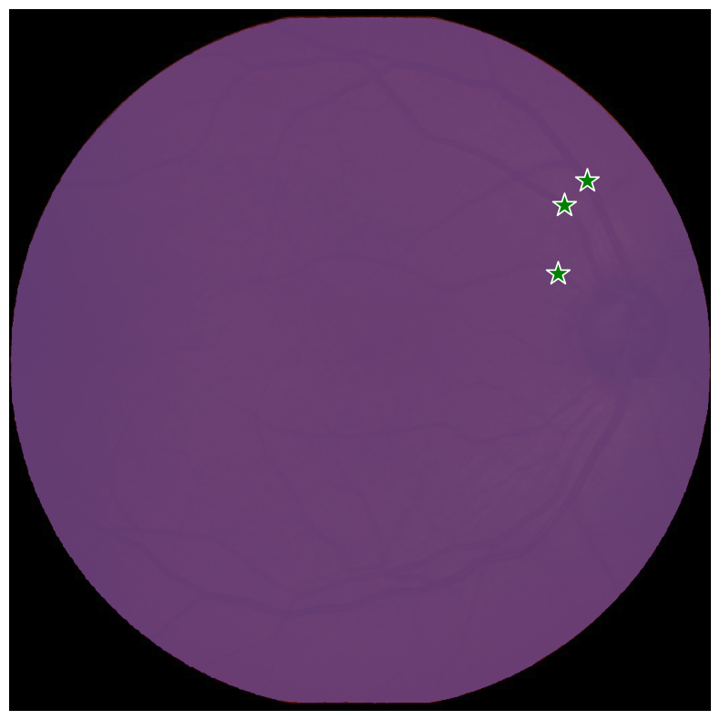}
        \caption{Segmentation.}
        \label{fig:Appendix:STAR38:segmentation}
    \end{subfigure}
    \caption[The output of the SAM 2 applied to a retinal image.]{The output of the SAM 2 \cite{ravi2024sam2} applied to a retinal image using to two different versions of the model; the auto masklet generator (cf.~Fig.~\ref{fig:Appendix:STAR38:autosegmentation}) and interactive segmentation that requires the (de)selection of pixels by the user in Fig.~\ref{fig:Appendix:STAR38:segmentation}. Three points (stars) were manually selected by the user (see Remark~\ref{remark:app:prediction}). Note that all three stars gave the same segmentation as output, visualized by the purple circle. See Table~\ref{tab:SAM:predictionScore} for the prediction scores of the segments in Fig.~\ref{fig:Appendix:STAR38:segmentation}.}
    \label{fig:Appendix:STAR38:original}
\end{figure*}

Fig.~\ref{fig:Appendix:STAR38:original} shows the output of the automatic and interactive model applied to a retinal image. Note that in both cases the model did not succeed in segmenting any blood vessels. We got the same results when trying the model on an enhanced version of the retinal image. 

\begin{figure*}
    \centering    
    \begin{subfigure}[t]{0.3\linewidth}
        \includegraphics[width=\textwidth]{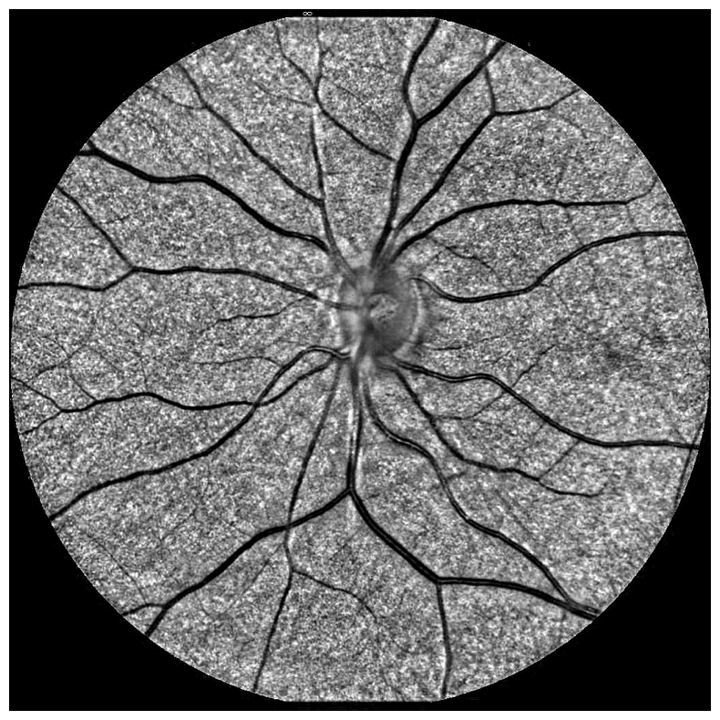}
        \caption{Enhanced retinal image.}
        \label{fig:Appendix:STAR48}
    \end{subfigure}\hfill
    \begin{subfigure}[t]{0.3\linewidth}
        \includegraphics[width=\textwidth]{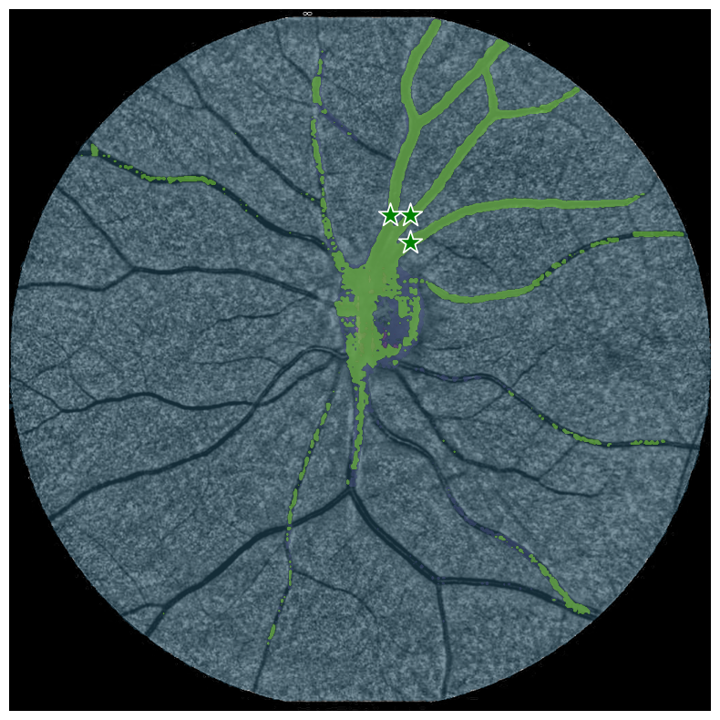}
        \caption{Best segmentation.}
        \label{fig:Appendix:STAR48:enhanced:segmentation1}
    \end{subfigure}\hfill
    \begin{subfigure}[t]{0.3\linewidth}
        \includegraphics[width=\textwidth]{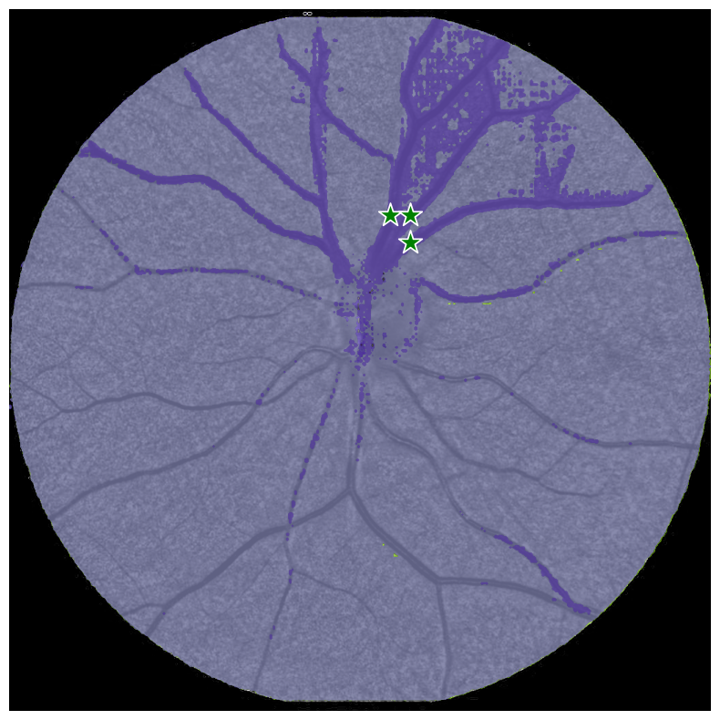}
        \caption{Third best segmentation.}
        \label{fig:Appendix:STAR48:enhanced:segmentation3}
    \end{subfigure}
    
    \centering    
    \begin{subfigure}[t]{0.3\linewidth}
        \includegraphics[width=\textwidth]{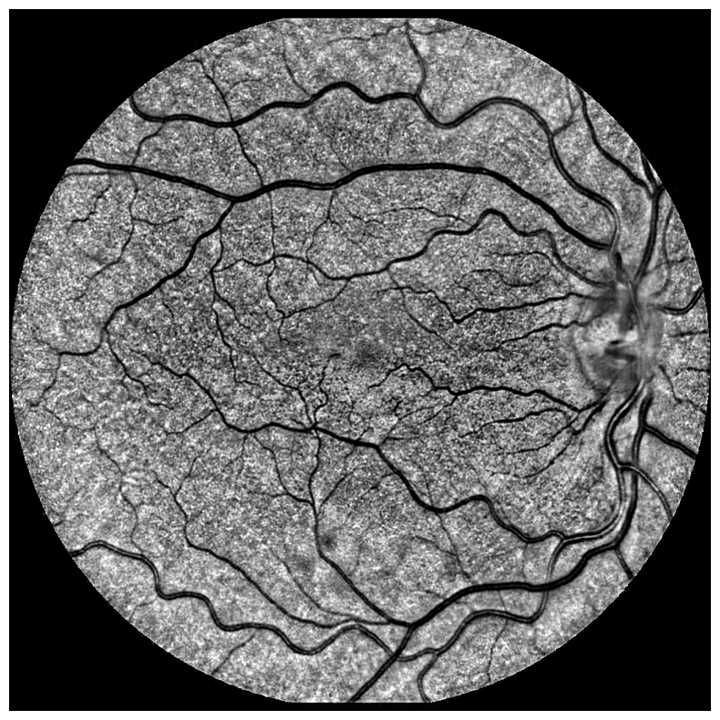}
        \caption{Enhanced retinal image.}
        \label{fig:Appendix:STAR34}
    \end{subfigure}\hfill
    \begin{subfigure}[t]{0.3\linewidth}
        \includegraphics[width=\textwidth]{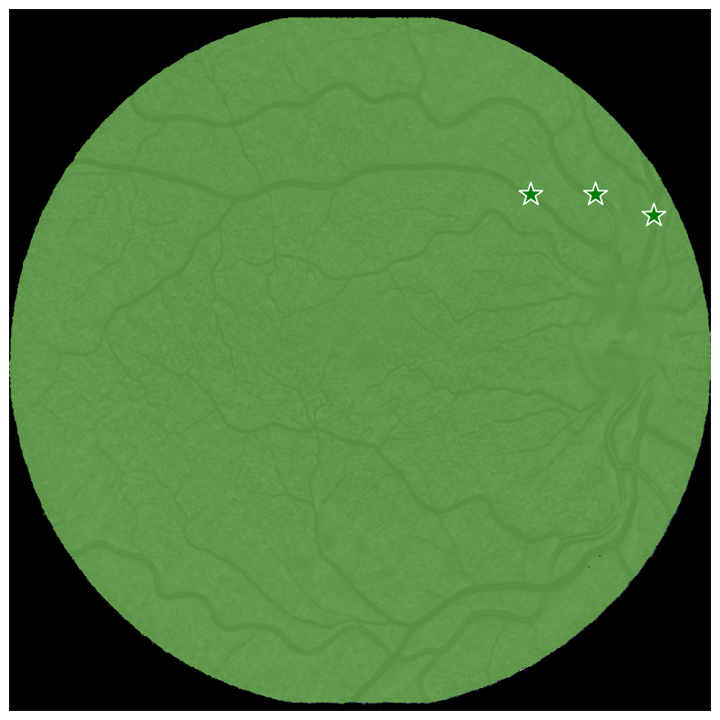}
        \caption{Best segmentation.}
        \label{fig:Appendix:STAR34:enhanced:segmentation1}
    \end{subfigure}\hfill
    \begin{subfigure}[t]{0.3\linewidth}
        \includegraphics[width=\textwidth]{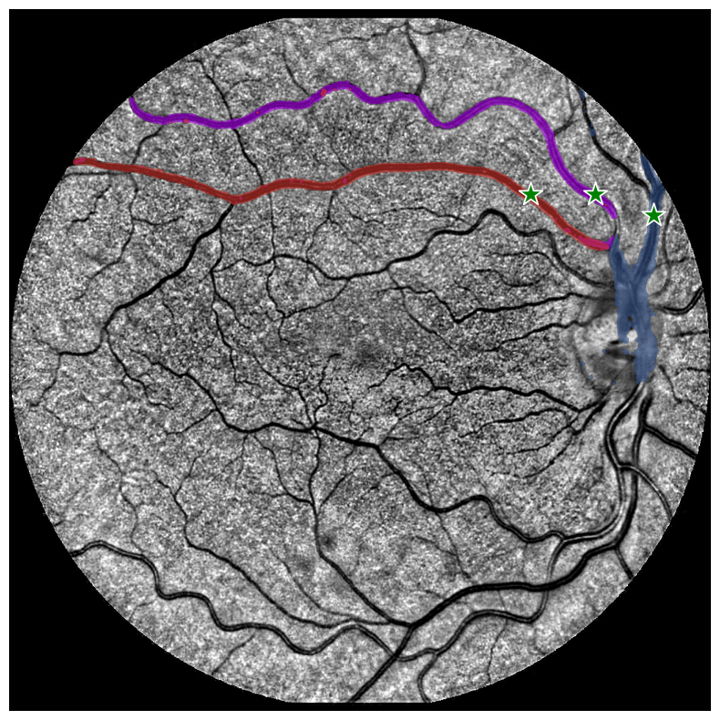}
        \caption{Third best segmentation.}
        \label{fig:Appendix:STAR34:enhanced:segmentation3}
    \end{subfigure}
        
    \centering
    \begin{subfigure}[t]{0.3\linewidth}
        \includegraphics[width=\textwidth]{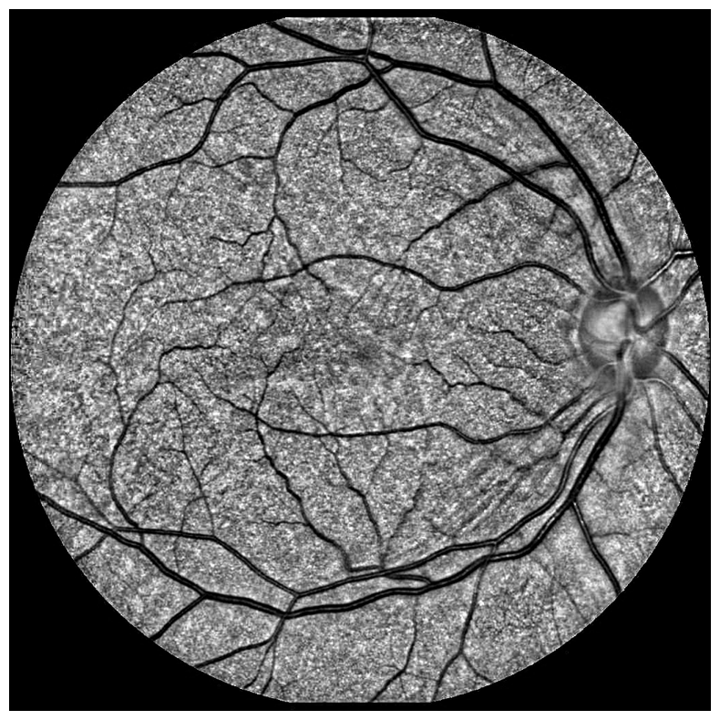}
        \caption{Enhanced retinal image.}
        \label{fig:Appendix:STAR38}
    \end{subfigure}\hfill
    \begin{subfigure}[t]{0.3\linewidth}
        \includegraphics[width=\textwidth]{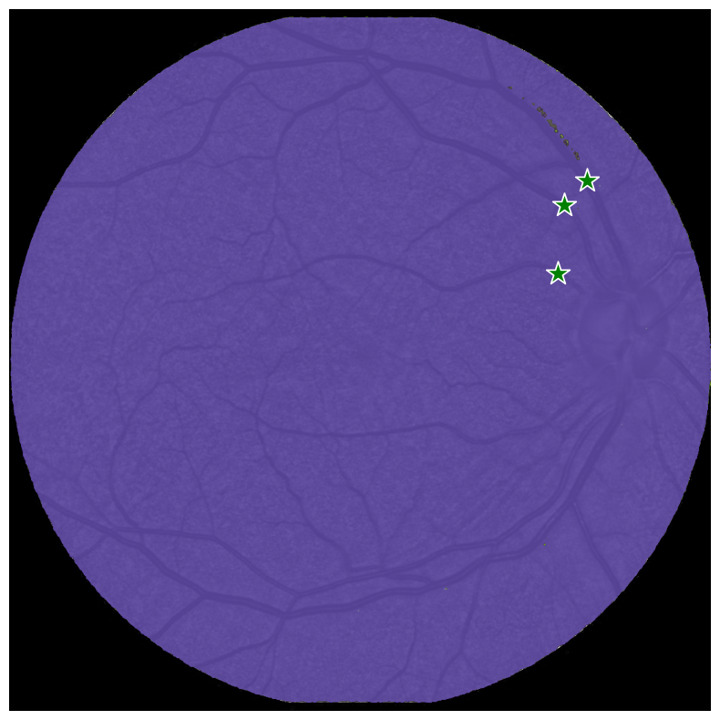}
        \caption{Best segmentation.}
        \label{fig:Appendix:STAR38:enhanced:segmentation1}
    \end{subfigure}\hfill
    \begin{subfigure}[t]{0.3\linewidth}
        \includegraphics[width=\textwidth]{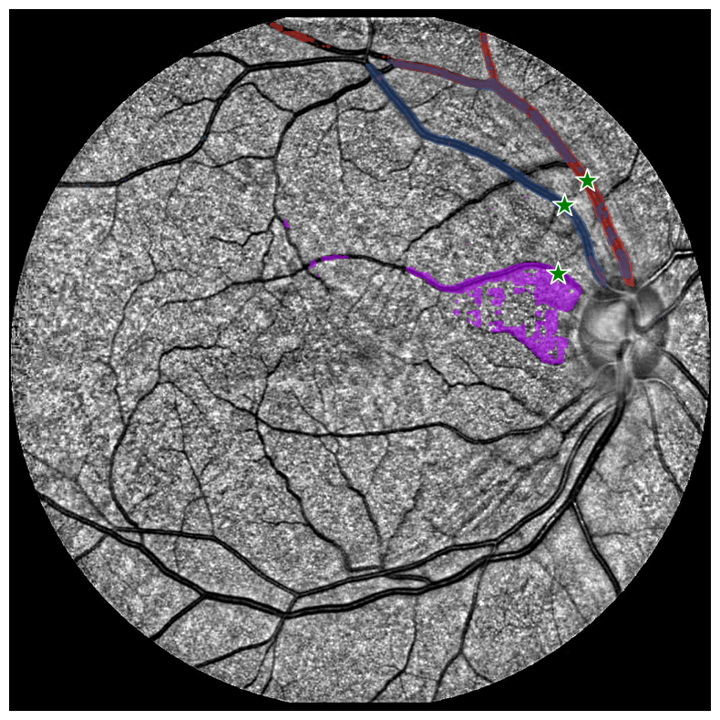}
        \caption{Third best segmentation.}
        \label{fig:Appendix:STAR38:enhanced:segmentation3}
    \end{subfigure}
    
    \caption[The output of the interactive version of the SAM 2 applied to three enhanced retinal images.]{The output of interactive version of the SAM 2 applied to three enhanced retinal images (STAR48, STAR34 and STAR38 in Figs.~\ref{fig:Appendix:STAR48},~\ref{fig:Appendix:STAR34},~\ref{fig:Appendix:STAR38}, respectively, to be compared to the \deltaConnected\ component output in Figs.~\ref{fig:CCIm5},~\ref{fig:CCIm21},~\ref{fig:CCIm24}, respectively). The middle images show the segmentation with the best quality prediction score generated by the model, often generating the same output for all stars. The left images show the segmentation with the third best score. Note that there are situations where the third best segmentation is actually better than the best segmentation according to the model.}
    \label{fig:Appendix:DifferentPredictionScores}
\end{figure*}

\begin{remark}\label{remark:app:prediction}
    The interactive version of SAM 2 produces multiple segmentations per user input. These segmentations are ranked by a quality prediction score generated by the model. We noted cases where the lower-ranked segmentations were closer to the desired output. However, this was not consistently the case. Fig.~\ref{fig:Appendix:DifferentPredictionScores} shows the difference between segmentations of higher and lower rank according to the model applied to different retinal images. Nevertheless, the output of SAM 2 could potentially be improved by adjusting the quality prediction score to our specific use case (i.e., segmentation of blood vessels).
\end{remark}

\begin{table}[]
    \centering
    \begin{tabular}{c|c c c}
         &  \multicolumn{3}{c}{\textbf{Segment Prediction scores}}\\
         Segment:& Top/left & Middle & Bottom/right \\\hline
         Fig.~\ref{fig:Appendix:Circles:segmentation}& 0.789 & 0.590 & 0.746\\
         Fig.~\ref{fig:Appendix:STAR38:segmentation} & 0.988 & 0.988 & 0.992\\
        Fig.~\ref{fig:Appendix:STAR48:enhanced:segmentation1} & 0.186 & 0.145 & 0.143 \\
        Fig.~\ref{fig:Appendix:STAR48:enhanced:segmentation3} & 0.053 & 0.093 & 0.028\\
        Fig.~\ref{fig:Appendix:STAR34:enhanced:segmentation1} & 0.311 & 0.914 & 0.391 \\
        Fig.~\ref{fig:Appendix:STAR34:enhanced:segmentation3} & 0.042 & 0.017 & 0.047 \\
         Fig.~\ref{fig:Appendix:STAR38:enhanced:segmentation1} & 0.193 & 0.938 & 0.945\\
         Fig.~\ref{fig:Appendix:STAR38:enhanced:segmentation3} & 0.079 & 0.032 & 0.030
    \end{tabular}
    \caption[The prediction scores for the segmentations of different images.]{The prediction scores for the segmentations of different images. For more information, see Remark~\ref{remark:app:prediction}}
    \label{tab:SAM:predictionScore}
\end{table}

\subsection{Topological Mode Analysis Tool (ToMATo)}\label{app:ToMATo}
\begin{figure*}
    \begin{subfigure}[t]{0.3\linewidth}
        \includegraphics[width=\textwidth]{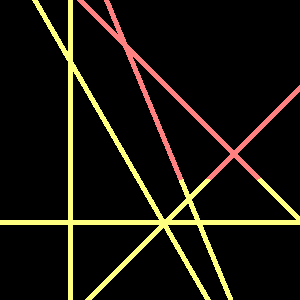}
        \caption{ToMATo applied to the $\bR^2$ image in Fig.~\ref{fig:20230922:comparisonR2vsSE2:Lines}, using \texttt{n\_clusters}=2, leading to 2 components.}
        \label{fig:ToMATo:linesDataR2}
    \end{subfigure}\hfill
    \begin{subfigure}[t]{0.3\linewidth}
        \includegraphics[width=\textwidth]{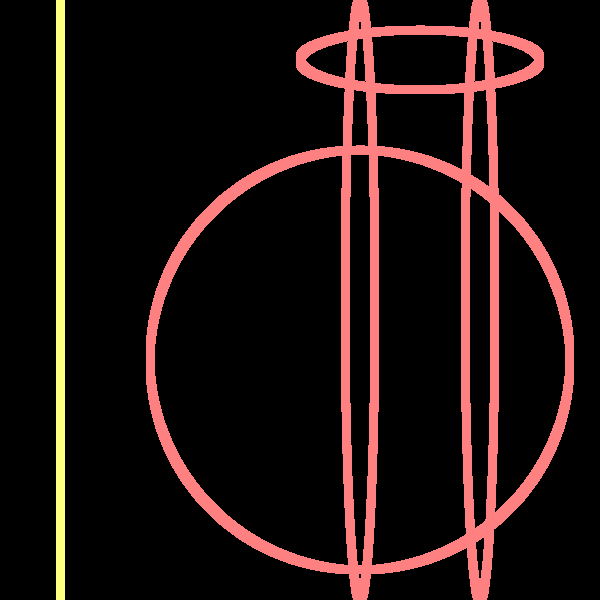}
        \caption{ToMATo applied to the $\bR^2$ image in Fig.~\ref{fig:20230922:comparisonR2vsSE2:chip}, using \texttt{n\_clusters}=2, leading to 2 components.}
        \label{fig:ToMATo:chipDataR2}
    \end{subfigure}\hfill
    \begin{subfigure}[t]{0.3\linewidth}
        \includegraphics[width=\textwidth]{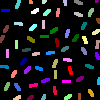}
        \caption{ToMATo applied to the $\bR^2$ image in Fig.~\ref{fig:InfluenceMetricTensorField}, using \texttt{n\_clusters}=3, leading to 57 components.}
        \label{fig:ToMATo:linesR2}
    \end{subfigure}\hfill
    \begin{subfigure}[t]{0.3\linewidth}
        \includegraphics[width=\textwidth]{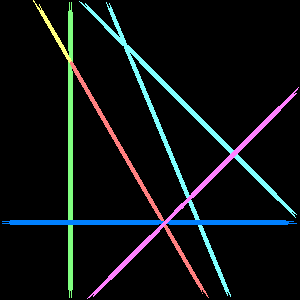}
        \caption{ToMATo applied to the $\SE{2}$ image in Fig.~\ref{fig:20230922:comparisonR2vsSE2:Lines}, using \texttt{n\_clusters}=6, leading to 6 components.}
    \end{subfigure}\hfill
    \begin{subfigure}[t]{0.3\linewidth}
        \includegraphics[width=\textwidth]{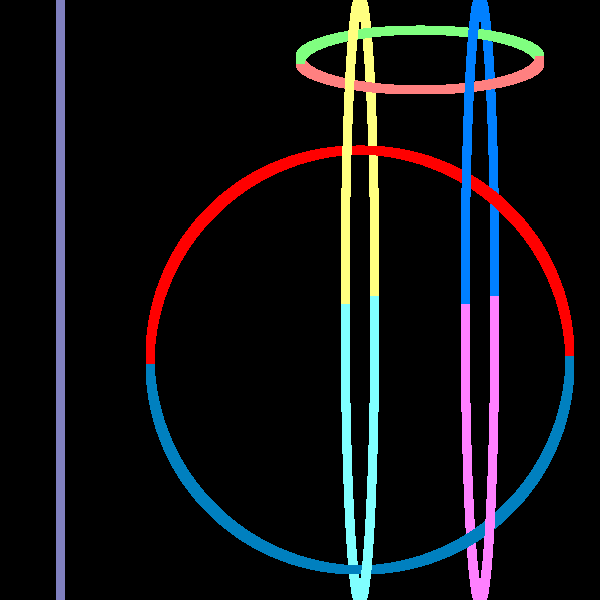}
        \caption{ToMATo applied to the $\SE{2}$ image in Fig.~\ref{fig:20230922:comparisonR2vsSE2:chip}, using \texttt{n\_clusters}=5, leading to 9 components.}
        \label{fig:ToMATo:chipDataSE2}
    \end{subfigure}\hfill
    \begin{subfigure}[t]{0.3\linewidth}
        \includegraphics[width=\textwidth]{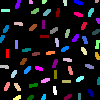}
        \caption{ToMATo applied to the $\SE{2}$ image in Fig.~\ref{fig:InfluenceMetricTensorField}, using \texttt{n\_clusters}=3, leading to 59 components.}
        \label{fig:ToMATo:linesSE2}
    \end{subfigure}\hfill
    \begin{subfigure}[t]{0.45\linewidth}
        \includegraphics[width=\textwidth]{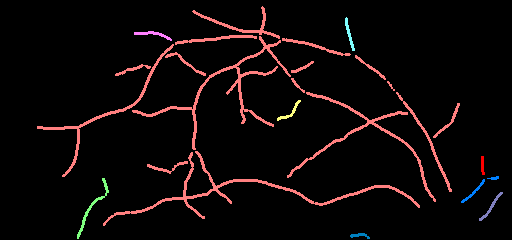}
        \caption{ToMATo applied to the $\bR^2$ image in Fig.~\ref{fig:CCIm24}, using \texttt{n\_clusters}=2, leading to 9 components.}
        \label{fig:ToMATo:vascularTreeR2}
    \end{subfigure}\hfill
    \begin{subfigure}[t]{0.45\linewidth}
        \includegraphics[width=\textwidth]{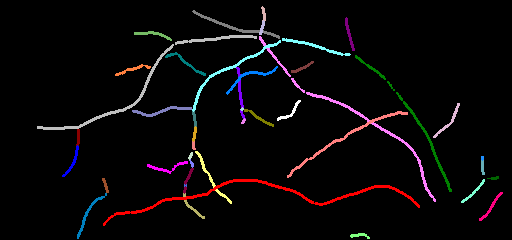}
        \caption{ToMATo applied to the $\SE{2}$ image in Fig.~\ref{fig:CCIm24}, using \texttt{n\_clusters}=5, leading to 42 components.}
        \label{fig:ToMATo:vascularTreeSE2}
    \end{subfigure}
        
    \caption{The output of ToMATo applied to multiple images, optimizing the basic parameter \texttt{n\_clusters} as described in \url{https://gudhi.inria.fr/python/latest/clustering.html}. The algorithm shows difficulties differentiating between crossing structures in $\bR^2$, whereas lifting the image to $\SE2$ does not respect the geometry of the lifted space.}
    \label{fig:ToMATo}
\end{figure*}

ToMATo \cite{Chazal2013persistenceBased} is a useful tool from topological data analysis that allows for the identification of clusters in images. In Fig.~\ref{fig:ToMATo}, we present some results where we applied this tool to some of the images in our introduction and experimental section. We carefully followed the parameter settings and the basic optimization steps described in \url{https://gudhi.inria.fr/python/latest/clustering.html}, always manually choosing the optimal number of clusters in the application of ToMATo. 

The method provides reasonable (sub-)segments, but does struggle to differentiate between different crossing structures in $\bR^2$, as can be seen in Figs.~\ref{fig:ToMATo:linesDataR2}~and~\ref{fig:ToMATo:chipDataR2} compared to Figs.~\ref{fig:20230922:comparisonR2vsSE2:Lines}~and~\ref{fig:20230922:comparisonR2vsSE2:chip}. Lifting the data to $\SE{2}$, and applying the algorithms to this lifted data, does not respect (nor take advantage of) the underlying geometry of the space. 

For instance, the periodicity in the angular direction is not respected, see Fig.~\ref{fig:ToMATo:chipDataSE2}. More importantly, akin to the results by Bekkers et al. \cite[Fig.11]{bekkers2015pde}, we see that (left-invariant) sub-Riemannian geometry (or anisotropic Riemannian geometry \cite[Fig.15, Thm.2]{duits2018optimal}) outperforms isotropic Riemannian geometry, in terms of perceptual organization and tracking of line elements.
Our \deltaConnected\ component algorithm takes advantage of such sub-Riemannian geometry on Lie groups and thereby provides well-aligned elements by design. 

ToMATo does provide a reasonable grouping of elements, but as it is not designed on the Lie group $\SE2$ (and does not include neurogeometrical Lie group models \cite{CittiSanguinetti,Barbieri,DuitsAssociationFieldsViaCusplessSRGeodesicsinSE2,Bellaard2023analysis,petitot2017elements}). Here, it struggles with perceptual grouping of lines as visible in Figs.~\ref{fig:ToMATo:linesR2}~and~\ref{fig:ToMATo:linesSE2}.

Lastly, we applied ToMATo to one of the retinal images in the experimental section. Again, applying the algorithm to the 2-dimensional image data leads to merging of separate vascular trees (cf.~Fig.~\ref{fig:ToMATo:vascularTreeR2}). Lifting the data to $\SE{2}$ results in the splitting of vessel segments (Fig.~\ref{fig:ToMATo:vascularTreeSE2}).

As ToMATo is a mathematically well-underpinned grouping algorithm to identify clusters in data, it could be very interesting to extend the methods to Lie groups such as $\SE2$. This requires a serious redesign of the existing ToMATo algorithm, which is beyond the scope of this paper.

%% file: Review_JMIV/Declarations.tex
\section*{Declarations}
\subsection*{Funding}
We gratefully acknowledge the Dutch Foundation of Science NWO for its financial support by Talent Programme VICI 2020 Exact Sciences (Duits, Geometric learning for Image Analysis, VI.C. 202-031) and the TKI-HTSM Project DELTAS nr. 24PPS053 at TU/e.

The EU is gratefully acknowledged for financial support through the REMODEL - MSCA-SE 101131557 project which supported the research directions and results of this work.

\subsection*{Availability of Data and Materials}
The STAR-dataset \cite{Zhang2016,AbbasiSureshjani} is, together with the \textit{Mathematica} and Python notebooks, available via \cite{Berg2024connectedNotebooks}. 